\documentclass[3pt]{elsarticle}
\pdfoutput=1
\usepackage{amsthm}
\usepackage[intlimits]{amsmath}
\usepackage{amssymb}
\usepackage{amsthm}
\usepackage{amsfonts}
\usepackage{bm}
\usepackage{mathdots}
\usepackage[normalem]{ulem}
\usepackage{setspace} 
\usepackage{lineno} 

\usepackage[utf8]{inputenc} 
\usepackage[T1]{fontenc}
\usepackage{times}
\usepackage{mathptmx}  
\DeclareSymbolFont{cmletters}{OML}{cmr}{m}{n}
\DeclareMathAlphabet{\mathcal}{OMS}{cmsy}{m}{n} 

\journal{Computer and Geotechnics}
\ifdefined\copyReadyOn
\usepackage{ecrc}
\volume{00}
\firstpage{1}
\runauth{G. Beer, Ch. Duenser, V. Mallardo}
\journalname{\journal}
\jid{CG}
\jnltitlelogo{CG}
\CopyrightLine{2015}{Published by Elsevier Ltd.}
\fi

\usepackage{tikz}
\usepackage{tikz-3dplot}
\usepackage{pgfplots}
\usepackage{pgfplotstable}
\usepackage{tikzscale}
\usepackage{standalone} 
\usepackage{hf-tikz}

\usepackage{forloop}
\usepackage{arrayjobx}

\usepackage{algorithmic}
\usepackage{algorithm}

\usepackage{booktabs} 
\usepackage{calc}
\usepackage{xcolor,pdfcomment,soul} 
\usepackage{xspace} 
\usepackage{multirow}

\usepackage{scrlayer-scrpage}
\usepackage{setspace}

\usepackage[labelformat=simple]{subcaption}   	
\usepackage{pdfpages} 

\usepackage{tabto} 

\usepackage{makeidx}
\makeindex

\usepackage{colortbl}

\usepackage{enumerate}

\usepackage{overpic}

\usepackage{listings} 

\usepackage[intlimits]{amsmath}
\usepackage{amssymb} 
\usepackage{amsthm}
\usepackage{amsfonts,mathrsfs}
\usepackage{bm}
\usepackage{mathtools} 

\usepackage{siunitx} 

\usepackage{xfrac}
\usepackage{listings}
\usepackage{color}
\definecolor{dkgreen}{rgb}{0,0.6,0}
\definecolor{gray}{rgb}{0.5,0.5,0.5}
\definecolor{mauve}{rgb}{0.58,0,0.82}

\usepackage[utf8]{inputenc} 
\usepackage[german,english]{babel}  
\usepackage[T1]{fontenc} 
\usepackage{lmodern}

\usepackage{geometry}
\usepackage{overpic}
\graphicspath{{pics/}{tmp/}}

\newboolean{inGray}
\setboolean{inGray}{false}

\newcommand{\addtoindex}[2][]{
    \ifthenelse { \equal{#1}{} }
    {#2\index{#2}\xspace}%
    {#2\index{#1}\xspace}%
}

\newcommand{\myVec}[1]{\mathbf{#1}}
\newcommand{\myVecGreek}[1]{\boldsymbol{#1}}

\newcommand\tens[2]{\mathsf{#1}} %
\newcommand{\trans}{\text{T}} 

\newcommand{\NURBS}{R} 			
\newcommand\uu{\xi} 			
\newcommand\vv{\eta} 			
\newcommand{\domain}{\Omega}

\newcommand{\boundary}{\Gamma}

\newcommand\primary{u} 			
\newcommand\dual{t} 			

\providecommand\url[1]{\emph{#1}}

\newcommand\fund[1]{\tens{#1}{2}}

\newcommand\pt[1]{\boldsymbol{#1}}
\newcommand\sourcept{\tilde{\pt{x}}}
\newcommand\fieldpt{\hat{\pt{x}}}

\newenvironment{mytable}[4]%
{
  \begin{table}[#1]
    \def\mycap{#2}
    \def\mylabel{#3}
    \centering
    \begin{tabular}{#4}
      \toprule
}
{
  \bottomrule
  \end{tabular}
  \caption{\mycap}
  \label{\mylabel}
  \end{table}
}
\newcommand{\mytableheader}[1]{
       #1 \\ \midrule
}

{
  \begin{algorithm}
    \caption{#1}
    \label{#2}
    \begin{algorithmic}[1]
}   
{
\end{algorithmic}
\end{algorithm}
}

\newtheoremstyle{myremark}
{3pt}
{3pt}
{}
{}
{\itshape}
{:}
{.5em}
{}
\theoremstyle{myremark}
\newtheorem*{remark}{Remark}

\newcommand{\myeqref}[1]{equation~(\ref{#1})}	
\newcommand{\myfigref}[1]{Figure~\ref{#1}}

\newcounter{footnoteNumber} 

\DeclareCaptionFormat{listing}{#1#2#3}
\newcommand{\myalignatsinglelabel}[3]%
{%
    \begin{equation}
        \label{#2}
        \begin{alignedat}{#1} 
            #3
        \end{alignedat}
    \end{equation}
}

\newcommand{\myalignat}[2]%
{%
  \begin{alignat}{#1} 
    #2
   \end{alignat}
}

\usepackage{pgf,pgffor,pgfmath}
\pgfplotsset{compat=1.10}

\usetikzlibrary{external}
\usetikzlibrary{positioning,calc,intersections,arrows,3d}
\usetikzlibrary{shapes,snakes,fit,spy}
\usepgfplotslibrary{patchplots}
\usetikzlibrary{matrix}

\usepgfplotslibrary{groupplots}
\usepgfplotslibrary{fillbetween}

\def\pgfplotfontsizetitle{\small}
\def\pgfplotfontsizelegend{\small}
\def\pgfplotfontsize{\small}
\def\pgfplotfontsizetiny{\scriptsize}

\def\tikzfontsizetiny{\scriptsize}

\pgfplotsset{
  mystyle/.style ={%
    grid = major,
    every tick label/.append style={font=\pgfplotfontsizetiny},
    every axis label/.append style={font=\pgfplotfontsize},
    legend style={font=\pgfplotfontsizelegend},
    label style={font=\pgfplotfontsize},
    title style={font=\pgfplotfontsizetitle},
    /pgf/number format/set thousands separator = {}, 
  }
}%

\makeatletter
\pgfplotsset{
    myIgnoreRowModulo2/.style args={#1}{%
        /pgfplots/x filter/.code={%
        \let\xValue\pgfmathresult 
        \pgfmathparse{int(mod(int(\coordindex),int(2))} \pgfmathresult 
        \ifnum#1=\pgfmathresult
            \def\pgfmathresult{} 
        \else
            \pgfmathparse{\xValue} \pgfmathresult 
        \fi
        }
    } 
}
\makeatother

\colorlet{drawblue}      {blue!80!white}
\colorlet{drawred}       {red!80!white}
\colorlet{drawgray}      {gray}
\definecolor{drawgreen}  {RGB}{44,162,95}
\colorlet{drawpurple}    {purple}
\colorlet{draworange}    {orange}
\colorlet{drawlime}      {lime!80!black}
\colorlet{drawartichoke} {yellow!60!black}
\colorlet{TUGgray}{black!15}
\definecolor{TUGred}{RGB}{247,1,70}
\definecolor{IFBblue}{RGB}{51,112,169}

\definecolor{basisColor1}{RGB}{59,76,192}
\definecolor{basisColor2}{RGB}{87,117,225}
\definecolor{basisColor3}{RGB}{119,154,247}
\definecolor{basisColor4}{RGB}{152,185,255}
\definecolor{basisColor5}{RGB}{184,208,249}
\definecolor{basisColor6}{RGB}{195,209,230}	
\definecolor{basisColor7}{RGB}{247,200,190}	
\definecolor{basisColor8}{RGB}{247,187,160}
\definecolor{basisColor9}{RGB}{244,154,123}
\definecolor{basisColor10}{RGB}{229,112,88}
\definecolor{basisColor11}{RGB}{203,62,56}
\definecolor{basisColor12}{RGB}{180,4,38}

\definecolor{basisColor8sw}{RGB} {189,189,189}
\definecolor{basisColor9sw}{RGB} {150,150,150}
\definecolor{basisColor10sw}{RGB}{115,115,115}
\definecolor{basisColor11sw}{RGB}{91,91,91}
\definecolor{basisColor12sw}{RGB}{37,37,37}

\colorlet{myblue}    {blue}
\colorlet{myred}     {red}
\colorlet{mygreen}   {drawgreen}
\colorlet{mypurple}  {purple}
\colorlet{myorange}  {orange}

\ifthenelse{\boolean{inGray}}{ 

\pgfplotscreateplotcyclelist{mycyclelist}{
  basisColor12sw,semithick, every mark/.append style={solid, fill=.},mark=*           \\%
  basisColor11sw,semithick, every mark/.append style={solid, fill=.},mark=square*     \\%
  basisColor10sw,semithick, every mark/.append style={solid, fill=.},mark=triangle*   \\%
  basisColor9sw , semithick,every mark/.append style={solid, fill=.},mark=diamond*    \\%
  basisColor8sw , semithick,every mark/.append style={solid, fill=.},mark=otimes    \\%
}
}{

\pgfplotscreateplotcyclelist{mycyclelist}{
  basisColor12,semithick,   every mark/.append style={solid, fill=.},mark=*           \\%
  basisColor11,semithick,   every mark/.append style={solid, fill=.},mark=square*     \\%
  basisColor10,semithick,   every mark/.append style={solid, fill=.},mark=triangle*   \\%
  basisColor9, semithick,   every mark/.append style={solid, fill=.},mark=diamond*    \\%
  basisColor8, semithick,   every mark/.append style={solid, fill=.},mark=otimes    \\%
}
}

\ifthenelse{\boolean{inGray}}{ 
\pgfplotscreateplotcyclelist{mycyclelistcompare}{
  basisColor12sw, semithick, loosely dashdotdotted     \\%
  basisColor11sw, semithick, loosely dashed      		\\%
  basisColor10sw, semithick, dashed \\
  basisColor9sw,  semithick, solid               	\\%
  basisColor12sw, only marks,   every mark/.append style={solid, fill=.},mark=otimes    \\%
  basisColor11sw, only marks,   every mark/.append style={solid, fill=.},mark=triangle*   \\%
  basisColor10sw, only marks,   every mark/.append style={solid, fill=.},mark=square*     \\%
  basisColor9sw,  only marks,   every mark/.append style={solid, fill=.},mark=*          \\%
}

}{
\pgfplotscreateplotcyclelist{mycyclelistcompare}{
  basisColor12, semithick, loosely dashdotdotted     \\%
  basisColor11, semithick, loosely dashed      		\\%
  basisColor10, semithick, dashed \\
  basisColor9,  semithick, solid               	\\%
  basisColor12, only marks,   every mark/.append style={solid, fill=.},mark=otimes    \\%
  basisColor11, only marks,   every mark/.append style={solid, fill=.},mark=triangle*   \\%
  basisColor10, only marks,   every mark/.append style={solid, fill=.},mark=square*     \\%
  basisColor9,  only marks,   every mark/.append style={solid, fill=.},mark=*          \\%
}
}

\ifthenelse{\boolean{inGray}}{ 
\tikzset{mycyclelistcompareReferenceA/.style={basisColor12sw,solid}}
\tikzset{mycyclelistcompareTestA/.style={basisColor12sw,only marks,mark=otimes}}
}{
\tikzset{mycyclelistcompareReferenceA/.style={basisColor12,solid}}
\tikzset{mycyclelistcompareTestA/.style={basisColor12,only marks,mark=otimes}}
}

\tikzset{helpline/.style={thin,dashed}}
\tikzset{labelline/.style={thin}}
\tikzset{referencePath/.style={dotted,very thick}}
\tikzset{helparrow/.style={thin,arrows={-latex}}}
\tikzset{axis/.style={thin,arrows={->}}}
\tikzset{force/.style={thick,arrows={->}}}
\tikzset{forceInverse/.style={thick,arrows={<-}}}
\tikzset{Gamma/.style={ultra thick}}
\tikzset{controlPoly/.style={draw=black}}
\tikzset{GammaFill/.style={fill=lightgray,fill opacity=0.5}}
\tikzset{colorDiri/.style={drawgreen}}
\tikzset{GammaFillDiri/.style={fill=drawgreen,fill opacity=0.5}}
\tikzset{initialgrid/.style={thin,gray}}
\tikzset{addgridline/.style={dashed,gray}}
\tikzset{trimmingcurve/.style={thick}}
\tikzset{boundingbox/.style={thick, dotted}}
\tikzset{parameterSpace/.style={ }}
\tikzset{basisfunction/.style={very thick,smooth}}
\tikzset{bspline/.style={very thick,smooth,red}}
\tikzset{intersectioncurve/.style={dashed,thick}}
\tikzset{integrationRegionEdge/.style={dashed}}
\tikzset{pointer/.style={arrows={-latex}}}

\tikzstyle{anode}= [circle, inner sep=1.3pt, draw, fill=black]
\tikzstyle{gausspoint}=[shape=circle,draw=black,fill=black,inner sep=1.1pt]
\tikzstyle{controlPoint}=[shape=circle,draw=black,fill=black,thin,inner sep=0pt,minimum size=1.5mm]
\tikzstyle{abscissaPoint}=[shape=circle,draw=black,fill=white,thin,inner sep=0pt,minimum size=1.5mm]
\tikzstyle{anchorPoint}=[shape=circle,draw=black,fill=black,thin,inner sep=0pt,minimum size=1.5mm]
\tikzstyle{anchorPointDeg}=[shape=circle,draw=black,fill=TUGred,thin,inner sep=0pt,minimum size=1.5mm]
\tikzstyle{anchorPointDegD}=[shape=cross out,thick,draw=black,inner sep=0pt,minimum size=1.5mm]
\tikzstyle{trimmingIntersectionPoint}=[shape=cross out,thick,draw=black,inner sep=0pt,minimum size=1.5mm]

\tikzset{%
  highlight/.style={rectangle,rounded corners,fill=red!60,draw,fill opacity=0.125,thick,inner sep=0pt}
}

\usepackage{colortbl}

\def\trianglecolor{black}
\newcommand{\upperSlopeTriangle}[4] 	
{
	\def\trianglecolor{black}
	\addplot[forget plot, domain=#3:#4,color=\trianglecolor,samples=2]{  #2 / (x^#1) } node (A1) [pos=1] {}; 
	\addplot[forget plot, domain=#3:#4,color=\trianglecolor,samples=2]{   #2  / (#3^#1)} node (A2) [pos=1] {} node [anchor=south,pos=0.5,black] {\tikzfontsizetiny $1$};
	\draw[color=\trianglecolor] (A1.center) -- (A2.center) node [anchor=west,pos=0.5,black] {\tikzfontsizetiny #1};
}

\newcommand{\lowerSlopeTriangle}[4] 	
{
	\def\trianglecolor{black}
	\addplot[forget plot, domain=#3:#4,color=\trianglecolor,samples=2]{  #2 / (x^#1) } node (A1) [pos=0] {}; 
	\addplot[forget plot, domain=#3:#4,color=\trianglecolor,samples=2]{   #2  / (#4^#1)} node (A2) [pos=0] {} node [anchor=north,pos=0.5,black] {\tikzfontsizetiny $1$};
	\draw[color=\trianglecolor] (A1.center) -- (A2.center) node [anchor=east,pos=0.5,black] {\tikzfontsizetiny #1};
}

\newcommand{\myaddgraphic}[5]
{
 \node[anchor=south west,inner sep=0] (image) {\phantom{\includegraphics[#2]{#1}}};
  \begin{scope}[x={(image.south east)},y={(image.north west)}]
      
      \begin{scope}
          
          #5
          
          \node[anchor=south west,inner sep=0] {\includegraphics[#2]{#1}};
      \end{scope} 
      
      #4
      
      \pgfmathparse{int(#3)} \let\gridIndicator\pgfmathresult
      \ifthenelse{ \gridIndicator = 1 }
      {
          \draw[help lines,xstep=.1,ystep=.1] (0,0) grid (1.001,1.001);
          \foreach \x in {1,...,9} { \node [anchor=north] at (\x/10,0) {\x};}
          \foreach \y in {1,...,9} { \node [anchor=east] at (0,\y/10) {\y};}
      }{}
      
  \end{scope}    
}

\tikzstyle{reverseclip}=[insert path={(current page.north east) --
    (current page.south east) --
    (current page.south west) --
    (current page.north west) --
    (current page.north east)}
]

\newcounter{itR}
\newcommand{\bsplinevalue}[5] 
{                    				
	\newarray\vKnots
	\newarray\vN
	\newarray\vNumeratorL
	\newarray\vNumeratorR
	\newarray\vSave

	\readarray{vKnots}{#1}
	\readarray{vN}{1}
	\readarray{vSave}{0}
	\readarray{vNumeratorL}{0}
	\readarray{vNumeratorR}{0}

	\foreach \j in {1,...,#2}
	{        
		\pgfmathparse{ int(#4+\j+1) } \checkvKnots(\pgfmathresult)
		\pgfmathsetmacro{\numR}{\cachedata-#3}
		
		\pgfmathparse{ int(#4-\j+2) } \checkvKnots(\pgfmathresult)
		\pgfmathsetmacro{\numL}{#3-\cachedata} 					

		\expandarrayelementtrue
		\pgfmathparse{ int(\j+1) }
		\vNumeratorL(\pgfmathresult)={\numL}
		\vNumeratorR(\pgfmathresult)={\numR}
		
		\forloop[1]{itR}{0}{\value{itR} < \j }
		{
			\pgfmathparse{ int(\theitR+1) } \let\tS\pgfmathresult  	
			\checkvSave(\tS)							
			\pgfmathsetmacro{\save}{\cachedata}  
		
			\pgfmathparse{int(\j-\theitR+1)} \checkvNumeratorL(\pgfmathresult)
			\pgfmathsetmacro{\tmpL}{\cachedata} 	
		    
			\pgfmathparse{ int(\theitR+2) } \checkvNumeratorR(\pgfmathresult)
			\pgfmathsetmacro{\tmpR}{\cachedata} 	       
	
			\pgfmathparse{ int(\theitR+1) } \checkvN(\pgfmathresult)
			\pgfmathparse{ \cachedata / (\tmpL+\tmpR) } \let\tmp\pgfmathresult
			
			\pgfmathparse{ \save + \tmpR * \tmp } 
			\vN(\tS)={\pgfmathresult}

			\pgfmathparse{ \tmpL * \tmp } \let\tmpsave\pgfmathresult
			\pgfmathparse{ int(\theitR+2) } \let\tS\pgfmathresult      
			\vSave(\tS)={\tmpsave}
		}
		
		\pgfmathparse{ int(\j+1) } \checkvSave(\pgfmathresult )
		\vN(\tS)={\cachedata}
	}

	\pgfmathparse{int( #2+1) } \let\lastIndex\pgfmathresult
	\checkvN(1) \pgfmathsetmacro{\first}{\cachedata}  
	\foreach \i [remember=\a as \lasta (initially \first)] in {2,...,\lastIndex}
	{
		\checkvN(\i) \def\a{\lasta,\cachedata} 
		\ifthenelse{\i=\lastIndex}{ \xdef#5{\a} }{}
	}
    
	\foreach \i in {1,...,\lastIndex}
	{
	    \clrarray{vNumeratorR}(\i)
	    \clrarray{vNumeratorL}(\i)
	    \clrarray{vN}(\i)
	    \clrarray{vSave}(\i)
	}
	\delarray\vN
	\delarray\vNumeratorL
	\delarray\vNumeratorR
	\delarray\vSave
}

\newcounter{countvalues}
\newcounter{getBasis}
\newcommand{\bsplinebasis}[4] 
{						
    \newarray\vKnots 	
    \readarray{vKnots}{#1}

    \pgfmathparse{#3}  \let\i\pgfmathresult
    \pgfmathparse{#2}  \let\p\pgfmathresult

    \setcounter{countvalues}{0}
    \setcounter{getBasis}{\p} 
    \pgfmathparse{int(\i+\p)}
    \foreach \knotspan in {\i,...,\pgfmathresult}
    {
	\pgfmathparse{ int(\knotspan+1+1) } \checkvKnots(\pgfmathresult)
	\pgfmathsetmacro{\tmpR}{\cachedata} 						
	
	\pgfmathparse{ int(\knotspan+1) } \checkvKnots(\pgfmathresult) 
	\pgfmathsetmacro{\tmpL}{\cachedata} 						

	\pgfmathparse{ \tmpR - \tmpL } \let\spansize\pgfmathresult  			
   
        \pgfmathparse{ \spansize > 0.0 } \let\bNonZero\pgfmathresult
        \ifthenelse{ \bNonZero = 1 }
        {
            \foreach \percentU in {0,10,...,100}
            {
		\pgfmathparse{\tmpL+\spansize*\percentU/100} \let\u\pgfmathresult

		\bsplinevalue{#1}{\p}{\u}{\knotspan}{\Basis} 
		\def\basisfuncarray{{\Basis}} 				
		
		\pgfmathparse{\basisfuncarray[\thegetBasis]} 
		
		\ifthenelse{\thecountvalues=0}
		{ 
			\xdef\nodeB{"\u,\pgfmathresult"}
		}{
			\xdef\nodeB{\nodeB,"\u,\pgfmathresult"}  
		}
		\addtocounter{countvalues}{1}
            }
        }{}
        \addtocounter{getBasis}{-1}
    }
    
	\xdef#4{\nodeB}

	\delarray\vKnots
}

\newcommand{\bsplinebasisspan}[6] 
{							
    \newarray\vKnots 	
    \readarray{vKnots}{#1}

    \pgfmathparse{#5}  \let\plotknotspan\pgfmathresult
    \pgfmathparse{#4}  \let\splineknotspan\pgfmathresult
    \pgfmathparse{#3}  \let\i\pgfmathresult
    \pgfmathparse{#2}  \let\p\pgfmathresult

    \setcounter{countvalues}{0}
    \setcounter{getBasis}{\i} 
    \foreach \knotspan in {\plotknotspan}
    {
	\pgfmathparse{ int(\knotspan+1+1) } \checkvKnots(\pgfmathresult)
	\pgfmathsetmacro{\tmpR}{\cachedata} 						
	
	\pgfmathparse{ int(\knotspan+1) } \checkvKnots(\pgfmathresult) 
	\pgfmathsetmacro{\tmpL}{\cachedata} 						

	\pgfmathparse{ \tmpR - \tmpL } \let\spansize\pgfmathresult  			
   
        \pgfmathparse{ \spansize > 0.0 } \let\bNonZero\pgfmathresult
        \ifthenelse{ \bNonZero = 1 }
        {
            \foreach \percentU in {0,10,...,100}
            {
		\pgfmathparse{\tmpL+\spansize*\percentU/100} \let\u\pgfmathresult

		\bsplinevalue{#1}{\p}{\u}{\splineknotspan}{\Basis} 
		\def\basisfuncarray{{\Basis}} 				
		
		\pgfmathparse{\basisfuncarray[\thegetBasis]} 
		
		\ifthenelse{\thecountvalues=0}
		{ 
			\xdef\nodeB{"\u,\pgfmathresult"}
		}{
			\xdef\nodeB{\nodeB,"\u,\pgfmathresult"}  
		}
		\addtocounter{countvalues}{1}
            }
        }{}
        \addtocounter{getBasis}{-1}
    }
    
	\xdef#6{\nodeB}

	\delarray\vKnots
}

\newcommand{\plotbsplinebasis}[4] 	
{							
	\bsplinebasis{#1}{#2}{#3}{\nodeOut}
	\def\nodearray{{\nodeOut}}

	\xdef\name{ }
	\addtocounter{countvalues}{-1}
	\foreach \i in {0,...,\thecountvalues}
	{
		\pgfmathparse{\nodearray[\i]}
		\coordinate (point\i) at (\pgfmathresult);	  
		\xdef\name{ \name (point\i)  }
	}
	
	\draw[#4] plot coordinates{ \name };
	
	\xdef\name{ }
}

\newcommand{\plotbsplinesegment}[6] 	
{								
	\bsplinebasisspan{#1}{#2}{#3}{#4}{#5}{\nodeOut}
	\def\nodearray{{\nodeOut}}

	\xdef\name{ }
	\addtocounter{countvalues}{-1}
	\foreach \i in {0,...,\thecountvalues}
	{
		\pgfmathparse{\nodearray[\i]}
		\coordinate (point\i) at (\pgfmathresult);	  
		\xdef\name{ \name (point\i)  }
	}
	
	\draw[#6] plot coordinates{ \name };
	
	\xdef\name{ }
}

\newcommand{\plotbsplineaccumulated}[5] 	
{						
    \newarray\vKnots 	
    \readarray{vKnots}{#1}
    \newarray\vSubCoef 	
    \readarray{vSubCoef}{#3}
    
    \pgfmathparse{#4}  \let\plotknotspan\pgfmathresult
    \pgfmathparse{#4}  \let\splineknotspan\pgfmathresult
    \pgfmathparse{#2}  \let\p\pgfmathresult
    
    \setcounter{countvalues}{0}
    \pgfmathparse{int( \p+1) } \let\lastIndex\pgfmathresult
    \foreach \knotspan in {\plotknotspan}
    {
        \pgfmathparse{ int(\knotspan+1+1) } \checkvKnots(\pgfmathresult)
        \pgfmathsetmacro{\tmpR}{\cachedata} 						
        
        \pgfmathparse{ int(\knotspan+1) } \checkvKnots(\pgfmathresult) 
        \pgfmathsetmacro{\tmpL}{\cachedata} 						
        
        \pgfmathparse{ \tmpR - \tmpL } \let\spansize\pgfmathresult  			
        
        \pgfmathparse{ \spansize > 0.0 } \let\bNonZero\pgfmathresult
        \ifthenelse{ \bNonZero = 1 }
        {
            \foreach \percentU in {0,10,...,100}
            {
                \pgfmathparse{\tmpL+\spansize*\percentU/100} \let\u\pgfmathresult
                
                \bsplinevalue{#1}{\p}{\u}{\splineknotspan}{\Basis} 
                \def\basisfuncarray{{\Basis}} 				
                
                \setcounter{getBasis}{0} 
                \pgfmathparse{\basisfuncarray[\thegetBasis]} 
                \let\basisValue\pgfmathresult
                
                \checkvSubCoef(1) \pgfmathsetmacro{\coef}{\cachedata}  
                \pgfmathparse{ \basisValue * \coef } \let\first\pgfmathresult
                
                \xdef\lastx{\first}
                \foreach \i in {2,...,\lastIndex}
                { 
                    \addtocounter{getBasis}{1}       
                    \pgfmathparse{\basisfuncarray[\thegetBasis]}
                    \let\basisValue\pgfmathresult
                    
                    \checkvSubCoef(\i) \pgfmathsetmacro{\coef}{\cachedata}  
                    \pgfmathparse{ \lastx + \basisValue * \coef } \let\sum\pgfmathresult
                    
                    \xdef\lastx{\sum}
                    
                    \ifthenelse{\i=\lastIndex}
                    {
                        \ifthenelse{\thecountvalues=0}
                        { 
                            \xdef\nodeBB{"\u,\sum"}
                        }{
                            \xdef\nodeBB{\nodeBB,"\u,\sum"}  
                        }
                        \addtocounter{countvalues}{1}
                    }{}
                    
                }
            }
        }{}
    }
    
    \delarray\vKnots
    \delarray\vSubCoef
    
    \def\nodearray{{\nodeBB}}
    
    \xdef\name{ }
    \addtocounter{countvalues}{-1}
    \foreach \i in {0,...,\thecountvalues}
    {
        \pgfmathparse{\nodearray[\i]}
        \coordinate (point\i) at (\pgfmathresult);                
        \xdef\name{ \name (point\i)  }
    }
    
    \draw[#5] plot coordinates{ \name };
    
    \xdef\name{ }
}

\begin{document}
    
\title{Efficient and realistic 3-D Boundary Element simulations of underground construction using isogeometric analysis.}
\begin{frontmatter}

\author[ifbaddr]{Gernot Beer\corref{cor1}}
\author[ifbaddr]{Christian Duenser}
\author[unife]{Vincenzo Mallardo}

\address[ifbaddr]{Institute of Structural Analysis, Graz University
  of Technology, Lessingstraße 25/II, 8010 Graz, Austria}

\address[unife]{Department of Architecture, University of Ferrara, Via Quartieri 8, 44121 Ferrara, Italy}

\cortext[cor1]{Corresponding author.
  Tel.: +43 316 873 6181, fax: +43 316 873 6185, mail: \url{gernot.beer@tugraz.at}, web: \url{www.ifb.tugraz.at}}

\begin{abstract}
The paper outlines some recent developments of the boundary element method (BEM) that makes it more user friendly and suitable for a realistic simulation in geomechanics, especially for underground excavations and tunnelling.  The innovations refer to the introduction of isogeometric concepts, elasto-plastic analysis and the simulation of ground support. The introduction of isogeometric concepts for the description of the excavation boundaries results in less user and analysis effort, since complex geometries can be modelled with few parameters and degrees of freedom. No mesh generation is necessary.
Heterogeneous and inelastic ground conditions are considered via general inclusions and rock bolts via linear inclusions.

A comparison of results of test examples with other numerical methods and analytical solutions confirm the efficiency and accuracy of the proposed implementation.
A practical example with a complex geometry is presented.

\end{abstract}

\end{frontmatter}

\section{Introduction}
Since the publication of the first paper on the topic \cite{Hughes2005a}, isogeometric analysis has gained increased popularity. 
The majority of applications have been with the Finite Element method  (FEM) and much less with the BEM.
However, the advantage of the BEM, that requires only the discretisation of the boundary, makes it an ideal companion to Computer Aided Design (CAD).
First applications of the isogeometric BEM (IGABEM) were published in elasticity in 2-D \cite{simpson2012two,simpson2013isogeometric} and in 3-D \cite{scott2013isogeometric}. Other applications followed (for example see \cite{An2018},\cite{Fang2020}).
In \cite{Marussig2015} the concept of a geometry independent field approximation, which involved a decoupling of the geometry definition and the approximation of the unknown,  was first introduced and was later adopted by others \cite{doi:10.1002/nme.5778}. The seamless integration of BEM and CAD was discussed in \cite{marussig2016b}.
In a recent book published on the isogeometric BEM \cite{BeerMarussig} it was shown how geometrical information can be taken directly from CAD data and that efficient and accurate simulations with very few unknowns can be obtained. 

One fact that has hampered the widespread use of the BEM is that fundamental solutions, on which the method is based, exist only for elastic material properties and homogeneous domains. Fundamental solutions can be obtained for anisotropic materials, but they are very complicated \cite{PanChow1976}.
To overcome the fact that the original BEM can only deal with homogeneous and elastic domains, several workarounds were introduced. Among them we mention the introduction of boundary element regions, to consider a piecewise heterogeneous domain \cite{BanerjeeRaveendra1986} and the coupling of the BEM with the Finite Element Method (FEM) where the FEM regions model non-linear behaviour \cite{Wendlandbook}. The topic of elastic inclusions was recently dealt with in \cite{Sun2020}. The concept of including non-linear effects by adding a volume integral was first introduced by Brebbia \cite{Brebbia1980} and Banerjee \cite{Banerjee} and later expanded in \cite{Gao2011}. 
 Various ways of avoiding the generation of a volume mesh were presented (see for example\cite{Tanaka1}), but their application is limited to finite domain problems. All of the mentioned solutions involve the introduction of errors or increase the discretisation effort.

In this paper we show several innovations that make the BEM suitable for realistic simulations in underground construction. This means that the BEM simulation can consider ground support, heterogenous ground conditions and non-linear material behaviour. 
We first introduce the theoretical background of the BEM with volume effects. Then the evaluation of the arising boundary and volume integrals is discussed in some detail. This includes the description of the excavation geometry with NURBS patches and the definition of subdomains, where material properties differ from the ones used to compute the fundamental solutions or which behave in an inelastic way. Rock bolts are also modelled as (linear) subdomains and their analytical integration allows many of them to be used in a simulation with a small increase in the numerical effort. 

Two test examples are included, which test the accuracy and efficiency of the simulation of ground support and elasto-plastic material behaviour. Finally a practical example with some complexity is presented.

It is emphasised that the simulation approach presented here does not involve any mesh generation. Instead, geometries are defined by NURBS patches using data generated by CAD programs or input data, in a similar data format, that are user generated. 

\section{Theory}

In the following we will use the word \textbf{inclusion} to specify part of the analysis domain that is not modelled by boundary elements, i.e. parts of the domain that has different material properties or behaves inelastically. This also applies to the ground support.

As will be explained, we use the concept of initial stresses inside inclusions to consider those volume effects, which are not considered by the boundary discretisation.

In the following we will use matrix algebra and it is therefore necessary to convert the stress and strain tensors $\sigma_{ij}, \epsilon_{ij}$ to pseudo-vectors $\myVecGreek{\sigma},  \myVecGreek{\epsilon}$ using \textit{Voigt} notation:

\begin{align}
\label{Voightnot}
\myVecGreek{\sigma}= \left\{\begin{array}{c}\sigma_{11} \\\sigma_{22} \\\sigma_{33} \\\sigma_{12} \\\sigma_{23}\\\sigma_{13}\end{array}\right\}   && \text{and} &&  \myVecGreek{\epsilon}= \left\{\begin{array}{c}\epsilon_{11} \\\epsilon_{22} \\\epsilon_{33} \\\epsilon_{12} + \epsilon_{21}\\\epsilon_{23} + \epsilon_{32}\\\epsilon_{13}+ \epsilon_{31}\end{array}\right\}  
\end{align}

The initial stress vector $ \myVecGreek{\sigma}_{0}$, due the fact that a point inside the inclusion has properties that are different to the ones used for computing the fundamental solutions, is given by:
\begin{equation}
\label{InitialS}
\myVecGreek{\sigma}_{0} =( \mathbf{ D} - \mathbf{ D}_{incl}) \myVecGreek{\epsilon} 
\end{equation}
where $\myVecGreek{\epsilon} $ is the total strain, $ \mathbf{ D}$  is the constitutive matrix for computing the fundamental solutions and $ \mathbf{ D}_{incl}$ is the corresponding matrix for the inclusion point. 

\begin{remark}
While the matrix $ \mathbf{ D}$ is restricted to an isotropic elasticity matrix, the constitutive matrix $ \mathbf{ D}_{incl}$ can be quite general, ranging from the sparsely populated isotropic elasticity matrix to a fully populated anisotropic elasticity matrix. Indeed, later on we will use an elasto-plastic constitutive matrix for $ \mathbf{ D}_{incl}$.
\end{remark}

In the following we first establish the governing integral equations and then discuss in detail how the arising volume and surface integrals are evaluated.

\subsection{Governing integral equations}
Consider a domain $\domain$ with a boundary $\boundary$, containing a subdomain $\domain_{0}$ where initial stresses $ \myVecGreek{\sigma}_{0}(\fieldpt)$ are present.
We apply the theorem by Betti and the \textbf{collocation method} to arrive at the governing integral equations.
This means that we set the work done on the boundary $\boundary$ by tractions $\fund{T}$ times displacements $\myVec{u}$ equal to the work done by displacements $\fund{U}$ times tractions $\myVec{t}$. We assume $\fund{T}$ and $\fund{U}$ to be fundamental solutions of the governing differential equation at $\fieldpt$ due to a source at $\sourcept_{n}$ and $\myVec{u}$, $\myVec{t}$ to be boundary values.
If initial stresses are present, additional work is done in the domain $\domain_{0}$ by the initial stresses $ \myVecGreek{\sigma}_{0}(\fieldpt)$ times the fundamental solution for strains $\fund{E} (\sourcept_{n},\fieldpt)$.
The integral equation can be written as (see \cite{Brebbia1983}, \cite{aliabadi},\cite{Banerjee1981}):
\begin{equation}
\label{Beq}
    \begin{aligned}
\int_{\boundary} \fund{T}(\sourcept_{n},\fieldpt)  \myVec{u}(\fieldpt) ) \ d\boundary(\fieldpt) ) &=& \int_{\boundary} \fund{U}(\sourcept_{n},\fieldpt) \ \myVec{t}(\fieldpt) \ d\boundary(\fieldpt) \\
+  \int_{\domain_{0}} \fund{E} (\sourcept_{n},\fieldpt)  \myVecGreek{\sigma}_{0} (\fieldpt) d \domain_{0} (\fieldpt)  .
\end{aligned}
\end{equation}
where $\sourcept_n$ are the coordinates of the \textbf{collocation point} $n$. To be able to solve the integral equations they have to be regularised.
The regularised integral equations are written as :
\begin{equation}
\label{BEM}
    \begin{aligned}
\int_{\boundary} \fund{T}(\sourcept_{n},\fieldpt) ( \myVec{u}(\fieldpt) - \myVec{u}(\sourcept_{n})) \ d\boundary(\fieldpt) - \mathbf{ A}_{n} \myVec{u}(\sourcept_{n}) &=& \int_{\boundary} \fund{U}(\sourcept_{n},\fieldpt) \ \myVec{t}(\fieldpt) \ d\boundary(\fieldpt) \\
+  \int_{\domain_{0}} \fund{E} (\sourcept_{n},\fieldpt)  \myVecGreek{\sigma}_{0} (\fieldpt) d \domain_{0} (\fieldpt)  .
\end{aligned}
\end{equation}
where $\mathbf{ A}_{n}=\mathbf{ 0}$ for finite domain problems and $\mathbf{ A}_{n}=\mathbf{ I}$ for infinite domain problems. The derivation of Eq. (\ref{BEM}) and the fundamental solutions $\fund{U}$ und $\fund{T}$ are presented in \cite{BeerMarussig}.
The fundamental solution $\fund{E}$ is given by:
\begin{equation}
\label{ }
E_{ijk}= \frac{-C}{r^{2}}\left[C_{3}(r_{,k}\delta_{ij} + r_{,j}\delta_{ik}) - r_{,i}\delta_{jk} + C_{4} \ r_{,i} r_{,j} r_{,k}\right]
\end{equation}
where $r=|\fieldpt - \sourcept|$, $r_{,i}=  \frac{r_{i}}{r}$ and $\delta_{ij}$ is the \textit{Kronecker Delta}. The constants are: $C=\frac{1}{16 \pi G (1- \nu})$, $C_{3}=1 - 2\nu$ and $C_{4}=3$ where $G$ is the shear modulus and $\nu$ the Poisson's ratio.

The tensor $E_{ijk}$ is converted to a matrix $\fund{E}$:
\begin{equation}
\label{eq_Voigt_convertion}
\fund{E}= \left[\begin{array}{cccccc}E_{111} & E_{122} & E_{133}  & E_{112}  + E_{121}& E_{123} + E_{132} & E_{113} +E_{131}  \\E_{211} & E_{222} & E_{233}  & E_{212} + E_{221}  & E_{223}  + E_{232} & E_{213} + E_{231}\\E_{311} & E_{322} & E_{333}  & E_{312} + E_{321} & E_{323}  + E_{332}& E_{313} + E_{331}\end{array}\right]
\end{equation}

\section{Discretisation of integral equations}
To be able to solve Equations (\ref{BEM}) we have to discretise them. This involves 2 steps:
\begin{itemize}
  \item  The subdivision of the boundary domain into patches and the volume domain into inclusions
  \item  The approximation of the unknown boundary values and the approximation of initial stresses.
\end{itemize}
This will be discussed in the subsequent sections separately for the boundary and volume integrals.
\subsection{Discretisation of the boundary integrals}
For the numerical solution of the boundary integral equations the integrals are expressed as sum of integrals over patches:
\begin{equation}
    \begin{aligned}
        \label{Regular2}
 \int_{\boundary} \fund{U}(\sourcept_{n},\fieldpt) \ \myVec{t}(\fieldpt) \ d\boundary(\fieldpt) - \int_{\boundary} \fund{T}(\sourcept_{n},\fieldpt) ( \myVec{u}(\fieldpt) - \myVec{u}(\sourcept_{n})) \ d\boundary(\fieldpt) + \mathbf{ A}_{n} \myVec{u}(\sourcept_{n})  = \\
        \sum_{e=1}^{E} \int_{\Gamma_{e}}  \fund{U}(\sourcept_{n},\fieldpt) \ \myVec{\dual}^{e}(\fieldpt)\ d\Gamma_{e}(\fieldpt)  
        - \sum_{e=1}^{E} \int_{\Gamma_{e}} \fund{T}(\sourcept_{n},\fieldpt)  \myVec{\primary}^{e}(\fieldpt) d \Gamma_{e}  \\
        +  \left[\sum_{e=1}^{E} \ \left(\int_{\Gamma_{e}} \fund{T}(\sourcept_{n},\fieldpt) d \Gamma_{e} \right) + \mathbf{A}_{n}\right] \myVec{\primary}(\sourcept_{n})
    \end{aligned}
\end{equation}
where $e$ specifies the patch number and $E$ is the total number of patches. In the following the geometry of patches is specified using NURBS basis functions. The advantage of this is that some geometrical shapes such as cylinder and spheres can be described exactly with few parameters.
For further information on NURBS and how (\ref{Regular2}) is obtained the reader is referred to  \cite{BeerMarussig}.

There are 3 types of patches that are useful for geomechanics modelling: finite, infinite and trimmed patches. In addition we introduce a special patch with a cutout that can be used for modelling intersections.

\subsubsection{Geometry definition of finite patches}
\begin{figure}
\begin{center}
\begin{overpic}[scale=0.5]{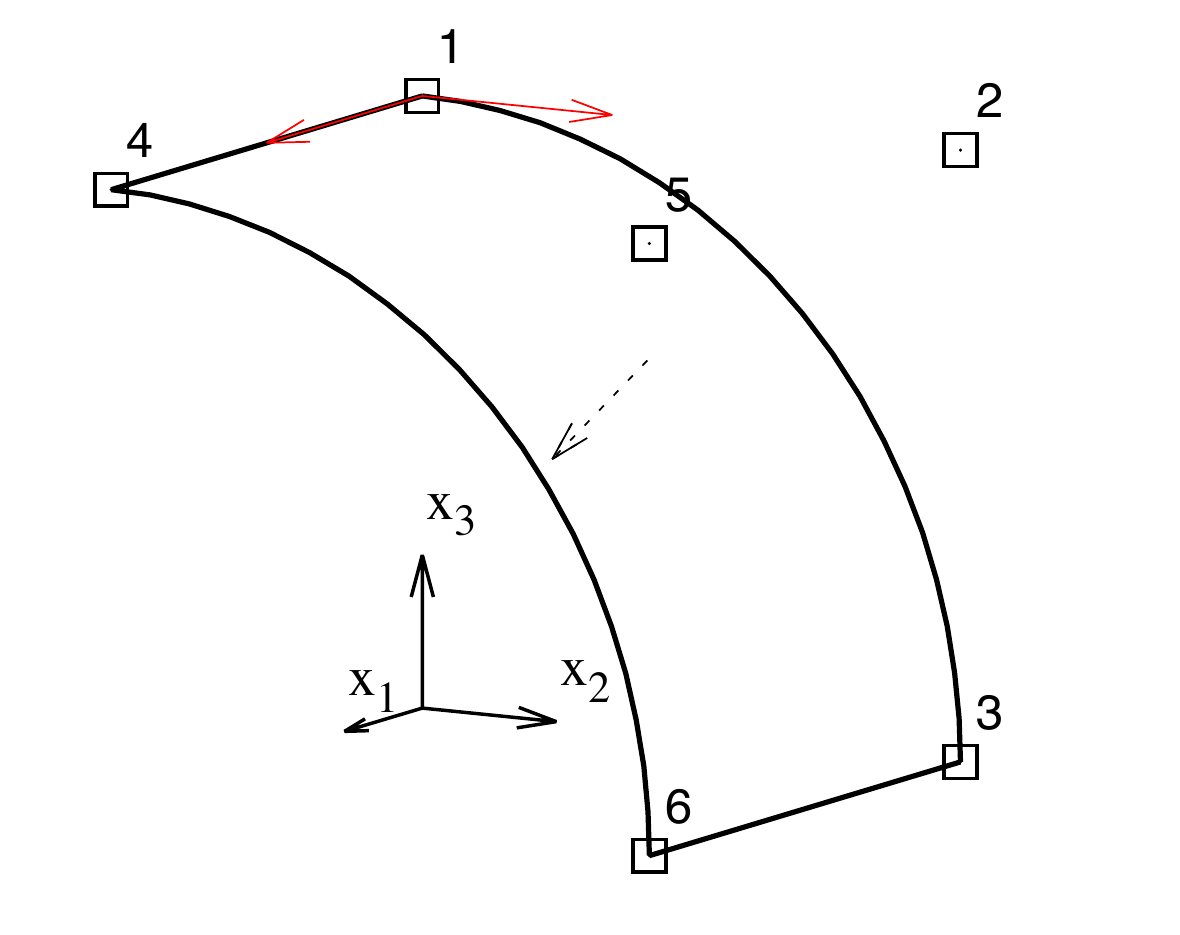}
 \put(50,75){$\uu$}
 \put(20,70){$\vv$}
\end{overpic}
\begin{overpic}[scale=0.5]{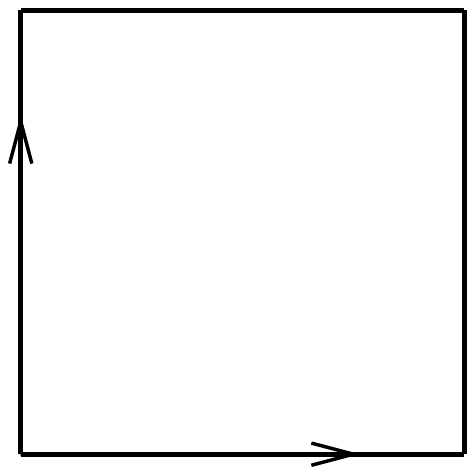}
 \put(70,25){$\uu$}
 \put(30,50){$\vv$}
\end{overpic}
\caption{A finite patch with control points (numbered squares). Left: in the global, right: in the local coordinate system. Also shown is the ``outward normal''.}
\label{fig6:Patch3D}
\end{center}
\end{figure}

In \myfigref{fig6:Patch3D} we show an example of a finite patch.
The mapping from the local $\myVecGreek{\xi}(\uu,\vv)$ to the global $\pt{x}$ coordinate system is given by
\begin{equation}
\pt{x} (\xi,\eta)= \sum_{i=1}^{I} \NURBS_{i} (\xi, \eta) \pt{x}_{i}.
\end{equation}
where $\NURBS_{i} (\xi, \eta)$ are NURBS basis functions and the control points (coordinates $\pt{x}_{i}$) are numbered consecutively, first in the $\uu$- and then in the $\vv$-direction.

The vectors tangential to the surface are given by

\begin{align}
\mathbf{ v}_{\xi}= \frac{\partial \pt{x}}{\partial \xi}= \left(\begin{array}{c}\frac{\partial x_{1}}{\partial \xi} \\ \\ \frac{\partial x_{2}}{\partial \xi} \\ \\ \frac{\partial x_{3}}{\partial \xi}\end{array}\right)
&& \text{and} && \mathbf{ v}_{\eta}= \frac{\partial \pt{x}}{\partial \eta}= \left(\begin{array}{c}\frac{\partial x_{1}}{\partial \eta} \\ \\ \frac{\partial x_{2}}{\partial \eta} \\ \\ \frac{\partial x_{3}}{\partial \eta} \end{array}\right)
\end{align}
and the unit vector normal is
\begin{equation}
\mathbf{ n}= \frac{\mathbf{ v}_{\xi} \times  \mathbf{ v}_{\eta}}{J}.
\end{equation}
The Jacobian is
\begin{equation}
J= | \mathbf{ v}_{\xi} \times  \mathbf{ v}_{\eta} |.
\end{equation}
 The  direction of the ``outward normal'' depends on how the control points are numbered.

\subsubsection{Geometry definition of infinite patches}
\label{sec:InfinitePatches}
\begin{figure}
\begin{center}
\begin{overpic}[scale=0.5]{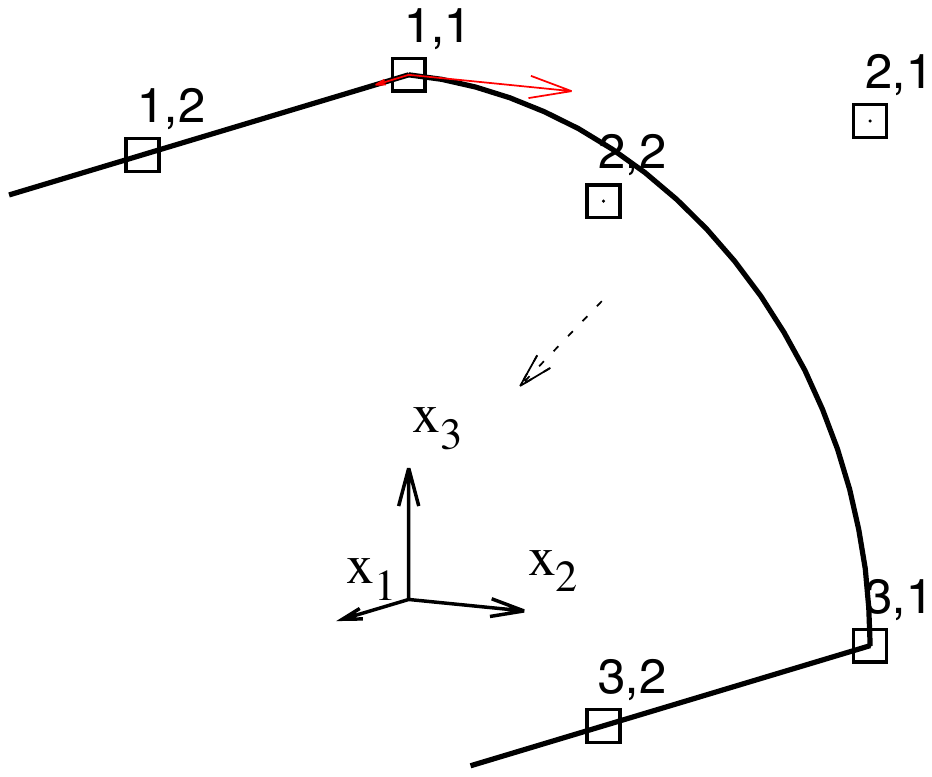}
 \put(50,70){$\uu$}
 \put(30,68){$\vv$}
\end{overpic} 
\begin{overpic}[scale=0.5]{pics/Patch3Dst.pdf}
 \put(70,25){$\uu$}
 \put(30,50){$\vv$}
\end{overpic}
\caption{Example of an infinite patch. Left in the global and right in the local coordinate system}
\label{fig6:Ptch3Dinf}
\end{center}
\end{figure}

Here we introduce a patch definition that is useful for the simulation in geomechanics where one sometimes has to consider a surface that tends to infinity \cite{Beer2015b}. In this case we define an infinite patch as shown in \myfigref{fig6:Ptch3Dinf}.
The mapping for a patch that extends to infinity in the $\eta$-direction is given by
\begin{equation}
\pt{x}= \sum_{j=1}^{2}\sum_{i=1}^{I} \NURBS^{\infty}_{ij}(\xi,\eta) \pt{x}_{ij}
\end{equation}
where

\begin{equation}
\NURBS^{\infty}_{ij}(\xi,\eta)= \NURBS_{i}(\xi) M_{j}^{\infty}(\eta)
\end{equation}

and the special infinite basis functions are

\begin{align}
 M_{1}^{\infty} =  \frac{1 - 2\eta}{1-\eta} && \text{and} &&  M_{2}^{\infty}  =  \frac{\eta}{1-\eta}. 
\end{align}

The vectors in the tangential directions are given by

\begin{eqnarray}
\mathbf{ v}_{\xi}= \frac{\partial \pt{x}}{\partial \xi}=\sum_{j=1}^{2}\sum_{i=1}^{I} \frac{\partial \NURBS_{i}(\xi)}{\partial \xi}  M_{j}^{\infty}(\eta) \pt{x}_{ij} \\
\mathbf{ v}_{\eta}= \frac{\partial \pt{x}}{\partial \eta}=\sum_{j=1}^{2}\sum_{i=1}^{I} \NURBS_{i}(\xi)  \frac{\partial M_{j}^{\infty}(\eta)}{\partial \eta} \pt{x}_{ij} 
\end{eqnarray}

where

\begin{align}
\frac{\partial M_{1}^{\infty}}{\partial \eta}= \frac{-1}{(1-\eta)^{2}} && \text{and} && \frac{\partial M_{2}^{\infty}}{\partial \eta}= \frac{1}{(1-\eta)^{2}}.
\end{align}

The unit vector normal is computed as for the finite patch. It is noted that the Jacobian $J$ tends to infinity as $\vv$ tends to 1.

\subsubsection{Trimmed patches}
\begin{figure}
\begin{center}
\begin{overpic}[scale=0.45]{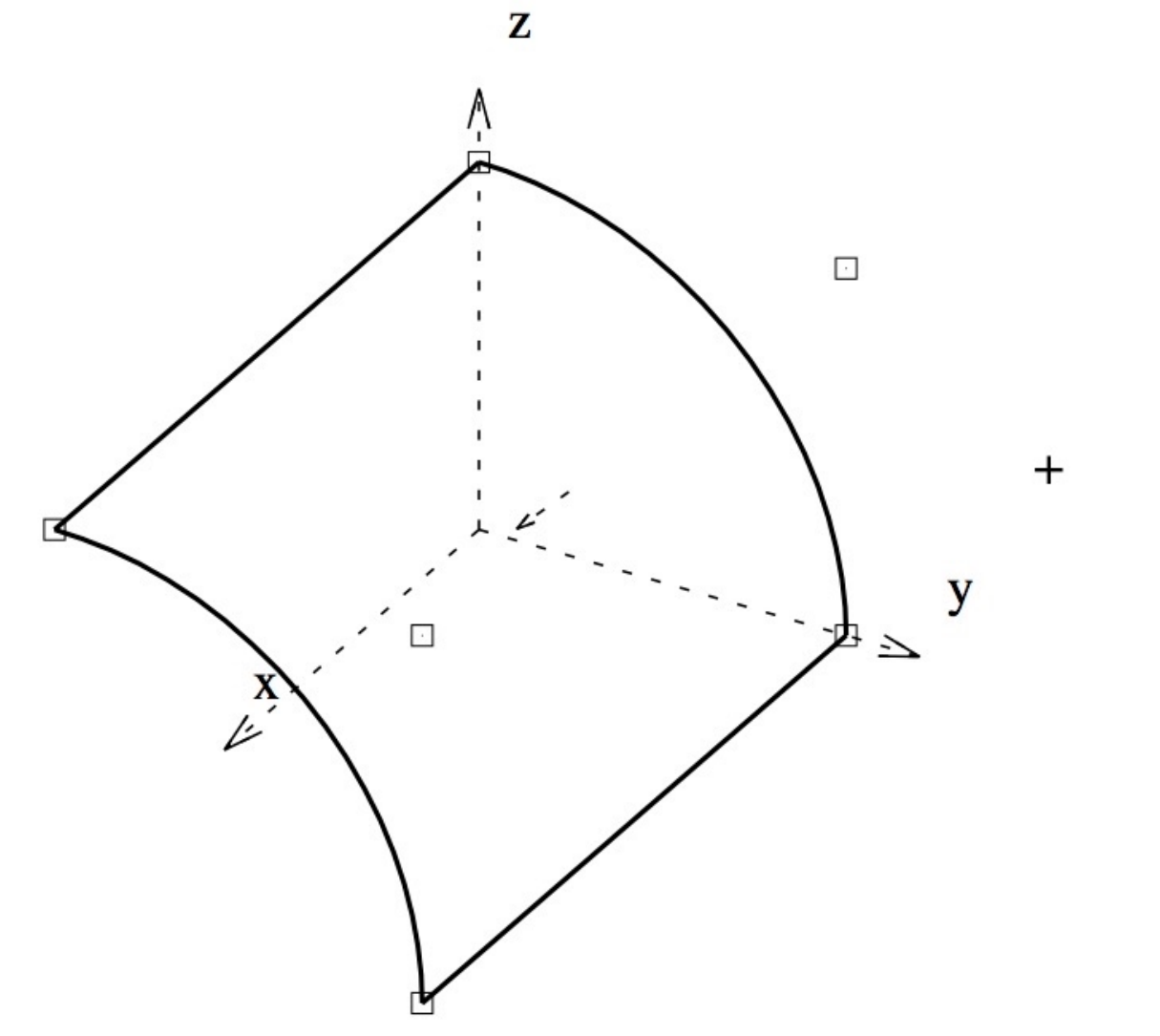}
\end{overpic}
\begin{overpic}[scale=0.35]{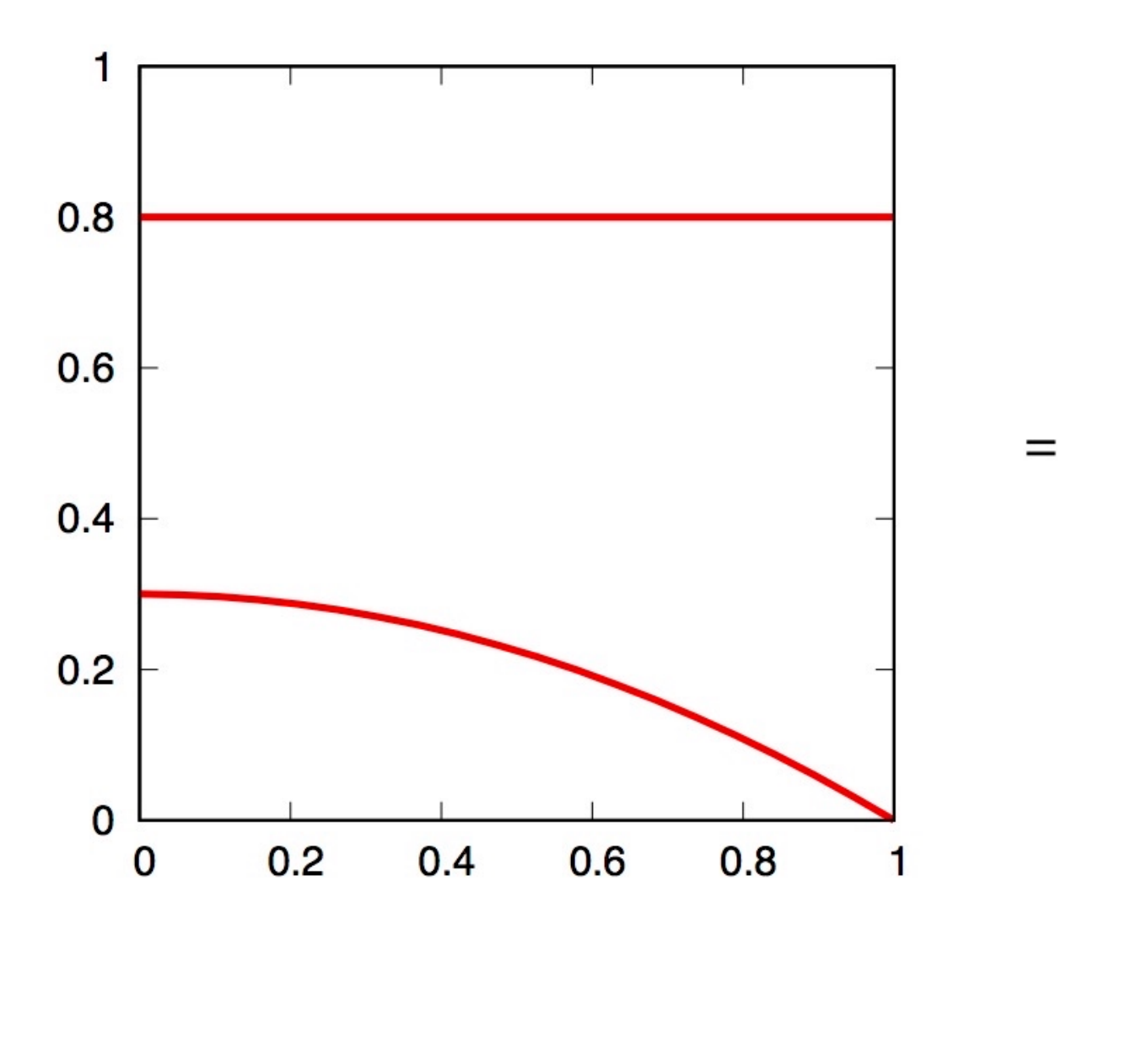}
 \put(40,10){$\uu$}
 \put(0,50){$\vv$}
\end{overpic}
\begin{overpic}[scale=0.25]{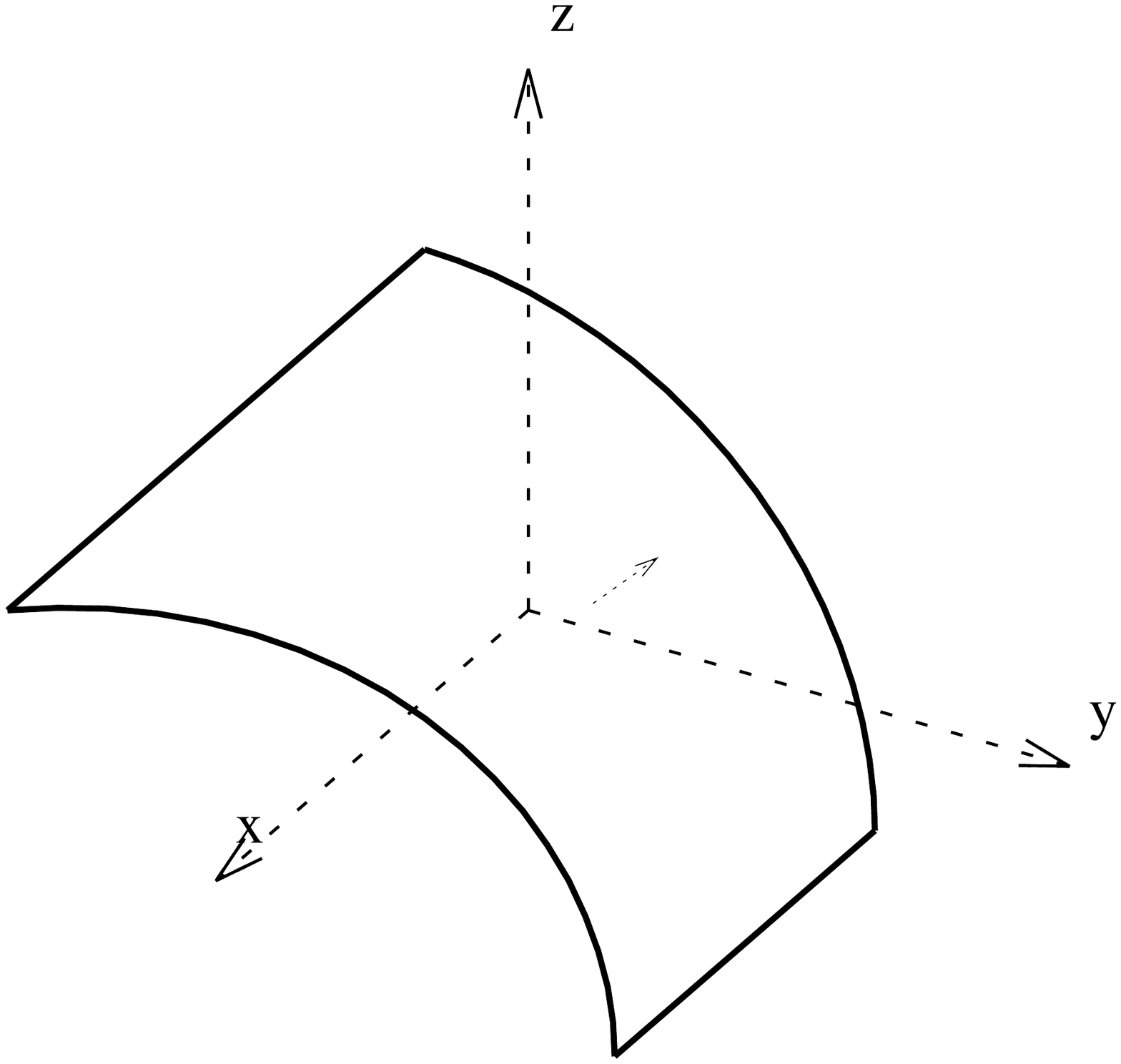}
\end{overpic}
\caption{Figure showing a normal patch, trimming curves and the resulting trimmed patch.}
\label{Trim}
\end{center}
\end{figure}

Patches can be trimmed  using trimming curves, resulting in more complex geometries. The trimming curves are defined in patch coordinates $\uu,\vv$.
A trimmed patch is shown in \myfigref{Trim}.
More information about trimming can be found in \cite{Beer2015b}.

\subsubsection{Special patches}
\begin{figure}
\begin{center}
\begin{overpic}[scale=0.6]{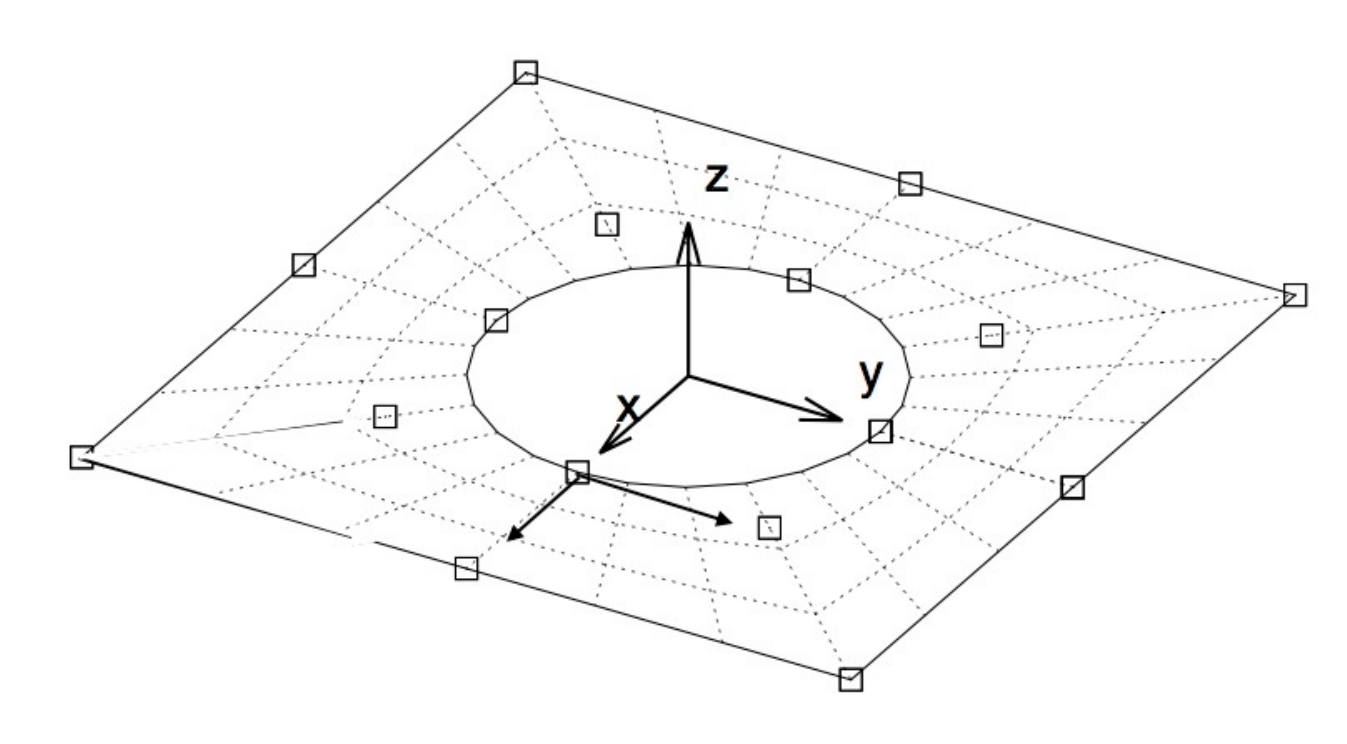}
 \put(35,17){$\vv$}
 \put(55,17){$\uu$}
\end{overpic}
\caption{Example of a special patch with a circular cutout.}
\label{SpecialPatch}
\end{center}
\end{figure}

Special patches can be used to model intersections of curved surfaces with flat surfaces with a minimum of effort. 
A special patch is shown in \myfigref{SpecialPatch}.

Its geometry is defined by:
\begin{equation}
\label{ }
\pt{x} (\xi,\eta)= (1-\vv) \pt{x}^{I} + \vv\pt{x}^{II} 
\end{equation}

where

\begin{align}
\pt{x}^{I} (\xi)= \sum_{j=1}^{J} \NURBS_{j} (\xi) \pt{x}^{I}_{j} 
&& \text{and} && \pt{x}^{II} (\xi)= \sum_{j=1}^{J} \NURBS_{j} (\xi) \pt{x}^{II}_{j}
\end{align}
The superscript \textit{II} refers to the inner and \textit{I} to the outer bounding curve.

The vectors in $\uu$ and $\vv$ directions are given by:

\begin{align}
\mathbf{ v}_{\xi}= (1-\vv) \mathbf{v}^{I} + \vv \mathbf{v}^{II} 
&& \text{and} && \mathbf{ v}_{\vv}= \pt{x}^{I} - \pt{x}^{II}
\end{align}

where

\begin{align}
\mathbf{v}^{I} (\xi)= \sum_{j=1}^{J} \frac{\partial \NURBS_{j}(\xi)}{\partial \xi}\pt{x}^{I}_{j} 
&& \text{and} && \mathbf{v}^{II} (\xi)= \sum_{j=1}^{J} \frac{\partial \NURBS_{j}(\xi)}{\partial \xi}\pt{x}^{II}_{j}
\end{align}
The outward normal and the Jacobian are computed the same way as for finite patches.

\subsubsection{Defining geometry with NURBS}
\label{DefNURBS}
CAD programs use NURBS to describe geometrical shapes (cylinder, spheres or general smooth shapes).
If two shapes intersect trimming is applied. Since CAD programs are mainly designed for visualisation, the resulting intersection geometry may have small gaps.
 The interested reader may consult \cite{BeerMarussig}, where a whole chapter is devoted to the topic of how data from CAD can be used for a BEM simulation.
However, instead of asking the CAD program to compute the intersection geometry and then extract the necessary data, which is complicated, one may write a MATLAB function that computes the intersection geometry without gaps. For the practical example we have developed a function that does this.

Here we show on an example how easy it is to define complex geometrical shapes with NURBS.
\begin{figure}
\begin{center}
\includegraphics[scale=0.25]{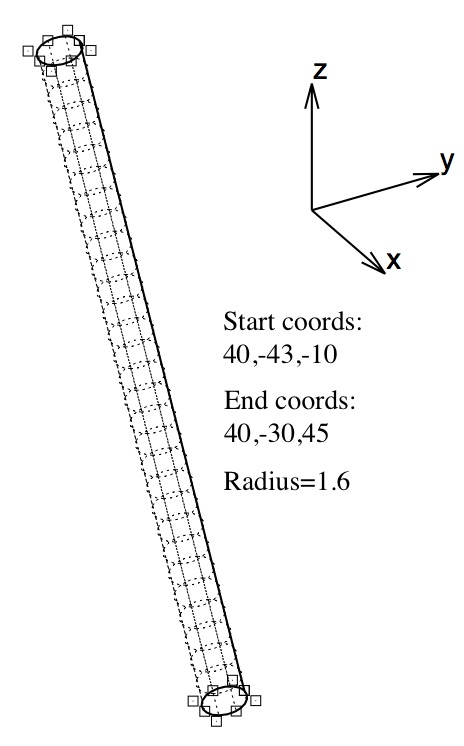}
\includegraphics[scale=0.25]{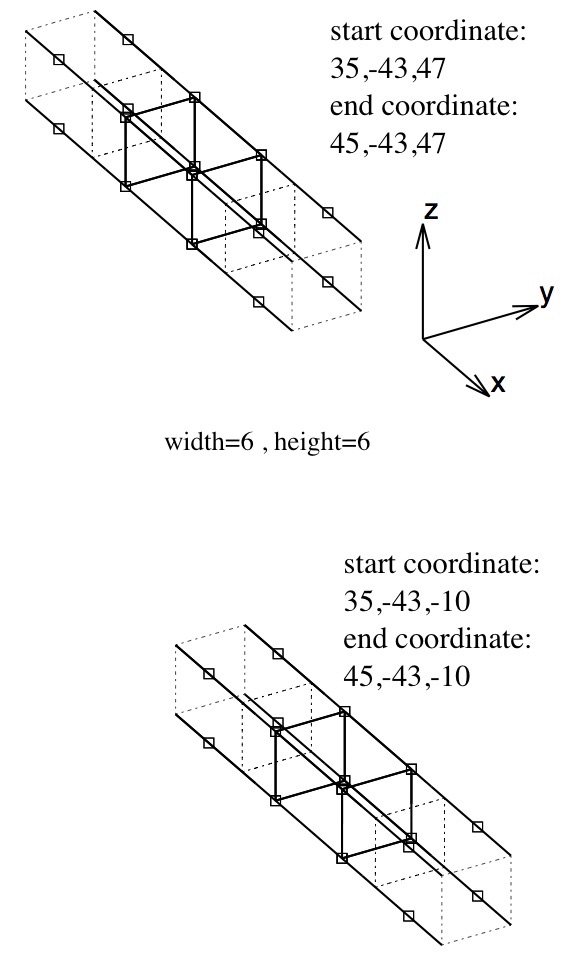}
\includegraphics[scale=0.5]{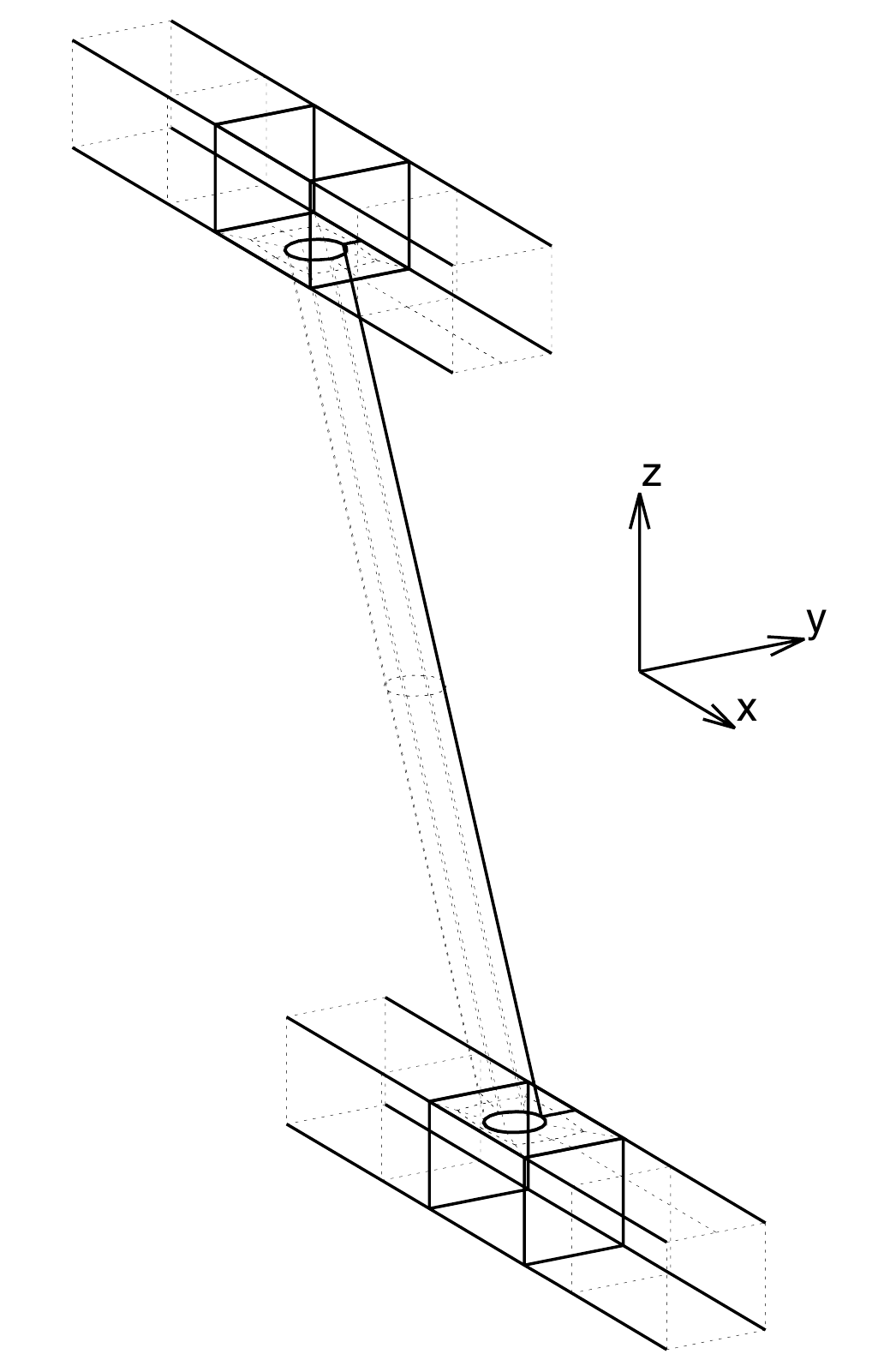}
\caption{Example of definition of geometry with NURBS, showing input data required. left: circular excavation, middle: cross-passages, right: combined geometry.}
\label{IGACAD}
\end{center}
\end{figure}
We start with the definition of a circular excavation and of 2 cross-passages. Only a few lines of input data (shown in \myfigref{IGACAD} on the left and middle) are required. The two geometries are then intersected resulting in the geometry definition on the right of \myfigref{IGACAD} consisting of a trimmed patch, normal patches, infinite patches and special patches. Note that no mesh generation is necessary. The dotted lines indicate integration regions (see the section on integration below). 

\subsubsection{Approximation of boundary values}
To be able to solve the patch integrals in \myeqref{Regular2}, the boundary values must be approximated.
For the approximation we also use NURBS basis functions.
For normal patches the unknown boundary values are approximated by
\begin{equation}
\begin{aligned}
\label{eq6:ApproxSecond}
 \myVec{\hat{\primary}}^{e} (\xi,\eta) &=  \sum_{k=1}^{K}   \hat{ \NURBS}_{k} (\xi,\eta) \   \myVec{\hat{\primary}}_{k}^{e}\\
 \myVec{\hat{\dual}}^{e}(\xi,\eta)  &=   \sum_{k=1}^{K}    \hat{\NURBS}_{k} (\xi,\eta) \  \myVec{\hat{\dual}}_{k}^{e}.
\end{aligned}
\end{equation}
where $\hat{ \NURBS}_{k} (\xi,\eta)$ are NURBS basis functions (the hat indicating that they may be different to the ones used for describing the geometry) and $ \myVec{\hat{\primary}}_{k}^{e}, \myVec{\hat{\dual}}_{k}^{e}$ are parameter values. It should be noted that in contrast to Lagrange polynomials, parameter values do not in general represent real values.

For infinite patches we have 2 choices for the displacements:
\begin{itemize}
  \item  \textit{Plane strain}: Displacements are constant to infinity:
 \begin{equation}
 \begin{aligned}
 \label{Approxinf}
  \myVec{\hat{\primary}}^{e} (\xi,\eta) =  \myVec{\hat{\primary}}^{e} (\xi,\eta=0)&=  \sum_{k=1}^{K^\infty}   \hat{ \NURBS}_{k} (\xi) \    \myVec{\hat{\primary}}_{k}^{e\infty}\\
 \end{aligned}
 \end{equation}
  \item \textit{Decay}: Displacements decay to zero as infinity is approached
   \begin{equation}
  \begin{aligned}
  \label{Approxinf}
  \myVec{\hat{\primary}}^{e} (\xi,\eta) =  (1- \vv)\myVec{\hat{\primary}}^{e} (\xi,\eta=0)&=  (1-\vv)\sum_{k=1}^{K^\infty}   \hat{ \NURBS}_{k} (\xi) \     \myVec{\hat{\primary}}_{k}^{e\infty}\\
 \end{aligned}
 \end{equation}
\end{itemize}
where $K^\infty$ is the number of parameters and $ \myVec{\hat{\primary}}_{k}^{e\infty}$ are the parameter values on the finite boundary.

Our refinement philosophy is to take the NURBS functions that define the geometry of the problem and refine them as necessary using knot insertion and order elevation.

Known values are defined by
\begin{equation}
\begin{aligned}
\label{}
 \myVec{\bar{\primary}}^{e}(\xi,\eta)  &=  \sum_{k=1}^{\bar{K}}    \bar{\NURBS}_{k}(\xi,\eta) \   \myVec{\bar{\primary}}_{k}^{e}\\
 \myVec{\bar{\dual}}^{e}(\xi,\eta)  &=   \sum_{k=1}^{\bar{K}}    \bar{\NURBS}_{k}(\xi,\eta) \  \myVec{\bar{\dual}}_{k}^{e}.
\end{aligned}
\end{equation}
where $\bar{\NURBS}_{k}$ are basis functions, which may be different from the ones defining the geometry and the unknown values. 

Inserting the approximations into the patch integrals allows the boundary parameters can be taken outside: 
\begin{eqnarray}
\label{ }
\int_{\Gamma_{e}}  \fund{U}(\sourcept_{n},\fieldpt) \ \myVec{\dual}^{e}(\fieldpt)\ d\Gamma_{e}(\fieldpt) = \int_{\Gamma_{e}}  \fund{U}(\sourcept_{n},\fieldpt) \left( \sum_{k=1}^{K}    \NURBS_{k} (\xi,\eta) \  \myVec{\hat{\dual}}_{k}^{e} \right)d\Gamma_{e}(\fieldpt)   \\
\nonumber
= \sum_{k=1}^{K} \left( \quad \int_{\Gamma_{e}}  \fund{U}(\sourcept_{n},\fieldpt)    \NURBS_{k} (\xi,\eta) \  d\Gamma_{e}(\fieldpt)  \right) \myVec{\hat{\dual}}_{k}^{e} \\
\int_{\Gamma_{e}}  \fund{T}(\sourcept_{n},\fieldpt) \ \myVec{\primary}^{e}(\fieldpt)\ d\Gamma_{e}(\fieldpt) = \int_{\Gamma_{e}}  \fund{U}(\sourcept_{n},\fieldpt) \left( \sum_{k=1}^{K}    \NURBS_{k} (\xi,\eta) \   \myVec{\hat{\primary}}_{k}^{e} \right) d\Gamma_{e}(\fieldpt) \\
\nonumber
= \sum_{k=1}^{K} \left( \quad \int_{\Gamma_{e}}  \fund{T}(\sourcept_{n},\fieldpt)    \NURBS_{k} (\xi,\eta) \  d\Gamma_{e}(\fieldpt)  \right) \myVec{\hat{\primary}}_{k}^{e}
\end{eqnarray}
where the hat and overbar has been omitted, because this depends if the values are known or unknown. This requires only the integration of fundamental solutions times the basis functions, which will be discussed later.

\subsection{Discretisation of volume integral.}
The volume integral is solved numerically by dividing the volume into inclusions, defining each one geometrically. The integral is replaced by a sum of integrations over inclusions:
\begin{equation}
\label{Vol}
 \begin{aligned}
\int_{\domain_{0}} \fund{E} (\sourcept_{n},\fieldpt)  \myVecGreek{\sigma}_{0} (\fieldpt) d \domain_{0} (\fieldpt) = \sum_{ni=1}^{Ni} \quad \int_{\domain_{ni}} \fund{E} (\sourcept_{n},\fieldpt)  \myVecGreek{\sigma}_{0} (\fieldpt) d \domain_{ni} (\fieldpt).
\end{aligned}
\end{equation}
where $Ni$ is the number of inclusions and $ \domain_{ni} $ specifies the inclusion domain. 
For the numerical treatment  an approximation of the initial stress $\myVecGreek{\sigma}_{0}$ is assumed inside the inclusion and the integrals are solved numerically or analytically as will be explained later.

Two types of inclusions are considered: General inclusion that represent a volume of material that has different elastic properties or behaves inelastically and linear inclusions to represent rock bolts. 
It should be noted that a $C^{0}$ continuity of displacements exists at the interface between the inclusion and the domain.

\subsubsection{Geometrical definition of general inclusion}
\begin{figure}[h]
\begin{center}
\begin{overpic}[scale=0.7]{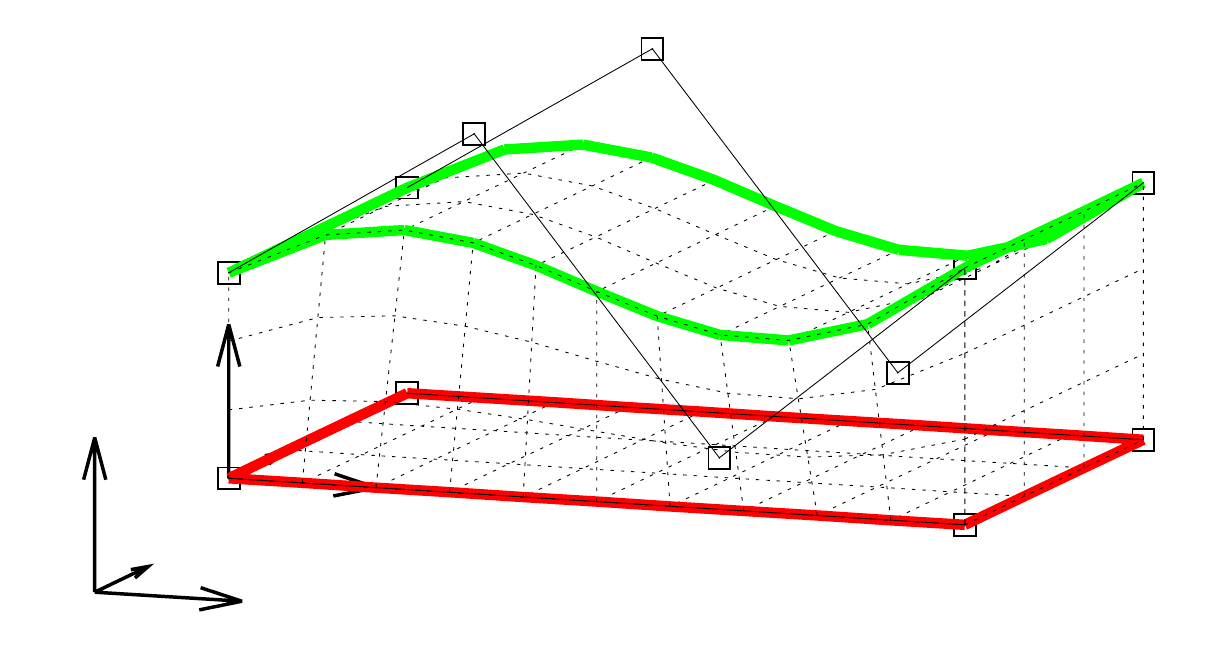}
 \put(20,5){$x$}
 \put(6,20){$z$}
 \put(10,10){$y$}
  \put(20,28){$r$}
 \put(30,15){$s$}
 \put(23,18){$t$}
\end{overpic}
\begin{overpic}[scale=0.5]{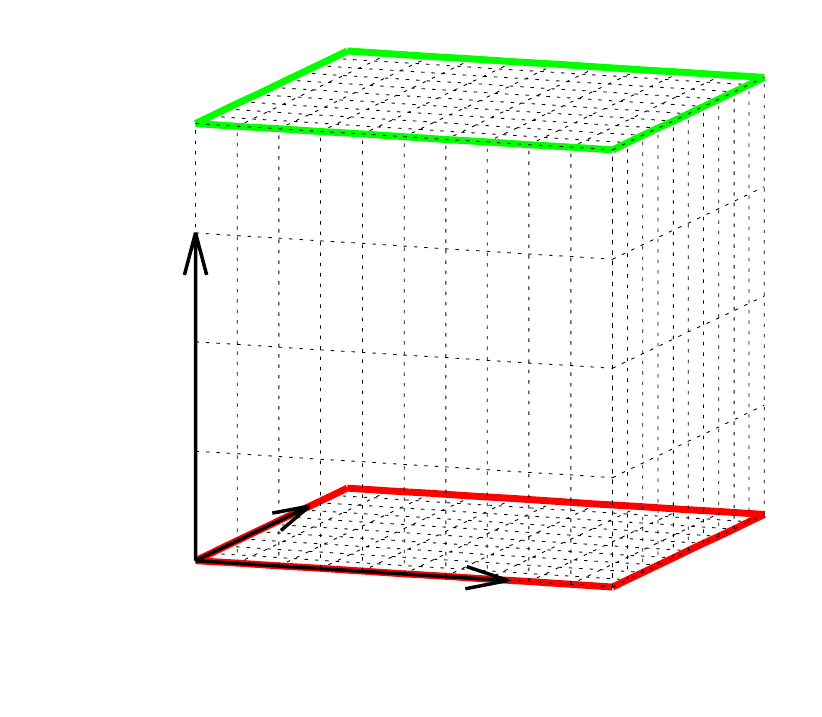}
  \put(60,20){$s$}
 \put(35,30){$t$}
 \put(15,60){$r$}
\end{overpic}
\caption{Mapping of 3-D inclusion showing the bottom and top NURBS surfaces and the associated control points defining the inclusion: Left in global $\pt{x}$, right in local $\pt{s}$ space. Also shown are subregions for the volume integration.}
\label{fig:Incl3D}
\end{center}
\end{figure}
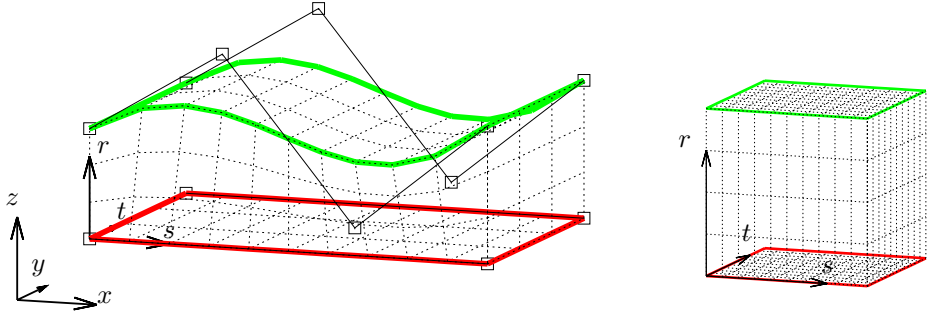

General inclusions are defined by bounding NURBS surfaces.
We establish a local coordinate system $\pt{s}=(s,t,r)^{\trans}=[0,1]^3$ as shown in \myfigref{fig:Incl3D} and map from local $\pt{s}$ coordinates to global $\pt{x}$ coordinates.
The global coordinates of a point $\pt{x}$ with the local coordinates $\pt{s}$ are given by
\begin{equation}
\pt{ x}({s,t,r})= (1-r) \ \pt{ x}^{I}(s,t) + {r} \ \pt{ x}^{II}({s,t})
\end{equation}

where

\begin{align}
  \pt{ x}^{I}({s,t})=\sum_{k=1}^{K^{I}} R_{k}^{I}({s,t}) \ \pt{ x}_{k}^{I} && \text{and} &&\pt{ x}^{II}({s,t})=\sum_{k=1}^{K^{II}} R_{k}^{II}({s,t}) \ \pt{ x}_{k}^{II} .
\end{align}

The superscript $I$ relates to the bottom (red) surface and $II$ to the top (green) bounding surface and $ \pt{ x}_{k}^{I} $, $ \pt{ x}_{k}^{II} $ are control point coordinates. $K^{I}$ and $K^{II}$ represent the number of control points, $R_{k}^{I}({s,t})$ and $R_{k}^{II}({s,t})$ are NURBS basis functions. Note that there is a one to one mapping between the local surface coordinates $\uu,\vv$ and the local coordinates $s,t$.

The derivatives are given by
\begin{equation}
\begin{aligned}
  \frac{\partial \pt{ x}({s,t,r})}{\partial {s}}&=& (1-{r}) \ \frac{\partial \pt{ x}^{I}({s,t})}{\partial {s}} &+&& {r} \ \frac{\partial \pt{ x}^{II}({s,t}) }{\partial {s}} \\
    \frac{\partial \pt{ x}({s,t,r})}{\partial {t}}&=& (1-{r}) \ \frac{\partial \pt{ x}^{I}({s,t})}{\partial {t}} &+&& {r} \ \frac{\partial \pt{ x}^{II}({s,t}) }{\partial {t}} \\
  \frac{\partial \pt{ x}({s,t,r})}{\partial {r}}&=& -\pt{ x}^{I}({s,t}) &+&&
  \ \pt{ x}^{II}({s,t})
\end{aligned}
\end{equation}

where for example:

\begin{align}
  \frac{\partial\pt{ x}^{I}({s,t})}{\partial {s}}=\sum_{k=1}^{K^{I}} \frac{\partial R_{k}^{I}({s,t})}{\partial {s}} \ \pt{ x}_{k}^{I} &&\text{and} &&
  \frac{\partial\pt{ x}^{II}({s,t})}{\partial {}s}=\sum_{k=1}^{K^{II}}
  \frac{\partial R_{k}^{II}({s,t})}{\partial {s}} \ \pt{ x}_{k}^{II}  .
\end{align}
The Jacobi matrix of this mapping is
\begin{equation}
\label{eq:jac3D}
\mathbf{J}=
\begin{pmatrix}
  \frac{\partial \pt{x}}{\partial {s}}   \\ \\
  \frac{\partial \pt{x}}{\partial {t}}  \\ \\
  \frac{\partial \pt{x}}{\partial {r}}  
\end{pmatrix}
\end{equation}
and the Jacobian is $J=| \mathbf{ J} |$.

\subsubsection{Geometry definition of linear inclusion}
\begin{figure}
\begin{center}
\begin{overpic}[scale=0.6]{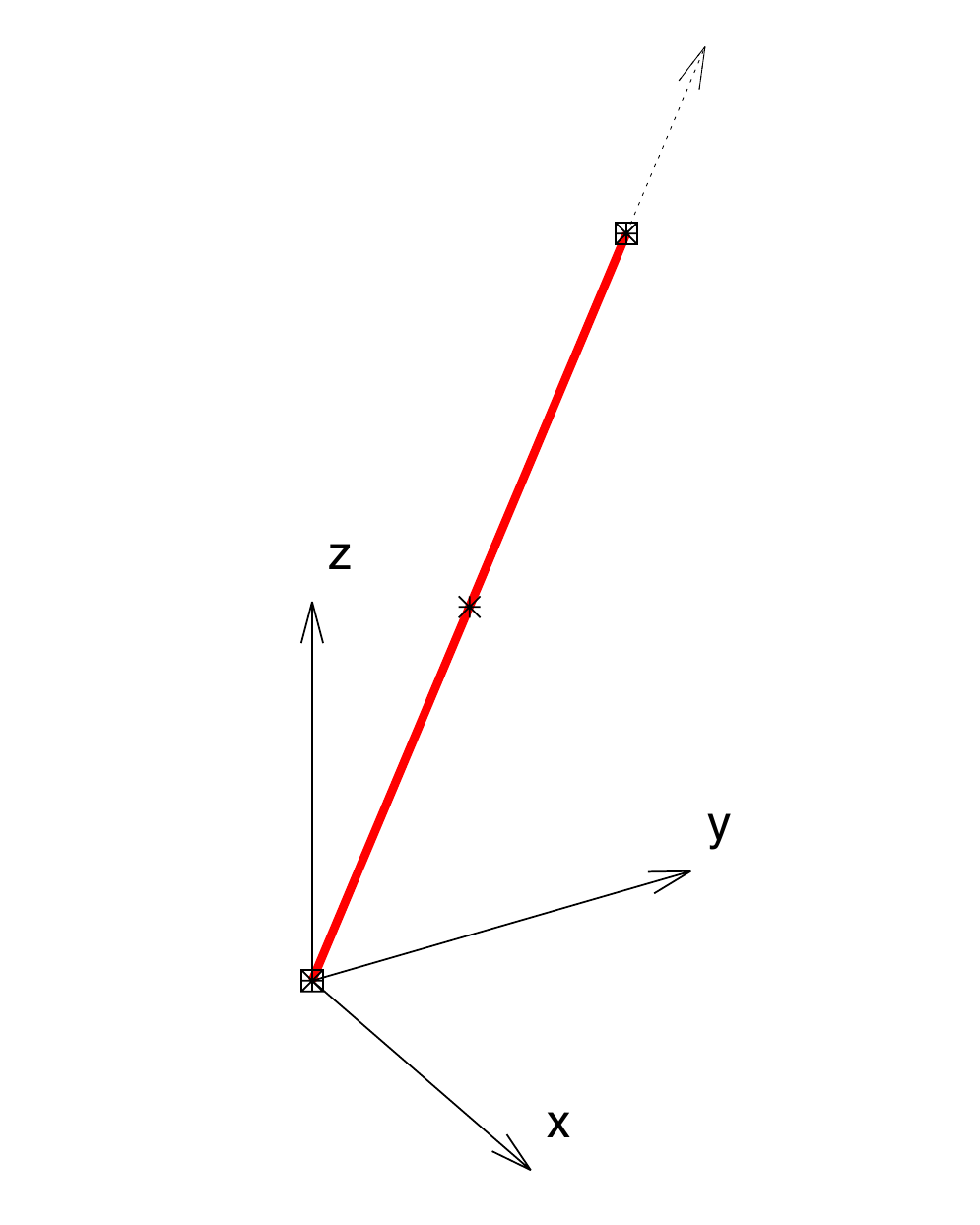}
\put(50,95){$z^{\prime}$}
\end{overpic}
\begin{overpic}[scale=0.6]{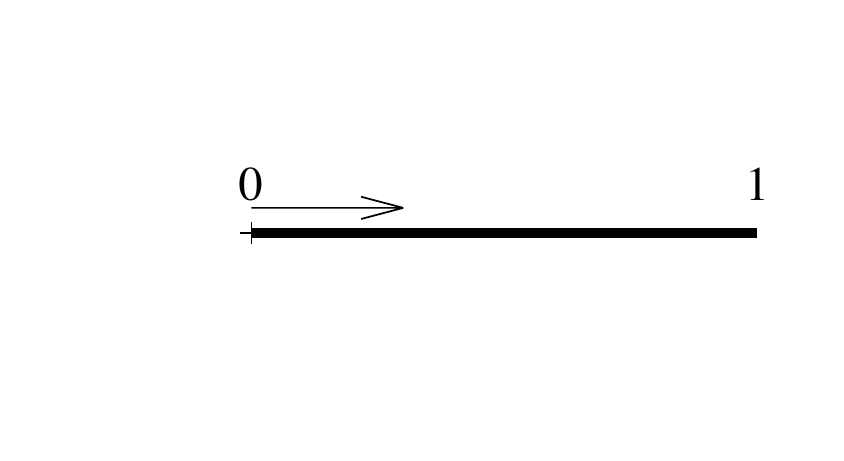}
\put(40,35){$s$}
\end{overpic}
\caption{Definition of linear inclusion by a NURBS curve with control points as hollow squares in global (left) and local (right) coordinates.}
\label{Rebar}
\end{center}
\end{figure}

This type of inclusions is used to model cables and rock bolts. Here we assume that the geometry is defined by a linear NURBS curve and that the bar has a circular cross-section with radius $R$ over which the stress and strain are assumed constant. \emph{The assumption is that the area of the cross-section of the inclusion is significantly smaller than that of the medium it is embedded in, allowing simplifications to be introduced for the integration.}
We establish a local coordinate system $s=\left[0,1\right]$ as shown on the right in \myfigref{Rebar}.
The global coordinates of a point $\pt{x}$ with the local coordinate $s$ are given by
\begin{eqnarray}
  \pt{ x} ({s})=\sum_{k=1}^{K} R_{k}(s) \ \pt{ x}_{k}  
\end{eqnarray}
where $K$ is the number of control points, $R_{k}(s)$ are NURBS basis functions and $\pt{ x}_{k}$ are control point coordinates.
We also define a local coordinate system whereby the $z^{\prime}$ axis is along the bar, specified by unit vector $\mathbf{ v}_{z^{\prime}}$. 
The vector along the bar is given by
\begin{eqnarray}
\mathbf{ V}_{z^{\prime}}= \frac{\partial\pt{ x}({s})}{\partial {s}}&=&\sum_{k=1}^{K} \frac{\partial R_{k}({s})}{\partial {s}} \ \pt{ x}_{k} 
\end{eqnarray}
The Jacobian is
\begin{equation}
\label{ }
J= \sqrt{\mathrm{V}_{z^{\prime}_x}^{2} + \mathrm{V}_{z^{\prime}_y}^{2} + \mathrm{V}_{z^{\prime}_z}^{2}}
\end{equation}
The unit vector in $z^{\prime}$ direction is given by
\begin{equation}
\label{ }
\mathbf{ v}_{z^{\prime}}= \mathbf{ V}_{z^{\prime}}/J
\end{equation}

\subsubsection{Approximation of initial stress}

For the numerical integration, introduced below, we need the values of the initial stress at Gauss points. 
It is inefficient to compute $\myVecGreek{\sigma}_{0}$ at every Gauss point, whose location changes according to the location of $\sourcept_{n}$.
Instead we compute the initial stress at a fixed grid of points inside the inclusion.
The value of initial stress at a point with the local coordinates $\pt{s}$ $=(s,t,r)^\mathrm{T}$ for general inclusions and $\pt{s}=s$ for linear inclusions can be obtained by interpolation between grid points:
\begin{equation}
\label{eq11:Interpol}
 \myVecGreek{\sigma}_{0}(\pt{s})= \sum_{l=1}^{L} M_{l}^{\sigma}(\pt{s})  \myVecGreek{\sigma}_{0l}
\end{equation}
where $  \myVecGreek{\sigma}_{0l} $ is the initial stress vector at grid point $l$ with the local coordinate $\pt{s}_{l}$. $L$ is the total number of inclusion points and $M_{l}^{\sigma}(\pt{s})$ are linear or constant basis functions, which will be shown later. 

\section{Numerical integration of boundary integrals}
The boundary integrals to be solved are:
\begin{eqnarray}
\label{ }
\nonumber
 \mathbf{ U}_{nk}^{e}= \int_{\Gamma_{e}}  \fund{U}(\sourcept_{n},\fieldpt)   \hat{ \NURBS}_{k} (\xi,\eta) \  d\Gamma_{e}(\fieldpt) \\
\mathbf{ T}_{nk}^{e} = \int_{\Gamma_{e}}  \fund{T}(\sourcept_{n},\fieldpt)    \hat{\NURBS}_{k} (\xi,\eta) \  d\Gamma_{e}(\fieldpt) \\
\nonumber
\mathbf{ T}_{n}^{e} = \int_{\Gamma_{e}}  \fund{T}(\sourcept_{n},\fieldpt)   \  d\Gamma_{e}(\fieldpt)
\end{eqnarray} 
They are evaluated numerically, using Gauss Quadrature
The integration scheme now depends on the location of the collocation point. If it is outside the patch we use regular integration otherwise we have to use singular integration.

Initially we divide the patch into \textbf{integration regions} depending on the following:
\begin{itemize}
  \item The location of the collocation points. They should be on integration region boundaries.
  \item The aspect ratios of each integration region. It should be moderate and this is particularly important for singular integration.
\end{itemize}

\subsection{Regular integration}
For regular integration we have to consider that the value of integrand tends to infinity as the collocation point is approached. To maintain an adequate precision of integration is crucial to the quality of the results. Therefore the number of Gauss points has to be increased near the collocation point. There is no analytical formula to determine the number of Gauss points required for a certain precision, but estimates have been worked out in \cite{BeerMarussig}. The number depends on the size of the integration region and the proximity of the collocation point. The best strategy is to limit the number of Gauss points available and to subdivide the integration region into subregions.
A Quadtree method that increases the number of Gauss points near the collocation point is most efficient.

Gauss Quadrature requires limits which range from -1 to +1. Therefore we introduce new local coordinates inside each subregion  $\bar{\myVecGreek{\xi}}=(\bar{\xi},\bar{\eta})^{\mathrm{T}}=[-1,1]^2$. The transformation to the patch coordinate system $\myVecGreek{\xi}=(\xi,\eta)^{\mathrm{T}}=[0,1]^2$ is given by:
\begin{eqnarray}
\uu  &=& \frac{ \triangle \uu_{s}}{2} (1+\bar{\uu}) + \uu_{s1}  \\
\vv  &=& \frac{ \triangle \vv_{s}}{2} (1+\bar{\vv}) + \vv_{s1} 
\end{eqnarray}
where $ \triangle \uu_{s} \times \triangle \vv_{s}$ is the size of the subregion and $\uu_{s1}, \vv_{s1}$ are the starting coordinates. 

The integration can be written as:
\begin{eqnarray}
\label{eq:regularNumericalIntegration}
 \nonumber
 \mathbf{ U}_{nk}^{e}& = & \sum_{s=1}^{S} \sum_{i=1}^{g_{\uu}(s)}   \sum_{j=1}^{g_{\vv}(s)} \fund{U}\left( \sourcept_{n}, \fieldpt(\bar{\uu}_{i}, \bar{\vv}_{j}) \right) \ \NURBS_{k} (\uu (\bar{\uu}_{i}), \vv (\bar{\vv}_{j}))  \frac{\triangle \uu_{s}}{2} \frac{\triangle \vv_{s}}{2} J(\bar{\uu}_{i}, \bar{\vv}_{j}) W_{i} W_{j}\\
 \mathbf{ T}_{nk}^{e}& = & \sum_{s=1}^{S} \sum_{i=1}^{g_{\uu}(s)}   \sum_{j=1}^{g_{\vv}(s)} \fund{T}\left( \sourcept_{n}, \fieldpt(\bar{\uu}_{i}, \bar{\vv}_{j}) \right) \ \NURBS_{k} (\uu (\bar{\uu}_{i}), \vv (\bar{\vv}_{j}))  \frac{\triangle \uu_{s}}{2} \frac{\triangle \vv_{s}}{2} J(\bar{\uu}_{i}, \bar{\vv}_{j}) W_{i} W_{j}\\
  \nonumber
  \mathbf{ T}_{n}^{e}& = & \sum_{s=1}^{S} \sum_{i=1}^{g_{\uu}(s)}   \sum_{j=1}^{g_{\vv}(s)} \fund{T}\left( \sourcept_{n}, \fieldpt(\bar{\uu}_{i}, \bar{\vv}_{j}) \right)  \  \frac{\triangle \uu_{s}}{2} \frac{\triangle \vv_{s}}{2} J(\bar{\uu}_{i}, \bar{\vv}_{j}) W_{i} W_{j}
 \end{eqnarray}
 where $g_{\uu}(s)$, $g_{\vv}(s)$  is the number of  Gauss points in $\uu, \vv$ directions, $W_{i}, W_{j}$ are Gauss weights, $S$ is the number of subegions and $J$ is the Jacobian of the transformation from global coordinates $\pt{x}$ to local $\myVecGreek{\xi}$ coordinates.
 
\subsection{Singular integration}
\begin{figure}
\begin{center}
\begin{overpic}[scale=0.6]{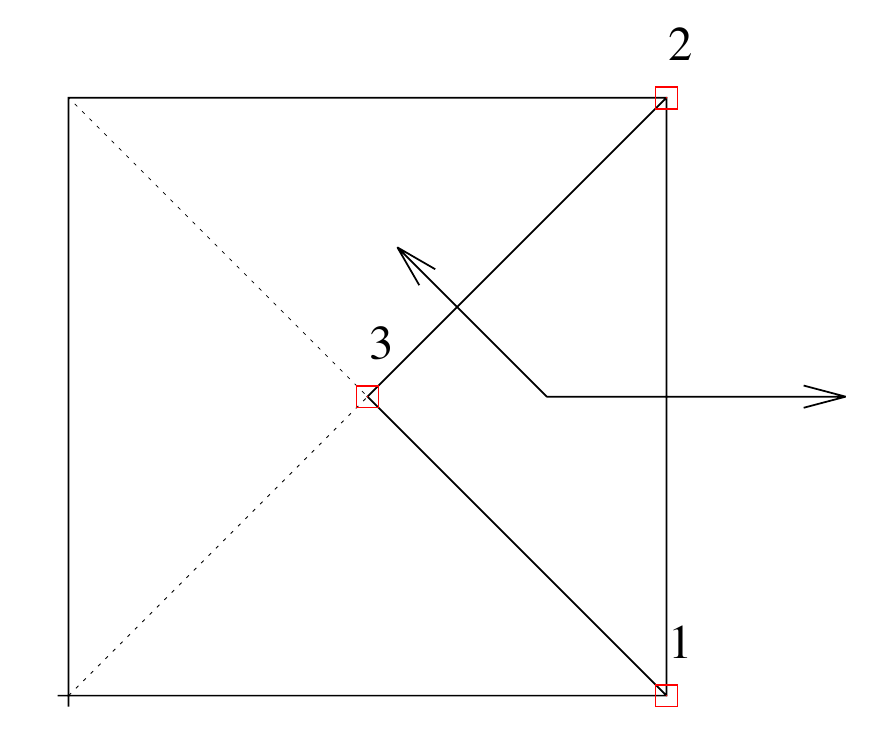}
\put(100,40){$\bar{\uu}$}
\put(50,60){$\bar{\vv}$}
\end{overpic}
\caption{Definition of a triangular subregion. The collocation point is located at point 3.}
\label{SingB}
\end{center}
\end{figure}

If the collocation point is part of the subregion then the integral involving $\fund{U}$ is weakly singular. The integral is solved by subdividing the integration region into triangular subregions with the collocation point at the apex (\myfigref{SingB}). This means that the Jacobian tends to zero as the collocation point is approached.
The singular integration can now be written as
\begin{equation}
 \triangle \mathbf{ U}_{nk}^{e} =  \sum_{n_{t}=1}^{N_{t}} \sum_{i=1}^{g_{\uu}}   \sum_{j=1}^{g_{\vv}} \fund{U}\left( \sourcept_{n}, \fieldpt(\bar{\uu}_{i}, \bar{\vv}_{j}) \right)  \NURBS_{k} (\uu (\bar{\uu}_{i}), \vv (\bar{\vv}_{j}))   0.25 \ J_{\uu,n_{t}}(\bar{\uu}_{i}, \bar{\vv}_{j})   J(\bar{\uu}_{i}, \bar{\vv}_{j}) W_{i} W_{j}
\end{equation}
where $N_{t}$ is the number of triangles. There are now two Jacobians involved, one for the transformation from the patch coordinates to triangular coordinates ($J_{\uu,n_{t}}(\bar{\uu}_{i}, \bar{\vv}_{j})$) which tends to zero and one for the transformation from patch to global coordinates ($J(\bar{\uu}_{i}, \bar{\vv}$)).

\section{Integration of integrals over $\domain_{0}$}
For the inclusion integrals we use numerical integration for general inclusions and analytical integration for linear inclusions.
For the numerical integration we subdivide the inclusion into integration regions with the same criteria as used for the boundary integration and apply Gauss quadrature.
When point $n$ is part of the integration region we have to invoke singular integration, if it is not a regular one.

\subsection{General inclusions}
The integral to be solved for each inclusion is :
\begin{equation}
\mathbf{ B}^{ni}_{0n}  = \int_{\domain_{ni}}  \fund{E} (\sourcept_{n},\fieldpt)   \myVecGreek{\sigma}_{0} (\fieldpt) d \domain_{ni} (\fieldpt) 
\end{equation}
Introducing the interpolation of initial stresses (\ref{eq11:Interpol}) we have:
\begin{equation}
\label{ }
\mathbf{ B}^{ni}_{0n} =  \int_{\domain_{ni}} \fund{E} (\sourcept_{n},\fieldpt)  \sum_{j=1}^{J} M_{j}^{\sigma}(\fieldpt) \myVecGreek{\sigma}_{0j} = \sum_{j=1}^{J} \mathbf{ B}^{ni}_{0nj} \myVecGreek{\sigma}_{0j}
\end{equation}
where
\begin{equation}
\mathbf{ B}^{ni}_{0nj}  = \int_{\domain_{ni}}  \fund{E} (\sourcept_{n},\fieldpt)  M_{j}^{\sigma}(\fieldpt) d \domain_{ni} (\fieldpt) 
\end{equation}

\subsubsection{Regular integration}
To maintain adequate precision of integration we subdivide the integration region into subregions depending on the size of the integration region and the proximity of point $n$. 
For sub-region $n_{s}$ the transformation from the inclusion ($\pt{s}$) coordinates to the coordinates used for Gauss integration $\bar{\myVecGreek{\xi}}=(\bar{\xi},\bar{\eta},\bar{\zeta})^{\mathrm{T}}=[-1,1]^3$ is  given by
\begin{eqnarray}
\nonumber
s & = \frac{\Delta s_{n}}{2} (1+\bar{\xi}) + s_{n_s} \\
t & = \frac{\Delta t_{n}}{2} (1+\bar{\eta}) + t_{n_s} \\
\nonumber
r& = \frac{\Delta r_{n}}{2} (1+\bar{\zeta}) + r_{n_s}
\end{eqnarray}
where $\Delta s_{n}\times \Delta t_{n} \times \Delta r_{n}$ denotes the size of the sub-region and $s_{n},t_{n},r_{n}$ are the edge coordinates.
The Jacobian of this transformation is $J_{\xi}^{n}=\frac{1}{8}\ \Delta s_{n} \  \Delta t_{n} \ \Delta r_{n}$.

We can write:
\begin{equation}
  \label{GaussSecond}
   \mathbf{ B}_{0nj}^{ni} = \sum_{n_{s}=1}^{N_{s}}\int_{-1}^{1} \int_{-1}^{1} \int_{-1}^{1} \fund{E} \left( \sourcept_{n},\bar{\pt{x}}(\bar{\xi},\bar{\eta}, \bar{\zeta}) \right)
 M_{j}^{\sigma} \left( \bar{\pt{x}} (\bar{\xi},\bar{\eta}, \bar{\zeta}) \right) J(\pt{s}) \ J_{\xi}^{n_{s}} \ d \bar{\xi} d \bar{\eta} d \bar{\zeta}
\end{equation}
where $ J(\mathbf{ s})$ is the Jacobian of the mapping between $\pt{s}$ and $\pt{x}$ coordinate systems.

Applying Gauss integration we have:
\begin{equation}
  \label{GaussintV}
   \mathbf{B}_{0nj}^{ni}  \approx \sum_{n_{s}=1}^{N_{s}} \sum_{g_{s}=1}^{G_{s}} \sum_{g_{t}=1}^{G_{t}}  \sum_{g_{r}=1}^{G_{r}}\fund{E}\left( \sourcept_{n},\bar{\pt{x}}(\bar{\xi}_{g_{s}},\bar{\eta}_{g_{t}}, \bar{\zeta}_{g_{r}}) \right) M_{j}^{\sigma} \left( \bar{\pt{x}}(\bar{\xi}_{g_{s}},\bar{\eta}_{g_{t}}, \bar{\zeta}_{g_{r}}) \right) J(\pt{s})  \ J_{\xi}^{n_{s}} \ W_{g_{s}} \ W_{g_{t}} \ W_{g_{r}}
\end{equation}
where $N_{s}$ is the number of subregions and $G_{s},G_{t}$ and $G_{r}$ are the number of Gauss points (which depends on the proximity of $n$) and $\bar{\xi}_{g_{s}},\bar{\eta}_{g_{t}}, \bar{\zeta}_{g_{r}}$ the Gauss point coordinates in $s, t$ and $r$ directions, respectively. $W_{g_{s}} \ W_{g_{t}} \ W_{g_{r}}$ are Gauss weights. 

\subsubsection{Singular integration}

If the integration region includes the point $\pt{x}_{n}$, then the integrand tends to infinity as the point is approached.
To deal with the integration involving the weakly singular Kernel we compute the Gauss points in a local coordinate system, where the Jacobian tends to zero as the singularity point is approached. Singular integration of general inclusions is discussed in detail in \cite{BeerMarussig}.

\subsection{Linear inclusion, reinforcement bar}
For the linear inclusions we can apply analytical integration.
We model the bar as a cylindrical region with radius $R$ and length $H$ and assume the initial stress to be in the local $z^{\prime}$ direction and to vary linearly along the bar. Since the initial stresses are computed from the strains, this means that the displacements along the bar can have a quadratic variation.
We consider two types of integration: one where point $\sourcept$ is outside the inclusion (regular integration) and one where it is not (singular integration). 

\subsubsection{Analytical computation of regular integral }

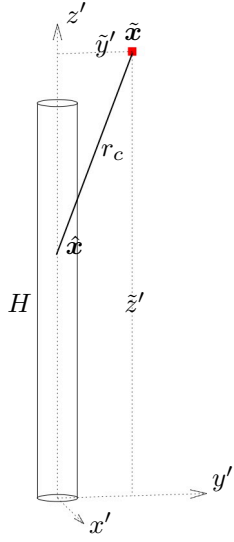
\begin{figure}
\begin{center}
\begin{overpic}[scale=0.5]{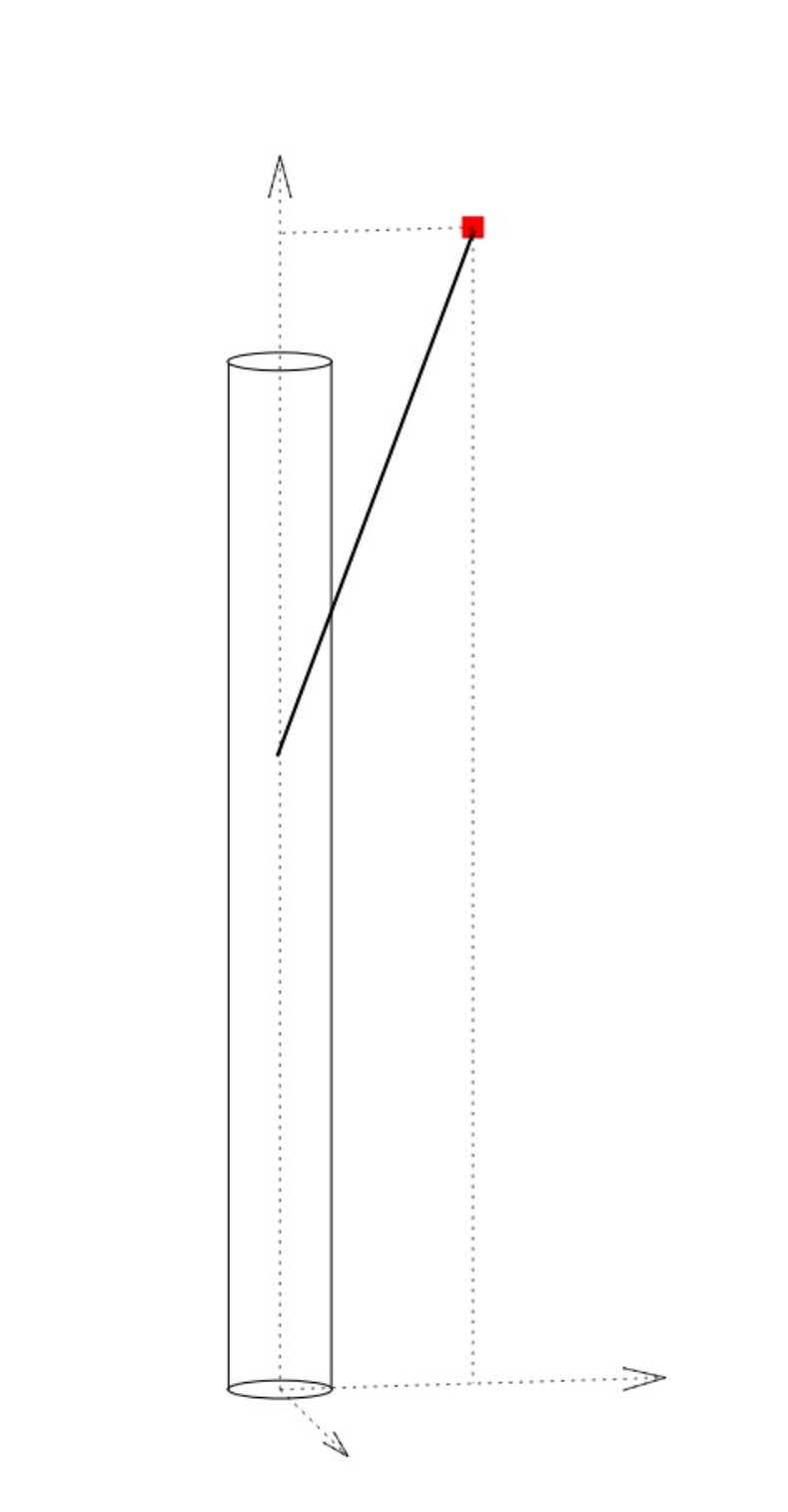}
\put(24,2){$x^{\prime}$}
\put(45,10){$y^{\prime}$}
\put(20,90){$z^{\prime}$}
\put(10,40){$H$}
\put(30,87){$\sourcept$}
\put(20,50){$\fieldpt$}
\put(30,40){$\tilde{z}^{\prime}$}
\put(25,85){$\tilde{y}^{\prime}$}
\put(26,67){$r_{c}$}
\end{overpic}
\caption{Analytical computation of regular integral for a subregion of length $H$ of a linear inclusion.}
\label{AnalR}
\end{center}
\end{figure}

Since we assume that the cross-sectional area is significantly smaller than the surrounding medium we can assume that $\fund{E}$ is constant over the cross-section. In addition we note that the result will multiply with the initial stresses in local directions ( $\{ \myVecGreek{\sigma^{\prime}}_{0} \}$). The integral to be solved is therefore:
\begin{equation}
\mathbf{ B}^{'}_{0nl}  = \int_{\domain} \frac{1}{r_{c}^{2}} \fund{\tilde{E}}^{\prime} M_{l}^{\sigma} d \domain_{i} (\fieldpt) 
\end{equation}
where the prime indicates that the result is computed in the local  $x^{\prime}, y^{\prime},z^{\prime}$ coordinate system (Fig. \ref{AnalR}). $r_{c}$ is the distance between the source point and a point on the axis of the inclusion .

The integral to be solved is:
\begin{equation}
\label{ }
\triangle \fund{E}^{\prime}_{ij}=   -C \int_V\ \frac{1}{r_{c}^{2}} \left[C_{3}(r^{\prime}_{,k}\delta_{ij} + r^{\prime}_{,j}\delta_{ik}) - r^{\prime}_{,i}\delta_{jk} + C_{4} \ r^{\prime}_{,i} r^{\prime}_{,j} r^{\prime}_{,k}\right] M_{l}^{\sigma}  d V
\end{equation}
We choose the local axes such that $\tilde{x}^{\prime}=0$ as follows. The vector pointing in the $x^{\prime}$ direction is given by:
\begin{equation}
\label{ }
\mathbf{V}_{x^{\prime}}= (\sourcept - \fieldpt) \times \mathbf{ v}_{z^{\prime}}
\end{equation}
and the one in $y^{\prime}$ direction is:
\begin{equation}
\label{ }
\mathbf{V}_{y^{\prime}}= \mathbf{ v}_{z^{\prime}} \times \mathbf{ v}_{x^{\prime}}
\end{equation}
where the capital letter indicates that the vector is not normalised.

If point $\sourcept$ is along the axis of the bar this computation does not work and then we assume
\begin{equation}
\label{ }
\mathbf{v}_{x^{\prime}}= \mathbf{ v}_{y} \times \mathbf{ v}_{z^{\prime}}
\end{equation}
where $\mathbf{ v}_{y}$ is a vector in global $y$-direction.

For the computation of the fundamental solution we have:
\begin{equation}
\label{ }
r_1 = 0 \quad r_2 = - \tilde{y}^{\prime} \quad r_3 = z^{\prime}-\tilde{z}^{\prime} \quad r = r_c = \sqrt{\tilde{y}^{\prime 2} + (z^{\prime}-\tilde{z}^{\prime})^2}
\end{equation}
and 
\begin{equation}
\label{ }
r^{\prime}_{,1}= 0 \quad r_{,2}= - \frac{\tilde{y}^{\prime}}{r_c} \quad r^{\prime}_{,3}= \frac{z^{\prime}-\tilde{z}^{\prime}}{r_c} \quad \quad dV = \pi R^2 d\tilde{z}^{\prime}
\end{equation}
The integral to be solved is:
\begin{equation}
\label{ }
\triangle \fund{E}^{\prime}_{ij}= \pi R^2 C \int_{z^{\prime}=0}^{H} \ \frac{1}{r_c^2}  \left[C_{3}(r_{,k}\delta_{ij} + r_{,j}\delta_{ik}) - r_{,i}\delta_{jk} + C_{4} \ r_{,i} r_{,j} r_{,k}\right]    M_{l}^{\sigma } (z^{\prime})
dz^{\prime} 
\end{equation}
where the linear interpolation functions are given by:

\begin{equation}
\label{ }
M_{1}^{\sigma }(z^{\prime}) = \frac{ z^{\prime}}{H} \quad  M_{2}^{\sigma }= 1 -\frac{z^{\prime}}{H}
\end{equation}

The analytical solution in Voigt notation is provided in Appendix A.

Since the result of the multiplication with $\myVecGreek{\sigma^{\prime}}_{0}$ has to be in global coordinates a transformation to the global system is necessary:
\begin{equation}
\label{ }
\mathbf{ B}_{0nj}= \mathbf{ T} \mathbf{ B}_{0nj}^{\prime}
\end{equation}
where $\mathbf{ T}$ is the transformation matrix given by:
\begin{equation}
\label{ }
\mathbf{ T}= \left(\begin{array}{ccc}v_{x^{\prime}_{x}} & v_{y^{\prime}_{x}} & v_{z^{\prime}_{x}} \\v_{x^{\prime}_{y}} & v_{y^{\prime}_{y}} & v_{z^{\prime}_{y}} \\v_{x^{\prime}_{z}} & v_{y^{\prime}_{z}} & v_{z^{\prime}_{z}}\end{array}\right)
\end{equation}

\subsubsection{Analytical computation of singular integral}
\begin{figure}
\begin{center}
\begin{overpic}[scale=0.5]{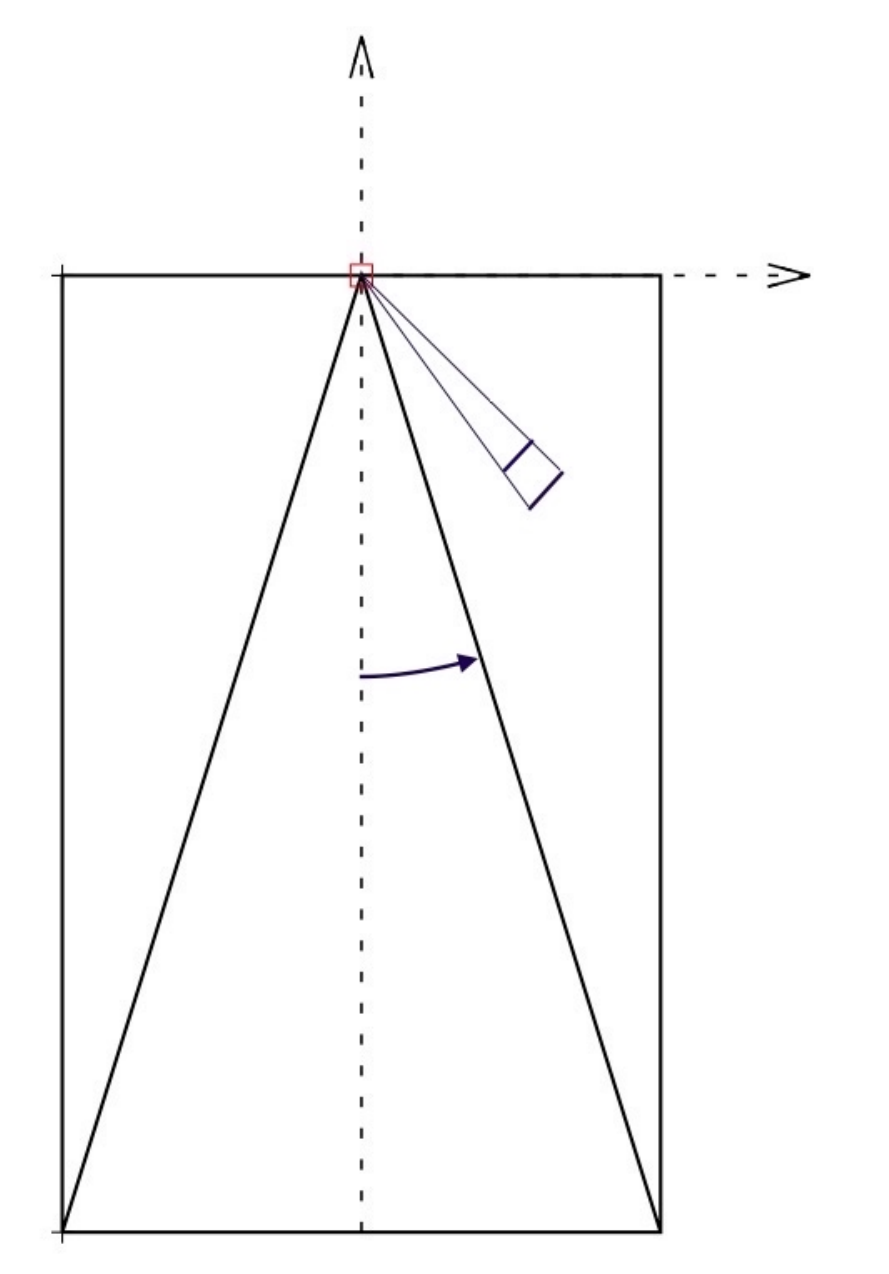}
\put(36,72){$r$}
\put(43,65){$dr$}
\put(31,55){$\theta$}
\put(40,5){$R$}
\put(44,60){$ r d \theta $}
\put(30,30){1}
\put(10,50){2}
\put(55,40){$H$}
\put(31,97){$z^{\prime}$}
\put(60,85){$y^{\prime}$}
\end{overpic}
\begin{overpic}[scale=0.4]{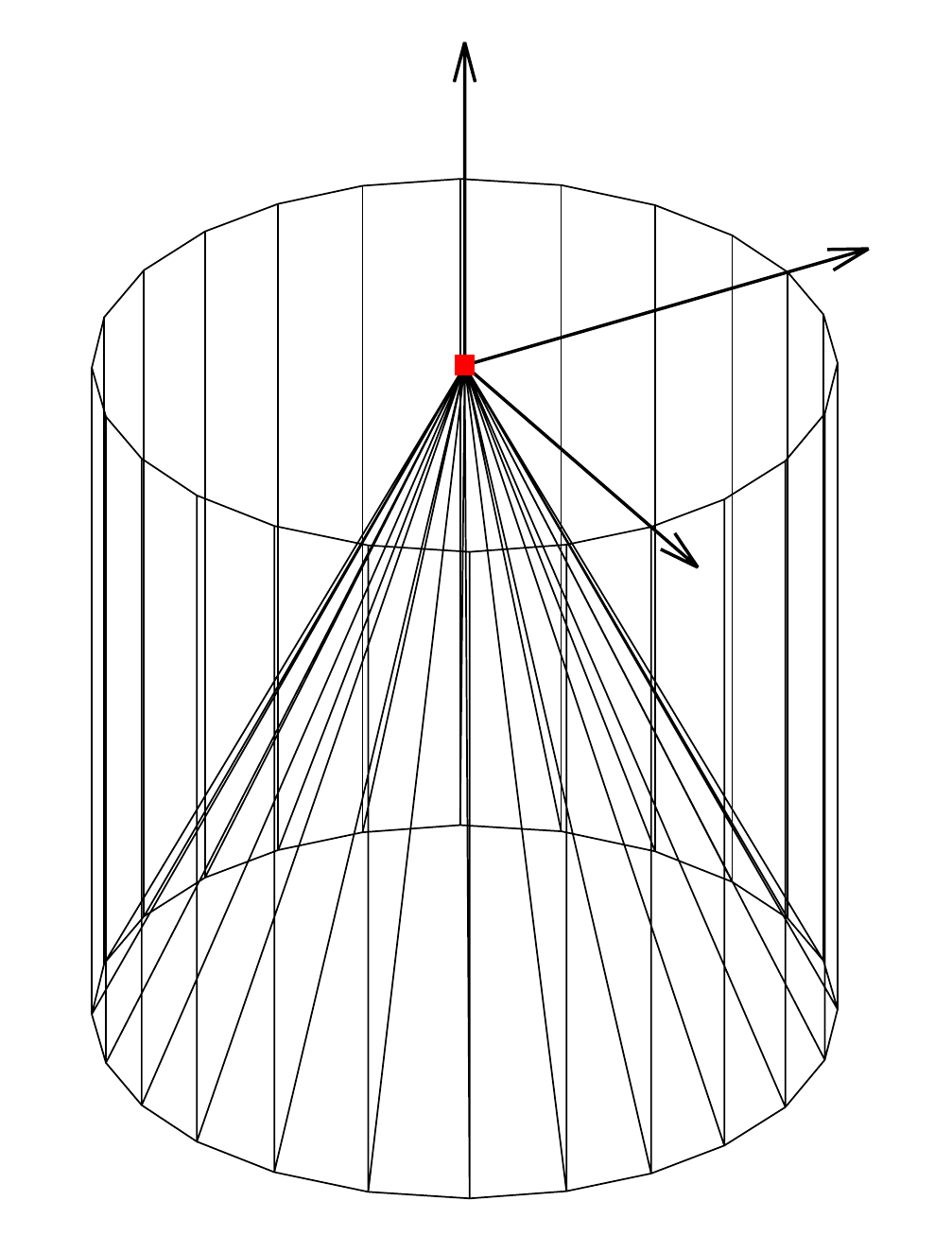}
\put(31,97){$z^{\prime}$}
\put(60,85){$y^{\prime}$}
\put(50,50){$x^{\prime}$}
\end{overpic}
\caption{Explanation of singular integration by subdivision into conical subregions. Left: section through bar, right: Axonometric view. Singular point is marked by red square}
\label{Sing2}
\end{center}
\end{figure}
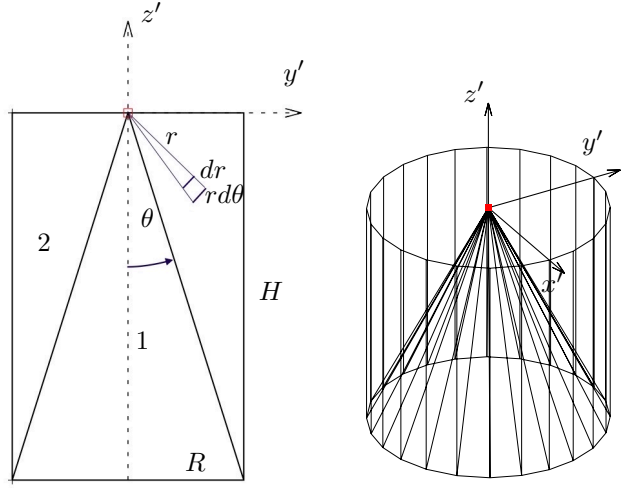

Referring to Fig. \ref{Sing2} we subdivide the bolt into 2 subregions and obtain the following integrals in polar coordinates:
\begin{equation}
\label{I1I2}
\triangle \fund{E^{\prime}}_1  =  \int_{\phi=0}^{2\pi} \int_{\theta=\pi-\tilde{\theta}}^{\pi} \int_{r=0}^{\frac{H}{\cos{(\pi-\theta)}}}  \frac{1}{r^{2}} \fund{\tilde{E}^{\prime}} \sin{\theta}  dr \ r^{2} M_{l}^{\sigma} d\theta d\phi 
\end{equation}
\begin{equation}
\triangle \fund{E^{\prime}}_2  =  \int_{\phi=0}^{2\pi} \int_{\theta=\pi/2}^{\pi-\tilde{\theta}} \int_{r=0}^{\frac{R}{\sin{(\pi-\theta)}}} \frac{1}{r^{2}}  \fund{\tilde{E}^{\prime}} \sin{\theta}  dr \ r^{2}  M_{l}^{\sigma}  d\theta d\phi \nonumber
\end{equation}
with $\tilde{\theta}= \arctan (R/H)$. It can be seen that the $r^{2}$ term cancels out which means that the integrand is no longer singular.

The shape functions can be expressed in cylindrical coordinates in two different formats, depending on the position of the collocation point, i.e.:
\begin{eqnarray}
M_{1}^{\sigma }(r,\theta) & = & - r \frac{\cos\theta}{H}  \\
M_{2}^{\sigma }(r,\theta) & = & 1 + r \frac{\cos\theta}{H}
\end{eqnarray}
if the collocation point is on the top of the cylinder,
\begin{eqnarray}
M_{1}^{\sigma }(r,\theta) & = &  r \frac{\cos\theta}{H}  \\
M_{2}^{\sigma }(r,\theta) & = & 1 - r \frac{\cos\theta}{H}
\end{eqnarray}
if the collocation point is on the bottom of the cylinder. In such a way we have that:
\begin{equation}
\triangle \fund{E^{\prime}}^{\sourcept_{n}=top}(j,k) = - \triangle \fund{E^{\prime}}^{\sourcept_{n}=bottom}(j,k)
\end{equation}

The terms  of $\triangle \fund{E}^{\prime} = \triangle \fund{E}^{\prime}_1+\triangle \fund{E}^{\prime}_2$, in Voigt notation, different from zero are provided in Appendix B.

As before a transformation to the global system is necessary:
\begin{equation}
\label{ }
\mathbf{ B}_{0nj}= \mathbf{ T} \mathbf{ B}_{0nj}^{\prime}
\end{equation}

\section{System of equations}
The discretised integral equations can be written in matrix form as:
\begin{equation}
\label{DisIE}
[\mathbf{ L}] \{\mathbf{ x}\} = \{\mathbf{ r}\} + [\mathbf{ B}_{0}] \{ \myVecGreek{\sigma}_{0}\} 
\end{equation}
where  $\mathbf{ L}$ is an assembled left hand side, $\mathbf{ x}$ is the vector of unknowns and $\mathbf{ r}$ is the assembled right hand side involving known boundary values (for details of derivation see\cite{BeerMarussig}). $ [\mathbf{ B}_{0}]$ is a matrix where the rows refer to the collocation points $n$ and which multiplies with vector $\{ \myVecGreek{\sigma}_{0}\}$ that gathers all initial stress components at grid points inside the inclusions.

\section{Computation of values at grid points inside the inclusion}
To compute the initial stresses at grid points we need to compute the strains at these points. Even though it is possible to compute the strains directly using integral equations this is rather cumbersome as it involves the integration of strongly singular functions and involves complicated methods to isolate the singularity ( see for example \cite{Gao2011} ). To avoid this we compute the displacements first and then compute the strains using a method that is used in the Finite Element Method.

\subsection{Computation of displacements}
The displacement vector $ \mathbf{ u} $ at a grid point $\pt{x}$ inside the inclusion is given by:
\begin{equation}
    \begin{aligned}
        \label{displ}
       \mathbf{ u}(\pt{x}) &= \int_{\boundary} \left[ \fund{U}(\pt{x},\fieldpt) \   \myVec{\dual} (\fieldpt)  - \fund{T}(\pt{x},\fieldpt) \   \myVec{\primary} (\fieldpt) \right] d\boundary (\fieldpt) \\
        &+ \int_{\domain_{0}} \fund{E} (\pt{x},\fieldpt) \myVecGreek{\sigma}_{0} (\fieldpt)  d \domain_{0} (\fieldpt)     
    \end{aligned}
\end{equation}
After the solution the values $ \myVec{\primary} (\fieldpt)$ and $ \myVec{\dual} (\fieldpt)$ are known, so the integrals can be evaluated numerically.

We gather displacement vectors at all grid points in a vector $\{\mathbf{ u} \}$ and obtain:
\begin{equation}
\label{ }
\{\myVec{u}\}= [\hat{\mathbf{ A}}] \mathbf{ x} + \{\bar{\mathbf{ c}}\} + [\bar{\mathbf{ B}}_{0}]\{ \myVecGreek{\sigma}_{0}\} 
\end{equation}
where $[\hat{\mathbf{ A}}]$ is an assembled matrix that multiplies with the unknown $\mathbf{ x}$ and $\{\bar{\mathbf{ c}}\}$ collects the displacement contribution due to given BC's. $ [\bar{\mathbf{ B}}_{0}]$ is similar to $ [\mathbf{ B}_{0}]$ except that the grid point coordinates $\pt{x}_{i}$ replace the source point coordinates $\sourcept_{n}$.

Because of the singularity of $\fund{T}$ the displacements can not be computed on the problem boundary. So if the inclusion point lies on a boundary patch we recover the displacement from the computed boundary values.
For points on a patch boundary ($\pt{x}_{k}$) we replace  Eq. (\ref{displ}) by:
\begin{equation}
\label{ }
\mathbf{ u}(\pt{x}_{k})=\sum_{i}^{I} R_{i}^{u}(\xi_{k},\eta_{k}) \mathbf{ u}_{i}^{e}
\end{equation}
where $R_{i}^{u}(\xi,\eta)$ are the NURBS basis functions used for approximating the displacements in patch $e$, that contains the point $\pt{x}_{k}$ and $\xi_{k},\eta_{k}$ are the local coordinates of the point.
The matrix $[\hat{\mathbf{ A}}]$ and the vector $\{\bar{\mathbf{ c}}\} $ have to be modified for these grid points, whereas $ [\bar{\mathbf{ B}}_{0}]$ will contain zero rows in this case.

\subsection{Computation of strains, general inclusions}

To compute the strains we interpolate the displacements between grid points and obtain for the displacement at a point with the local coordinate $\pt{s}$:
\begin{equation}
\label{Interpolu}
 \myVec{u}(\pt{s})= \sum_{n=1}^{N} M_{n}(\pt{s})  \myVec{u}_{n}
\end{equation}
where $\myVec{u}_{n}$ is the displacement vector at grid point $n$ and $N$ is the number of grid points.  
The interpolation functions $M_{n}(s,r,t)$ are obtained by considering grid coordinates $s,t,r$.

Rewriting (\ref{Interpolu}) for a point with the local coordinates $s,t,r$ in terms of local interpolation functions we have:
\begin{equation}
\label{ }
\mathbf{ u}(s,t,r)= \sum_{i=1}^{I} \sum_{j=1}^{J} \sum_{k=1}^{K} L_{i}(s) L_{j}(t) L_{k}(r) \mathbf{ u}_{n(i,j,k)}
\end{equation}
where $L_{i}(s)$, $L_{j}(t)$, $L_{k}(r)$ are piecewise constant, linear or quadratic interpolation functions of the local coordinates $s,t,r$ respectively and $I,J,K$ specify the span of the function in the local directions $s,t,r$, i.e. 1 for constant, 2 for linear and 3 for quadratic interpolation. $n(i,j,k)$ is the grid node number corresponding to $i,j,k$. The interpolation functions have zero values outside the span.

Replacing the 3 sums by one we have
\begin{equation}
\label{ }
\mathbf{ u}(s,t,r)= \sum_{n=1}^{N} M_{n(i,j,k)}(s,t,r) \mathbf{ u}_{n(i,j,k)}
\end{equation}
where $N$ is the total number of grid points and
\begin{equation}
\label{ }
 M_{n(i,j,k)}(s,t,r)= L_{i}(s) L_{j}(t) L_{k}(r) 
\end{equation}

The derivatives of the displacements are given by
\begin{eqnarray}
\frac{\partial \mathbf{ u}(s,t,r)}{\partial s} & = & \sum_{n=1}^{N}  \frac{\partial M_{n}}{\partial s}\mathbf{ u}_{n}\\
\frac{\partial \mathbf{ u}(s,t,r)}{\partial t} & = & \sum_{n=1}^{N}  \frac{\partial M_{n}}{\partial t}\mathbf{ u}_{n}\\
\frac{\partial \mathbf{ u}(s,t,r)}{\partial r} & = & \sum_{n=1}^{N}  \frac{\partial M_{n}}{\partial r}\mathbf{ u}_{n}\\
\end{eqnarray}
where
\begin{eqnarray}
\frac{\partial M_{n}}{\partial s} & = &  \frac{\partial L_{i}(s)}{\partial s} L_{j}(t) L_{k}(r) \\
\frac{\partial M_{n}}{\partial t} & = &   L_{i}(s) \frac{\partial L_{j}(t)}{\partial t} L_{k}(r)\\
\frac{\partial M_{n}}{\partial r} & = &   L_{i}(s) L_{j}(t) \frac{\partial L_{k}(r)}{\partial r} 
\end{eqnarray}

Unfortunately we can not use NURBS for the interpolation functions as they are based on parameter values instead of real values. Therefore they can not be used to interpolate the real displacement values at internal points. We use Lagrange polynomials instead.

The strains are given by:
\begin{eqnarray}
\epsilon_{x} & = & \frac{\partial u_{x}}{\partial x} = \sum_{n=1}^{N} \frac{\partial M_{n}}{\partial x} u_{xn}\\
\epsilon_{y} & = & \frac{\partial u_{y}}{\partial y} = \sum_{n=1}^{N} \frac{\partial M_{n}}{\partial y} u_{yn}\\
\epsilon_{z} & = & \frac{\partial u_{z}}{\partial z} = \sum_{n=1}^{N} \frac{\partial M_{n}}{\partial z} u_{zn}\\
\gamma_{xy} & = & \frac{\partial u_{x}}{\partial y} + \frac{\partial u_{y}}{\partial x}= \sum_{n=1}^{N} \frac{\partial M_{n}}{\partial x} u_{yn} + \sum_{n=1}^{N} \frac{\partial M_{n}}{\partial y} u_{xn}\\
\gamma_{zy} & = & \frac{\partial u_{z}}{\partial y} + \frac{\partial u_{y}}{\partial z}= \sum_{n=1}^{N} \frac{\partial M_{n}}{\partial z} u_{yn} + \sum_{n=1}^{N} \frac{\partial M_{n}}{\partial y} u_{zn}\\
\gamma_{xz} & = & \frac{\partial u_{x}}{\partial z} + \frac{\partial u_{z}}{\partial x}= \sum_{n=1}^{N} \frac{\partial M_{n}}{\partial x} u_{zn} + \sum_{n=1}^{N} \frac{\partial M_{n}}{\partial z} u_{xn}
\end{eqnarray}
The strains at grid point $k$ can be written in matrix notation:
\begin{equation}
\label{Bhat}
\myVecGreek{\epsilon}(\pt{x}_{k})= \hat{\mathbf{ B}}(\pt{x}_{k}) \{\myVec{u}\}
\end{equation}
where
\begin{equation}
\label{ }
\hat{\mathbf{ B}}(\pt{x}_{k})= \left(\begin{array}{ccc}\mathbf{ B}_{1} & \mathbf{ B}_{2} & \cdots\end{array}\right)
\end{equation}
and
\begin{equation}
\label{Bsub}
 \mathbf{ B}_{i}=\left(\begin{array}{ccc}\frac{\partial M_{i}}{\partial x} & 0 & 0 \\0 & \frac{\partial M_{i}}{\partial y} & 0 \\0 & 0 & \frac{\partial M_{i}}{\partial z}  \\ \frac{\partial M_{i}}{\partial y} & \frac{\partial M_{i}}{\partial x} & 0 \\ 0 & \frac{\partial M_{i}}{\partial z} & \frac{\partial M_{i}}{\partial y}  \\ \frac{\partial M_{i}}{\partial z} & 0 &\frac{\partial M_{i}}{\partial x} \end{array}\right)
\end{equation}

The global derivatives of $M_{n}$ are given by:
\begin{equation}
\label{ }
\left(\begin{array}{c}\frac{\partial M_{n}}{\partial x} \\ \\\frac{\partial M_{n}}{\partial y} \\ \\\frac{\partial M_{n}}{\partial z}\end{array}\right)= \mathbf{ J}^{-1} \left(\begin{array}{c}\frac{\partial M_{n}}{\partial s} \\ \\\frac{\partial M_{n}}{\partial t} \\ \\ \frac{\partial M_{n}}{\partial r}\end{array}\right)
\end{equation}
where $\mathbf{ J}$ is the Jacobian matrix Eq. (\ref{eq:jac3D}).  
For a linear inclusion we compute the strain in local directions as is shown later.

Gathering all strain vectors at grid points in $\{\myVecGreek{\epsilon}\}$ we can write:
\begin{equation}
\label{ }
\{\myVecGreek{\epsilon}\}= [\hat{\mathbf{ B}}] \{\myVec{u}\}
\end{equation}

After substitution of $\{\myVec{u}\}$ we obtain:
\begin{equation}
\label{Strain}
\{\myVecGreek{\epsilon}\}= [\hat{\mathbf{ B}}] \left( [\hat{\mathbf{ A}}] \mathbf{ x} + \{\bar{\mathbf{ c}}\} + [\bar{\mathbf{ B}}_{0}]\{ \myVecGreek{\sigma}_{0}\} \right)
\end{equation}

The initial stresses are computed by 
\begin{equation}
\label{ }
\{\myVecGreek{\sigma}_{0} \}=\left[ \mathbf{ D} - \mathbf{ D}_{incl}\right] \{ \myVecGreek{\epsilon} \}= \left[ \mathbf{ D} - \mathbf{ D}_{incl}\right][\hat{\mathbf{ B}}] \left( [\hat{\mathbf{ A}}] \mathbf{ x} + \{\bar{\mathbf{ c}}\} + [\bar{\mathbf{ B}}_{0}]\{ \myVecGreek{\sigma}_{0}\} \right)
\end{equation}
where $\left[ \mathbf{ D} - \mathbf{ D}_{incl}\right] $ is a matrix containing $\mathbf{ D} - \mathbf{ D}_{incl}$ as sub-matrices on the diagonal.

\subsubsection{Computation of strain for linear inclusions}

For linear inclusions it is convenient to work with the strain in local coordinates. If we assume the bolt to be fully grouted, i.e. no slip is allowed between the bolt and the domain it is embedded in and that the Poisson's ratio of the bolt has no effect, the only strain that has to be considered is the one along the bar\footnote{It should be noted that this restriction can be lifted, i.e. slip can be considered.}:
\begin{eqnarray}
\epsilon_{z^{\prime}} & = & \frac{\partial u_{z^{\prime}}}{\partial z^{\prime}} = \sum_{n=1}^{N} \frac{\partial M_{n}}{\partial z^{\prime}} u_{z^{\prime}n}= \sum_{n=1}^{N} \frac{\partial M_{n}}{\partial s} \frac{1}{J} ( \mathbf{ v}_{z^{\prime}}\cdot \mathbf{ u}_{n})
\end{eqnarray}
where $J$ is the Jacobian and $ \mathbf{ v}_{z^{\prime}} $ is a unit vector in $z^{\prime}$ direction.

 Eq. (\ref{Bsub}) now becomes
\begin{equation}
\label{ }
\mathbf{ B}_{n}=\frac{1}{J} \left(\begin{array}{ccc}0 & 0 & 0 \\0 & 0 & 0 \\\frac{\partial M_{n}}{\partial s}  v_{z^{\prime}_x} & \frac{\partial M_{n}}{\partial s}  v_{z^{\prime}_y}  & \frac{\partial M_{n}}{\partial s}  v_{z^{\prime}_z}  \\0 & 0 & 0 \\0 & 0 & 0 \\0 & 0 & 0\end{array}\right)
\end{equation}

The local initial stress vector is given by:
\begin{equation}
\label{ }
 \{ \myVecGreek{\sigma^{\prime}}_{0} \}= (\mathbf{ D}^{\prime} - \mathbf{ D}_{incl}^{\prime} )\{\myVecGreek{\epsilon}^{\prime} \}
\end{equation}
where
\begin{equation}
\label{ }
(\mathbf{ D}^{\prime} - \mathbf{ D}_{incl}^{\prime} )= \left(\begin{array}{cccccc}0 & 0 & 0 & 0 & 0 & 0 \\0 & 0 & 0 & 0 & 0 & 0 \\0 & 0 & E - E_{incl}  & 0 & 0 & 0 \\0 & 0 & 0 & 0 & 0 & 0 \\0 & 0 & 0 & 0 & 0 & 0 \\0 & 0 & 0 & 0 & 0 & 0\end{array}\right)
\end{equation}
where $E$ and $E_{incl}$ is the Young' modulus of the domain and the inclusion respectively and
\begin{equation}
\label{ }
 \{ \myVecGreek{\sigma^{\prime}}_{0} \}= \left\{\begin{array}{c}0 \\0 \\ \sigma_{0z^{\prime}} \\0 \\0 \\ 0 \end{array}\right\}
\end{equation}

\section{Solution procedure}
A solution that already includes the effect of inclusions that have different elastic properties is possible by combining equation (\ref{DisIE}) with (\ref{Strain}).
Eq. (\ref{Strain}) can be written in the following form:
\begin{equation}
\label{Strain_rev1}
\{\myVecGreek{\epsilon}\}= [\hat{\mathbf{ C}}]  \{ \mathbf{ x} \} + \{\bar{\bar{\mathbf{ c}}}\} + [\bar{\mathbf{ C}}_{0}]  ( [\mathbf{ D}] - [\mathbf{ D}_{incl}]) \{ \myVecGreek{\epsilon} \}
\end{equation}
\noindent where:
\begin{equation}
\label{}
[\hat{\mathbf{ C}}] = [\hat{\mathbf{ B}}] [\hat{\mathbf{ A}}] \hspace{10mm}
[\hat{\mathbf{ C}}_{0}] = [\hat{\mathbf{ B}}] [\bar{\mathbf{B}}_{0}] \hspace{10mm}
\{\bar{\bar{\mathbf{ c}}}\} = [\hat{\mathbf{ B}}] \{\bar{\mathbf{ c}}\}
\end{equation}
Eq. (\ref{Strain_rev1}) along with Eq. (\ref{DisIE}) form the following linear system of equations:
\begin{equation}
\label{Final_sys}
\begin{pmatrix}
  [\mathbf{ L}]        & - [\mathbf{ B}_{0}] ( [\mathbf{ D}] - [\mathbf{ D}_{incl}])  \\ \\
 - [\hat{\mathbf{ C}}] & [{\mathbf{I}}] - [\hat{\mathbf{ C}}_{0}] ( [\mathbf{ D}] - [\mathbf{ D}_{incl}])
\end{pmatrix}
\begin{pmatrix}
\{\mathbf{ x}\} \\ \\
\{\myVecGreek{\epsilon}\}
\end{pmatrix}
=
\begin{pmatrix}
\{\mathbf{ r}\} \\ \\
\{\bar{\bar{\mathbf{ c}}}\}
\end{pmatrix}
\end{equation}
that can be solved in terms of boundary unknowns and internal strains.

It is also possible to obtain a system of equations that only multiplies with the boundary unknown:
\begin{equation}
\label{onestep}
[\mathbf{ L}]^{\prime} \{\mathbf{ x}\} = \{\mathbf{ r}\}^{\prime} 
\end{equation}
where $[\mathbf{ L}]^{\prime}$ and $\{\mathbf{ r}\}^{\prime} $ are modified left and right hand sides.

We rewrite the strain vector as: 
\begin{equation}
\label{Strain_rev2}
\{\myVecGreek{\epsilon}\}= ([\mathbf{I}] - 
[\hat{\mathbf{ C}}_0] ([\mathbf{D}]-[\mathbf{D}_{incl}])^{-1}) ([\hat{\mathbf{C}}] \{\mathbf{ x}\} + \{\bar{\bar{\mathbf{ c}}}\}) = [\mathbf{A}] \{\mathbf{x}\} + \{\mathbf{b}\}
\end{equation}

\noindent where

\begin{equation}
[\mathbf{A}] = ( [\mathbf{I}] - [\hat{\mathbf{C}}_0] ([\mathbf{D}] - [\mathbf{D}_{incl}]))^{-1} [\hat{\mathbf{C}}] \hspace{5mm} \{\mathbf{b}\} = ( [\mathbf{I}] - [\hat{\mathbf{C}}_0] 
([\mathbf{D}] - [\mathbf{D}_{incl}]))^{-1} \{\bar{\bar{\mathbf{ c}}}\}
\end{equation}

Eq. (\ref{Strain_rev2}) can be inserted in Eq. (\ref{DisIE}) in order to obtain:
\begin{equation}
[\mathbf{L}] \{\mathbf{x}\} = \{\mathbf{r}\} + [\mathbf{B}_0] 
([\mathbf{D}] - [\mathbf{D}_{incl}]) ([\mathbf{A}] \{\mathbf{x}\} + \{\mathbf{b}\})
\end{equation}
and, hence, the following system of equations can be obtained
\begin{equation}
( [\mathbf{L}] - [\mathbf{B}_0]) ([\mathbf{D}] - [\mathbf{D}_{incl}]) [\mathbf{A}] ) 
\{\mathbf{x}\} = \{\mathbf{r}\} + [\mathbf{B}_0] 
([\mathbf{D}] - [\mathbf{D}_{incl}]) \{\mathbf{b}\}
\end{equation}

The matrices in Eq. (\ref{onestep}) are defined by:
\begin{eqnarray}
[\mathbf{L}]^{\prime} & = & ( [\mathbf{L}] - [\mathbf{B}_0]) ([\mathbf{D}] - [\mathbf{D}_{incl}]) [\mathbf{A}] )   \\
\{\mathbf{ r}\}^{\prime} & = &\{\mathbf{r}\} + [\mathbf{B}_0] 
([\mathbf{D}] - [\mathbf{D}_{incl}]) \{\mathbf{b}\} 
\end{eqnarray}

\subsection{Elasto-plasticity, Newton-Raphson method}
If the stress at an inclusion point exceeds the elastic limit incremental/iterative  elasto-plastic procedures, well known in the FEM,  are applied. A detailed discussion of these methods is beyond the scope of this paper. A very good description can be found in \cite{Smith}.

First we consider that $[\mathbf{D}_{incl}]=[\mathbf{D}_{e,incl}]$ for the case of elastic behaviour and $[\mathbf{D}_{incl}]=[\mathbf{D}_{ep,incl}]$ for the case of elasto-plastic behaviour, where $[\mathbf{D}_{e,incl}]$ is the elastic constitutive matrix and  $[\mathbf{D}_{ep,incl}]$ is the elasto-plastic constitutive matrix. 
We proceed in an incremental/iterative way and check after each increment if the yield function F($\myVecGreek{\sigma}$) is smaller or greater than zero.

The increment of stress in the plastic regime is:
\begin{equation}
\label{ }
\triangle \myVecGreek{\sigma}= \mathbf{ D}_{ep,incl} \triangle \myVecGreek{\epsilon}^{*} 
\end{equation} 
where $\triangle \myVecGreek{\epsilon}^{*}$ is the in-elastic strain increment (i.e. the one that occurs after the stress has reached a state where F($\myVecGreek{\sigma}$)=0) .
The in-elastic ($\triangle \myVecGreek{\epsilon}^{*}$) and elastic ($\triangle \myVecGreek{\epsilon}^{e}$) strain increments can be computed by:
\begin{equation}
\label{Plastr}
\triangle \myVecGreek{\epsilon}^{*}= f \triangle \myVecGreek{\epsilon} \quad , \quad \triangle \myVecGreek{\epsilon}^{e}= (1-f) \triangle \myVecGreek{\epsilon}
\end{equation}
where $\triangle \myVecGreek{\epsilon}$ is the total plastic strain increment and 
\begin{equation}
\label{ }
f=\frac{F_{new}}{F_{new}-F_{old}} \quad \text{ if $f$ > 0}  \quad \text{ $f$=0  otherwise}
\end{equation}
 $F_{new}$ is the value of $F$ at the end of the increment, $F_{old}$ is the value at the beginning.
 
The initial stress due to plasticity is given by:
\begin{equation}
\label{Plasigp}
\triangle \myVecGreek{\sigma}_{0}^{p}= ( \mathbf{ D} - \mathbf{ D}_{ep,incl}) \triangle \myVecGreek{\epsilon}^{*}=  ( \mathbf{ D} - \mathbf{ D}_{ep,incl}) f \triangle \myVecGreek{\epsilon}
\end{equation}

If the increment in total strain $\triangle \myVecGreek{\epsilon}$ has occurred while traversing the yield surface then:
\begin{equation}
\label{InitialS1}
\myVecGreek{\sigma}_{0}^{e} =( \mathbf{ D} - \mathbf{ D}_{e,incl}) (1-f) \triangle \myVecGreek{\epsilon} 
\end{equation}

The total initial stress increment, including plastic effects, is given by
\begin{equation}
\label{Plasig}
\triangle \myVecGreek{\sigma}_{0}= \triangle \myVecGreek{\sigma}_{0}^{e} + \triangle \myVecGreek{\sigma}_{0}^{p}= ( (1-f)(\mathbf{ D} - \mathbf{ D}_{e,incl}) + f(\mathbf{ D} - \mathbf{ D}_{ep,incl})) \triangle \myVecGreek{\epsilon} = \mathbf{ D}^{\prime} \triangle \myVecGreek{\epsilon} 
\end{equation}

where 
\begin{equation}
\label{Dprime}
 \mathbf{ D}^{\prime}= ( (1-f)(\mathbf{ D} - \mathbf{ D}_{e,incl}) + f(\mathbf{ D} - \mathbf{ D}_{ep,incl}))
\end{equation}

If the strain increment is totally elastic (f=0) we have:
\begin{equation}
\label{}
\triangle \myVecGreek{\sigma}_{0}=  \mathbf{ D}^{\prime} \triangle \myVecGreek{\epsilon} = (\mathbf{ D} - \mathbf{ D}_{e,incl}) \triangle \myVecGreek{\epsilon} 
\end{equation}

If the strain increment is totally plastic (f=1) we have:
\begin{equation}
\label{}
\triangle \myVecGreek{\sigma}_{0}= \mathbf{ D}^{\prime} \triangle \myVecGreek{\epsilon} = (\mathbf{ D} - \mathbf{ D}_{ep,incl})\triangle \myVecGreek{\epsilon} 
\end{equation}

To start the simulation we set $f$ in Eq. (\ref{Dprime}) equal to zero i.e. $\mathbf{ D}^{\prime}= \mathbf{ D}_{e,incl}$ and obtain the first result:
\begin{equation}
\label{}
[\mathbf{ L}]^{\prime} \{ \mathbf{ x}_{0}\} = \{\mathbf{ r}\}^{\prime} 
\end{equation}

With this result we compute the value of the yield function at internal points and compute a load factor $\lambda$ that reduces the load to the one where first yield occurred. We then reduce the results to first yield:
\begin{equation}
\label{ }
\{ \mathbf{ x}_{0}\}=> \lambda \{ \mathbf{ x}_{0}\}
\end{equation}
We also adjust the stresses at internal points to the new load level.
For the first load step ($n$=1) the matrix $\mathbf{ D}^{\prime}$ is updated with $f=1$ for the point where the first yield has occurred. 

We  apply the rest of the loading in $n$ steps. The increment in load is:
\begin{equation}
\label{ }
\triangle \{\mathbf{ r}\}^{\prime} = \frac{1 -\lambda}{n}  \{\mathbf{ r}\}^{\prime} 
\end{equation}
At each increment we solve:
\begin{equation}
\label{}
[\mathbf{ L}]^{\prime} \{\triangle \mathbf{ x}_{i}\} = \{\triangle \mathbf{ r}\}^{\prime} 
\end{equation}
For the second and subsequent load steps we update the matrix $\mathbf{ D}^{\prime}$ with $f$ according to the current state of stress.
Standard return alogrithms can be applied to ensure that the stresses stay on the yield surface.

\begin{remark}
Note that incremental/iterative procedures need only be applied for elasto-plastic behaviour. When inclusions are defined that have different elastic behaviour the solution is obtained without iteration. 
 \end{remark}

\section{Test Examples}
\begin{figure}
\begin{center}
\begin{overpic}[scale=0.6]{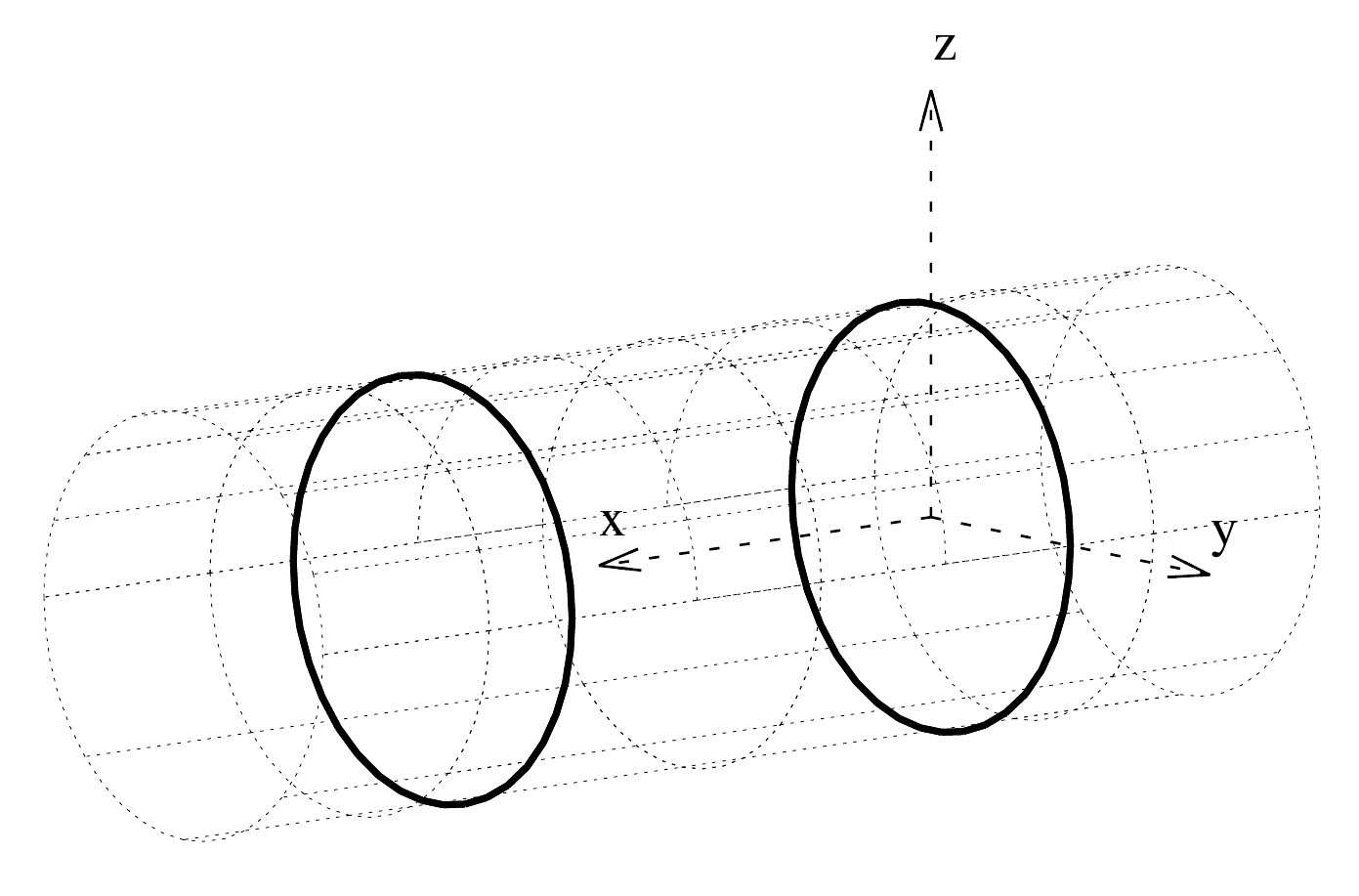}
\end{overpic}
\caption{Geometry of test example.}
\label{Exgeo}
\end{center}
\end{figure}

We test the implementation on an example of a circular, infinitely long, tunnel in an infinite domain subjected to a virgin stress (see \myfigref{Exgeo}).
Since the aim of the examples is to test the accuracy of the implementation and not to perform a realistic simulation, we use non-dimensional parameters. For the domain we assume E=1 and $\nu=0$ and for the radius of the tunnel, R=1.
The results of the IGABEM simulation are compared with a FEM analysis using the software PLAXIS and an analytical calculation.

\subsection{Elastic analysis with no ground support}
\label{ElastAnalysis}
In order to establish a comparable regime between the IGABEM and PLAXIS simulations we first conduct an elastic simulation without ground support and with a virgin stress $\sigma_{x}=0,\sigma_{y}=0,\sigma_{z}=-1$.

\subsubsection{Discretisation with PLAXIS}
\begin{figure}
	\begin{center}
		\begin{overpic}[scale=0.8]{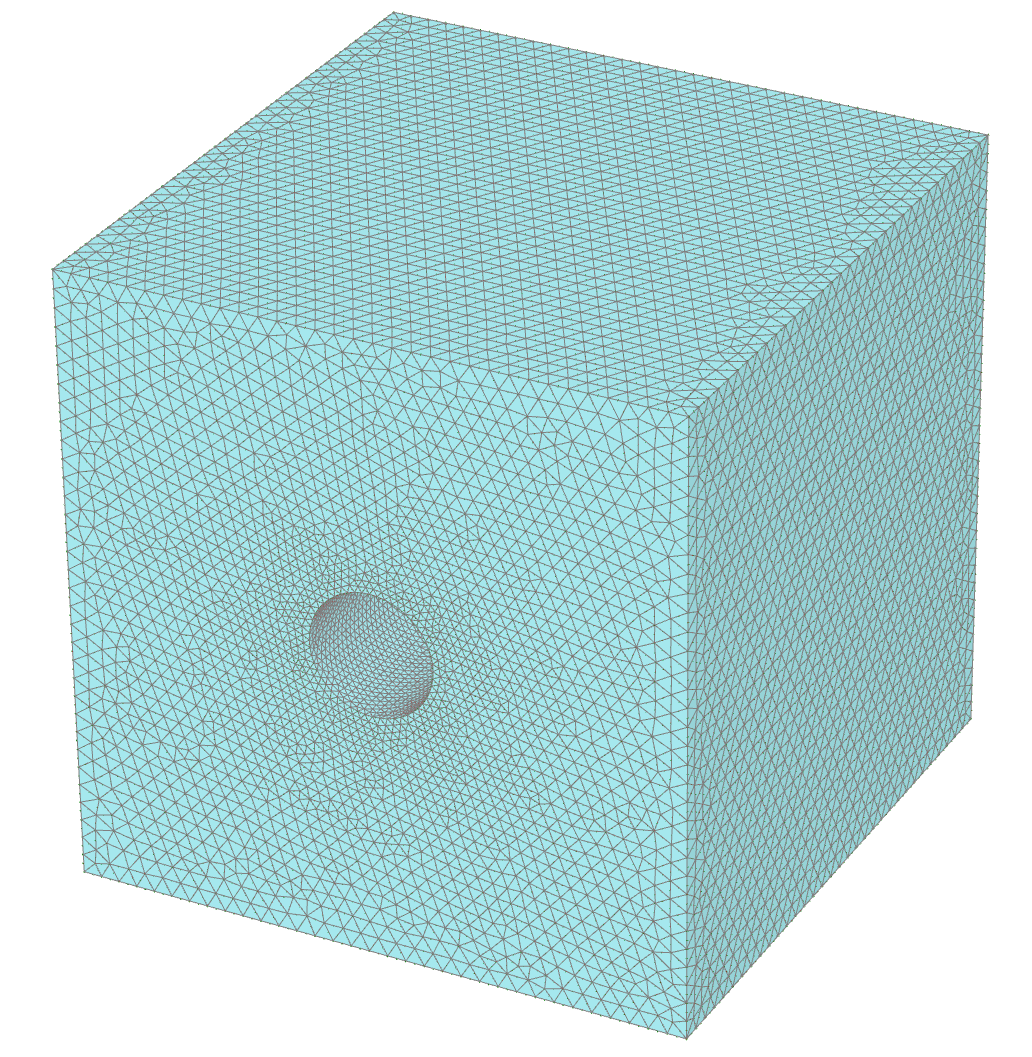}
		\end{overpic}
		\caption{Finite element mesh with distance 2.5D from the tunnel center to the outer boundary.}
		\label{FEM_Mesh}
	\end{center}
\end{figure}

The infinite domain is approximated by providing an artificial boundary at some distance from the tunnel. Since this distance affects the results we examine the error introduced by varying the distances from the tunnel centre to the artificial boundary to 2.5, 5 and 10 times the tunnel diameter (2.5D, 5D, 10D).
The boundary conditions at the outer boundary of the mesh (shown in \myfigref{FEM_Mesh}) are that displacements normal to the faces of the cuboid are set to zero. 
The fact that the tunnel is infinitely long can be modelled by extending the mesh along the tunnel axis to 10m (=5D). 
The mesh consists of quadratic 10-noded tetrahedral elements. For load case 0 the virgin stresses are assigned to all elements of the mesh and the excavation is then simulated by deactivating the elements inside the tunnel. This will give the same result as a much simpler plane strain simulation but the same mesh will be used for the second test example where a 3-D analysis is necessary.

\subsubsection{Discretisation with IGABEM}
Only the boundary of the tunnel is discretised as the infinite domain is explicitly considered by the fundamental solutions. The excavation boundary is defined by 16 control points and basis functions of order 2 (quadratic) along the tunnel walls and of order 1 (linear) along the tunnel axis are used. It should be noted that this geometrical description exactly represents a circular tunnel. To simulate the infinite extent of the tunnel we use plane strain infinite patches. 
\begin{figure}
\begin{center}
\begin{overpic}[scale=0.6]{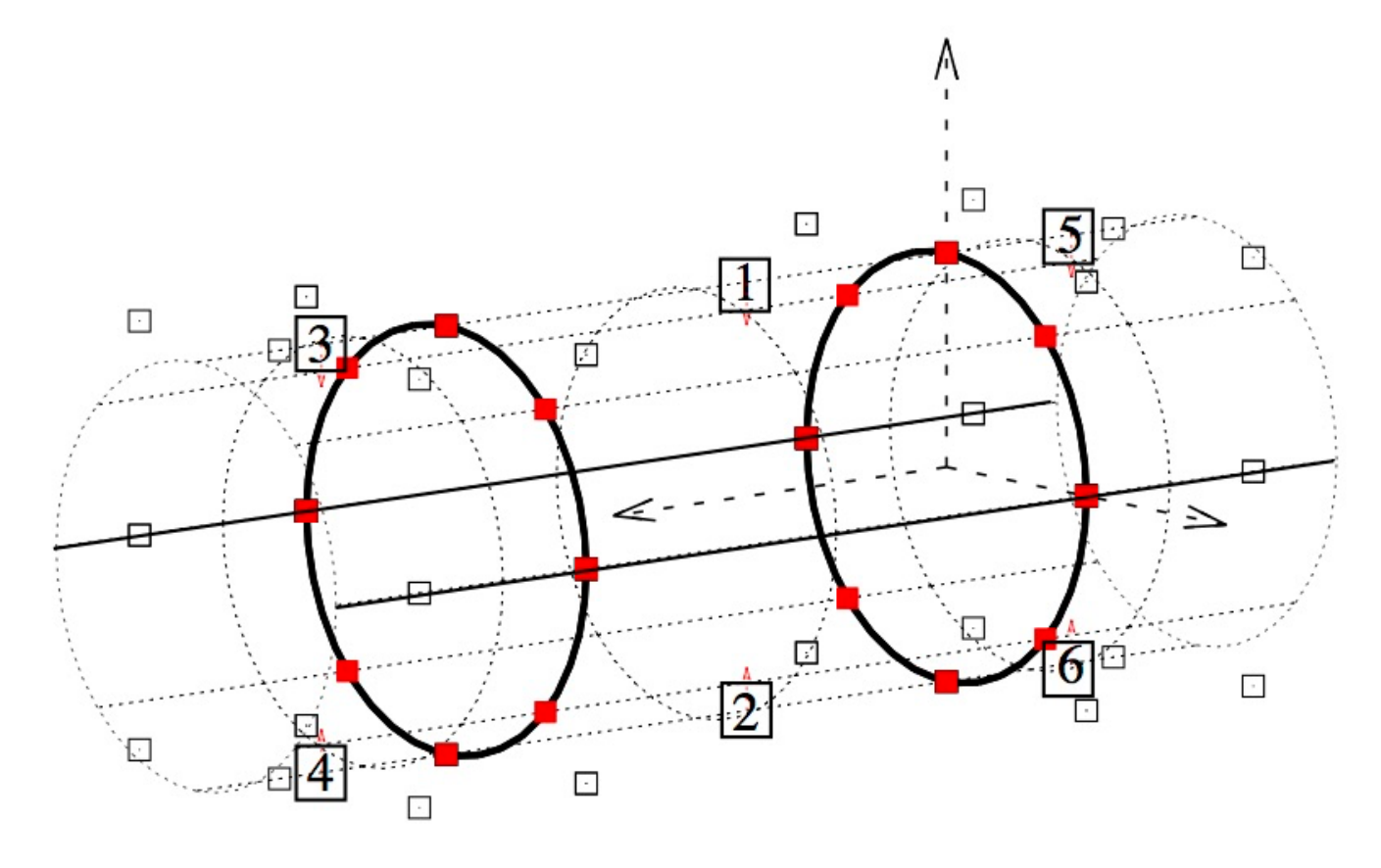}
\put(44,27){$x$}
\put(90,25){$y$}
\put(70,60){$z$}
\end{overpic}
\caption{Discretisation of tunnel into 6 patches. Patch numbers are shown. Patches 3 to 6 are infinite. Control points are shown as hollow squares, collocation points as filled squares}
\label{ExMesh}
\end{center}
\end{figure}
The excavation of the tunnel is simulated by assigning a virgin stress state. Excavation forces on the boundary are then automatically computed.
The discretisation is shown in \myfigref{ExMesh}. For the approximation of the displacements the same basis functions as for the description of the geometry are used, resulting in the collocation points shown. The discretisation has 48 degrees of freedom.

\subsubsection{Comparison of results}
To investigate the effect of the artificial boundary on the results of the FEM analysis we compare the values of vertical displacements along a vertical line above the tunnel with the exact result (Kirsch solution \cite{Kirsch}) and the IGABEM result.
\begin{figure}
\begin{center}
	\begin{overpic}[scale=0.6]{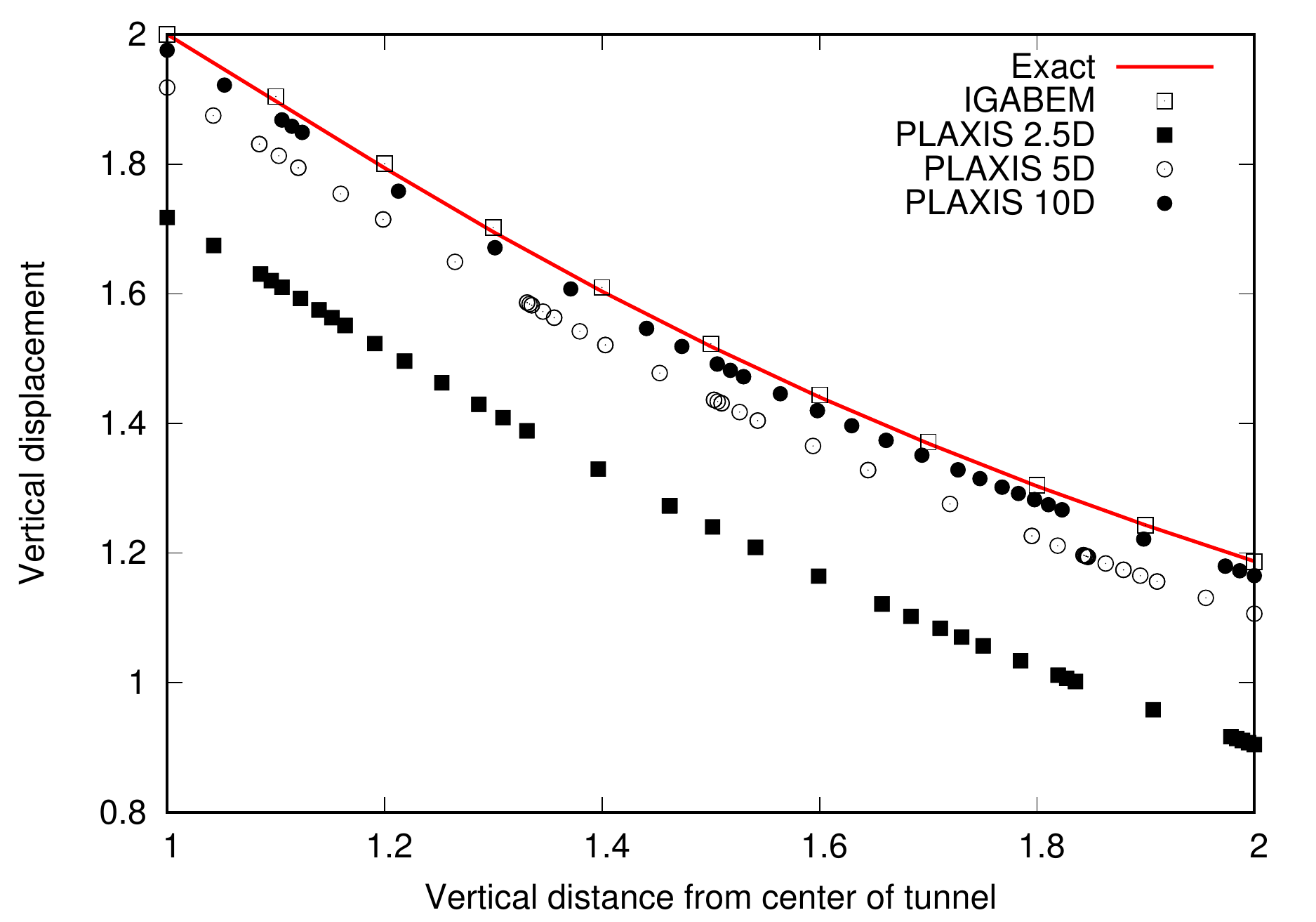}
\end{overpic}
\caption{Variation of the vertical displacement above the tunnel: Comparison of FEM results with the BEM result and the exact solution.}
\label{Uzconv}
\end{center}
\end{figure}
It can be seen in \myfigref{Uzconv} that the FEM results converge to the exact results, while the IGABEM result already is in excellent agreement. The computing time for the PLAXIS simulation was 3 minutes 53 seconds and for the IGABEM simulation 8 seconds. 

\subsection{Elastic simulation with rock bolts}
\begin{figure}
\begin{center}
\begin{overpic}[scale=0.5]{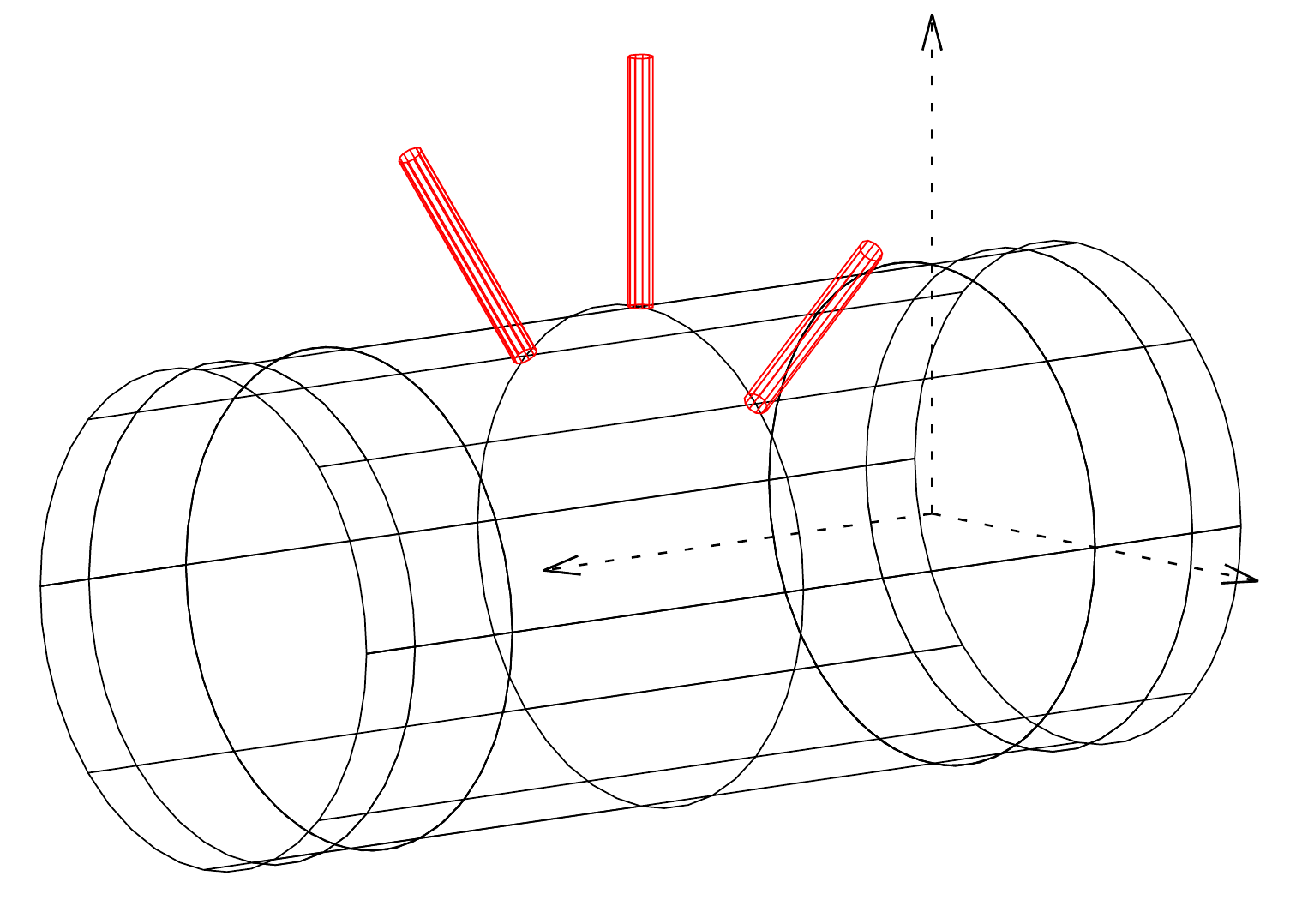}
\put(38,27){$x$}
\put(98,25){$y$}
\put(70,73){$z$}
\end{overpic}
\caption{Geometry of the elastic simulation with bolts}
\label{ExBoltsgeo}
\end{center}
\end{figure}

This example is designed to test the implementation of rock bolts. Three rock bolts are installed at the top of the tunnel (\myfigref{ExBoltsgeo}). The diameter of the bolts is 0.05 and the elastic modulus is twice the value of the domain (i.e. $\mathrm{E}_{bolt}=2$). It is noted that the consideration of the rock bolts renders the analysis three dimensional. The aim is to study the local influence of the rock bolts, assuming that further away plane strain conditions prevail.

\subsubsection{Discretisation with PLAXIS}
\begin{figure}
	\begin{center}
		\begin{overpic}[scale=0.25]{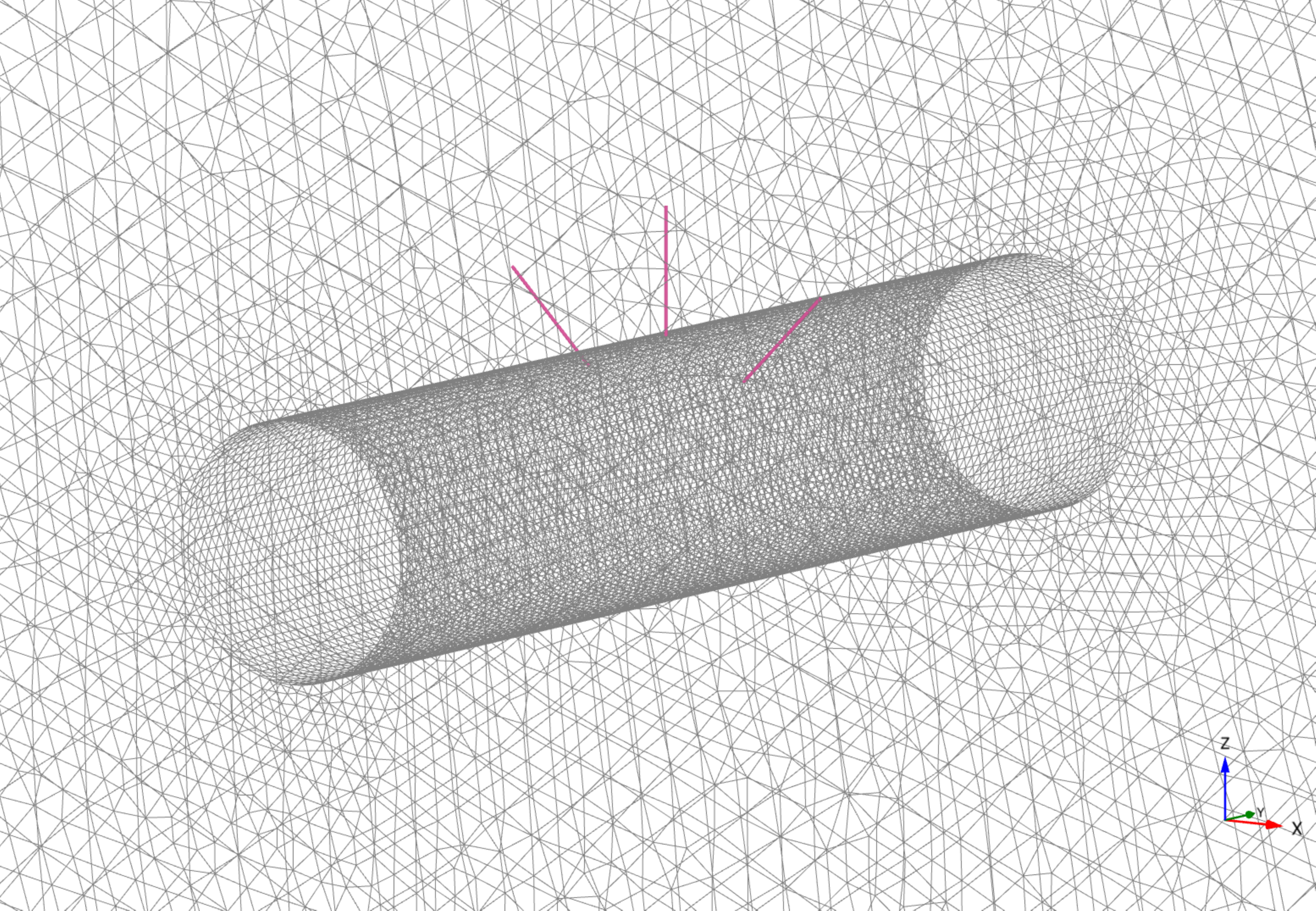}
		\end{overpic}
		\caption{Bolts simulated as embedded beams with PLAXIS}
		\label{FEMBeam}
	\end{center}
\end{figure}

In PLAXIS rock bolts are simulated as embedded beams. The geometry of the bolt can be arbitrarily located in the finite element mesh, thus nodes of the bolt do not need to coincide with nodes of the finite element mesh. As the beam cross section is small, compared with its length, the bending stiffness is almost zero, therefore it can be assumed that the beam acts as a bolt. In \myfigref{FEMBeam} the three bolts are shown embedded in the surrounding FEM mesh. 

\subsubsection{Discretsation with IGABEM}
\begin{figure}
\begin{center}
\begin{overpic}[scale=0.5]{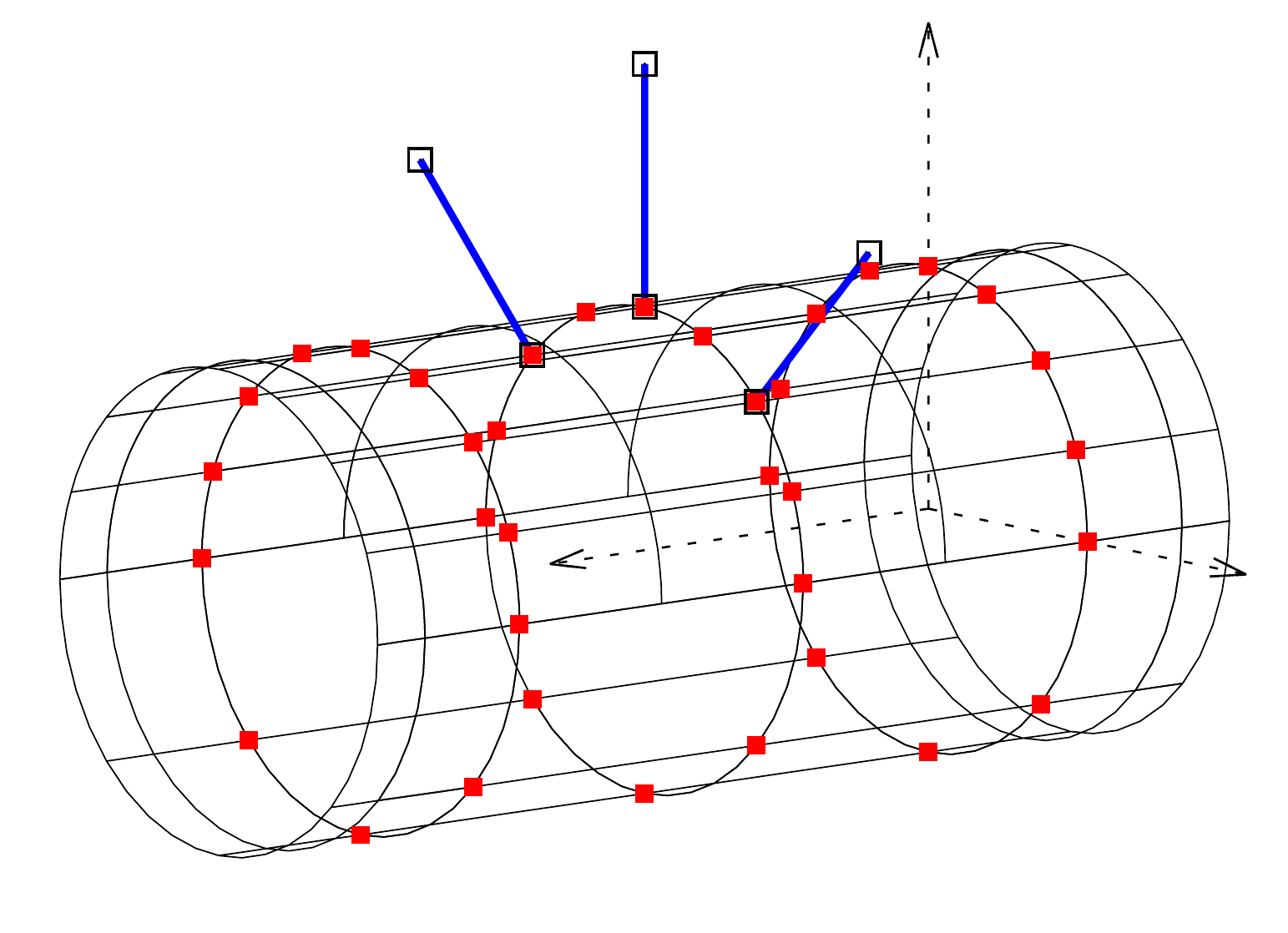}
\put(38,27){$x$}
\put(98,25){$y$}
\put(70,73){$z$}
\end{overpic}
\caption{Discretisation of the problem with IGABEM. Shown are the bolts in blue with the associated control points as hollow squares. The collocation points obtained after knot insertion are shown as red filled squares.}
\label{ExBoltsBmesh}
\end{center}
\end{figure}

The bolts are defined as linear inclusions.
In order to get results comparable to PLAXIS the continuity of displacements is changed to $C^0$ at the point where the rock bolts meet the boundary patch. This is done be inserting knots into the basis functions that approximate the unknown. This increases the number of collocation points and the degrees of freedom.
The IGABEM discretisation is shown in \myfigref{ExBoltsBmesh} and has 108 degrees of freedom.

\subsubsection{Comparison of results}
\begin{figure}
\begin{center}
\begin{overpic}[scale=0.5]{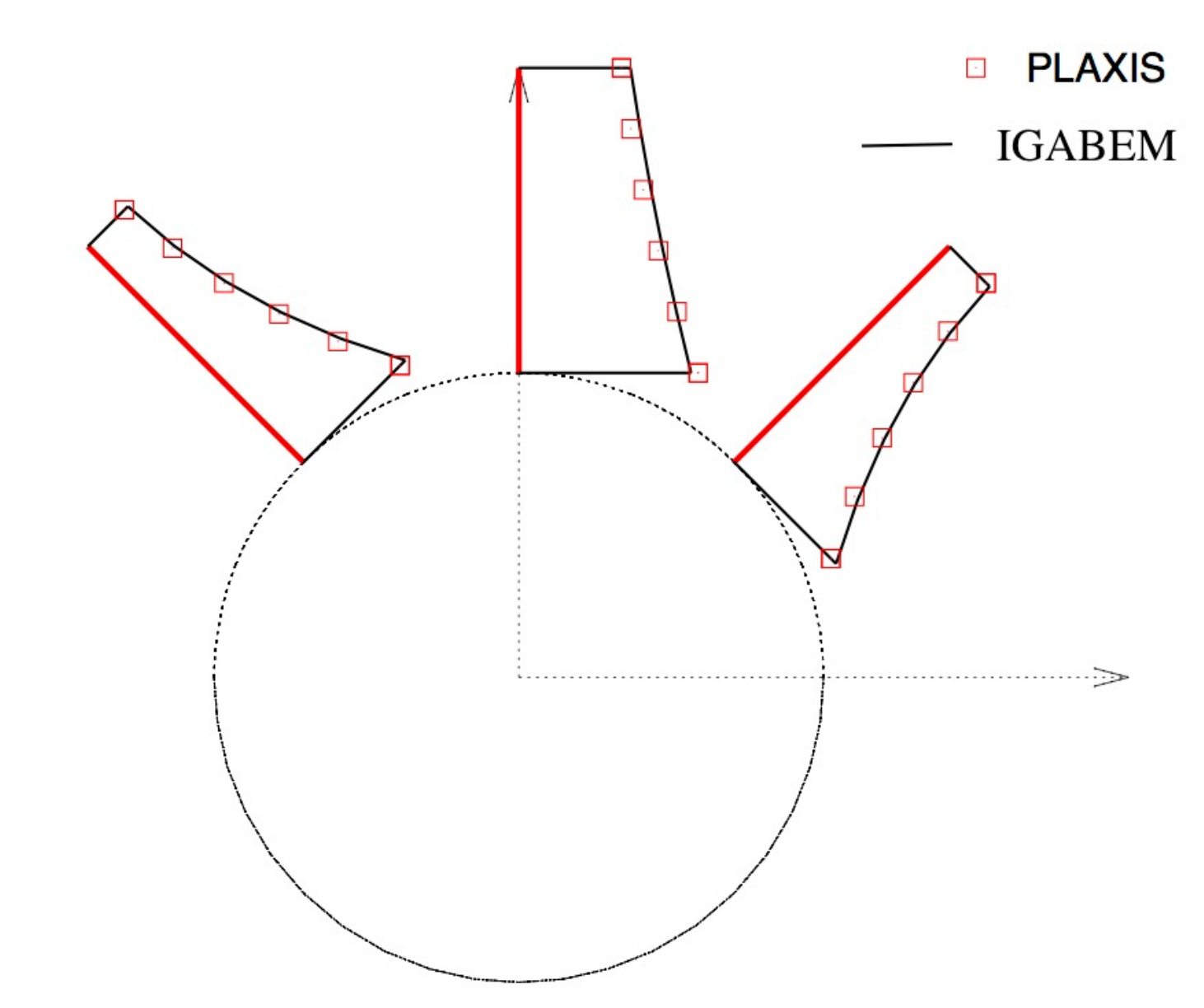}
\end{overpic}
\caption{Comparison of displacements along the rock bolts.}
\label{ExBoltsResult}
\end{center}
\end{figure}

We compare the displacements along the rock bolts in \myfigref{ExBoltsResult}. Good agreement can be observed.
The computation time of the PLAXIS simulation was 6 minutes and 12 seconds and for the IGABEM simulation 24 seconds.

\newpage

\subsection{Elasto-plastic simulation without rock bolts}
\begin{figure}
\begin{center}
\begin{overpic}[scale=0.6]{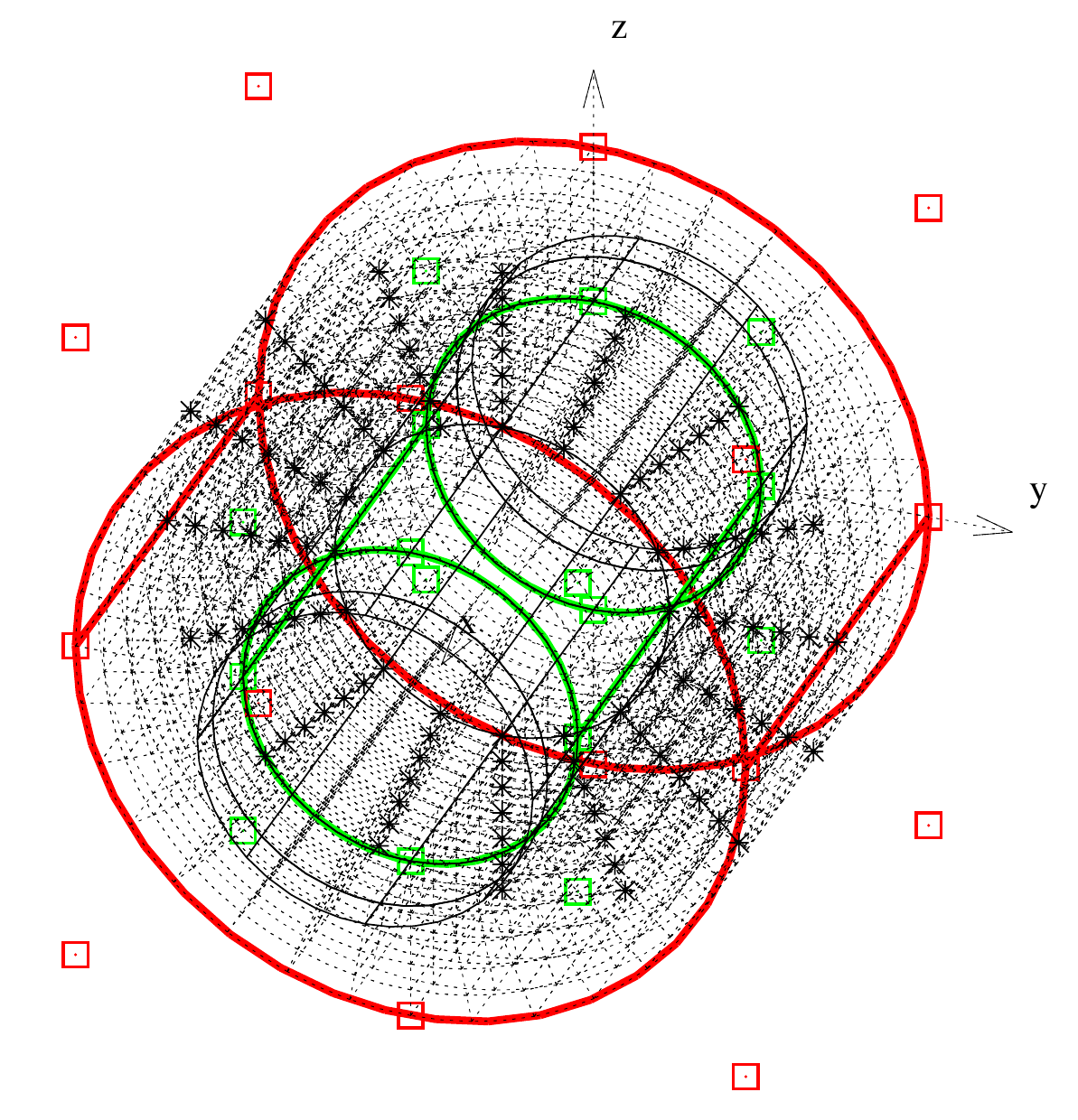}
\end{overpic}
\caption{Volume discretisation for the plasticity example. Shown are the two (green and red) surfaces defining the inclusion and the control points as hollow squares. Also shown are the internal points as stars and the subdivision into integration regions for the volume integration. These regions are automatically determined depending on the location of internal points and on the aspect ratio.}
\label{ExPlasmesh}
\end{center}
\end{figure}

For this test we can compare with an analytical solution in plane strain due to Duncan-Fama (see for example \cite{Hoek}). 
For a hydrostatic virgin stress of $p_0$  (compression positive) and a Mohr-Coulomb yield condition with a friction angle $\phi$ and cohesion c, the solution for the elasto-plastic radial displacement of the tunnel wall is given by:
\begin{equation}
\label{ }
u_{p}= \frac{R (1+ \nu)}{E} \left( 2(1-\nu)(p_{0} - p_{cr})\frac{r_{p}}{R}^2 - (1-2\nu)p_{0} \right)
\end{equation}
where the extent of the plastic zone is given by:
\begin{equation}
\label{Pext}
r_{p}=R \left( \frac{2p_{0}(k-1) + \sigma_{cm}}{(1+k) \sigma_{cm}} \right)^\frac{1}{k-1}
\end{equation}
and 
\begin{eqnarray}
 \sigma_{cm} & = & \frac{2\mathrm{c} \ cos\phi}{1-sin \phi} \\
k & = & \frac{1+sin\phi}{1-sin\phi} 
\end{eqnarray}

\subsubsection{Discretisation with IGABEM}
For this simulation the excavation is surrounded with a general inclusion extending one radius from the excavation surface as shown in \myfigref{ExPlasmesh}.
The properties assigned to the inclusion are $E=1, \ \nu=0, \ \phi=10^{\circ}, \  c=0.5$. A hydrostatic compressive virgin stress $p_{0}=1$ was applied.

\subsubsection{Comparison of results}
For the above input data the analytical solution for the extent of the plastic zone was computed as 1.3, which means it is within the inclusion.
The theoretical solution for the radial displacement was 1.262 which compares well with the IGABEM solution of 1.269. Convergence to 1\% of residual was achieved after 6 iterations.

\section{Practical example}
\begin{figure}
\begin{center}
\begin{overpic}[scale=0.2]{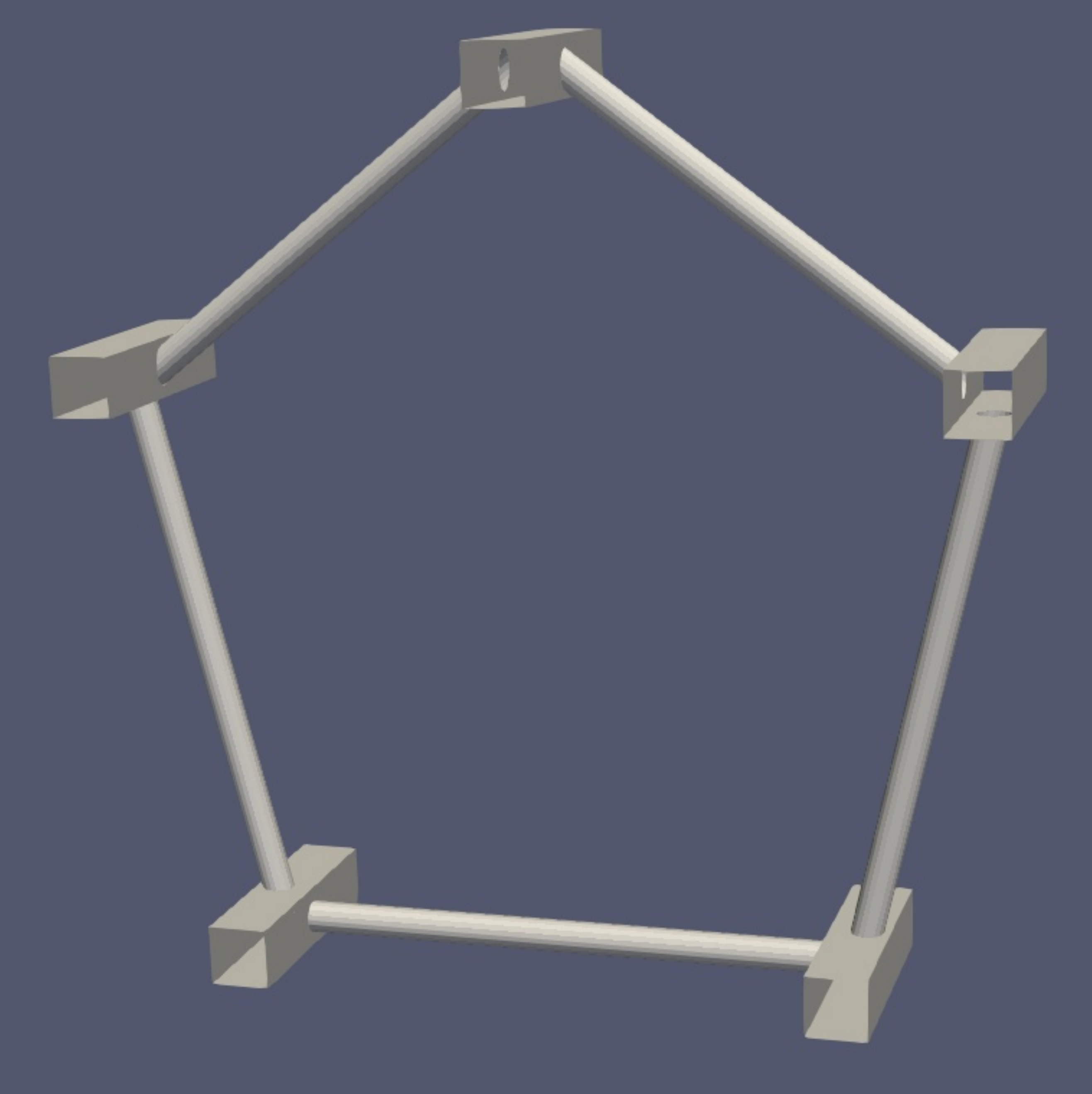}
\put(40,5){Stage 1}
\end{overpic}
\begin{overpic}[scale=0.208]{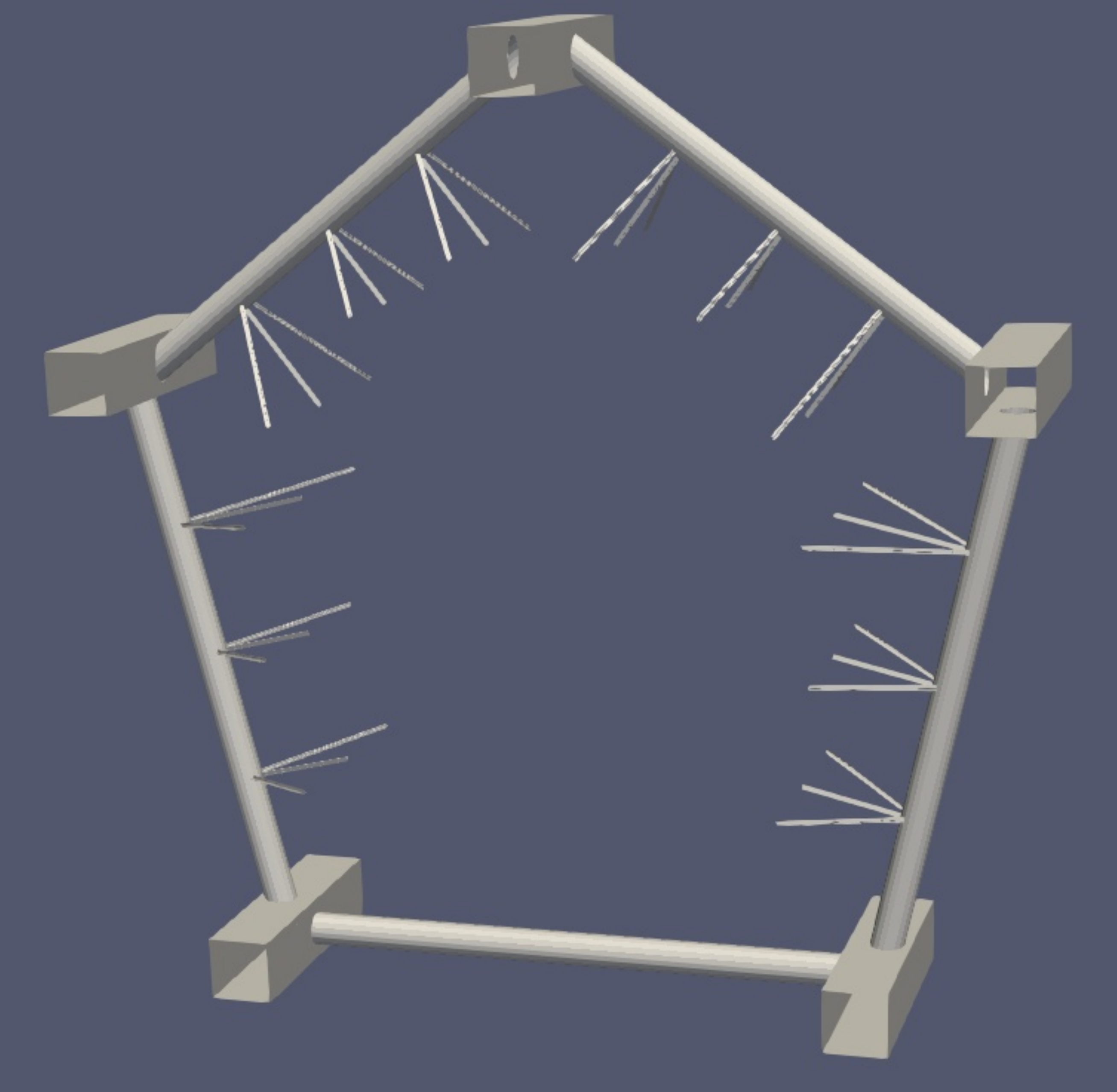}
\put(40,5){Stage 2}
\end{overpic}
\begin{overpic}[scale=0.2]{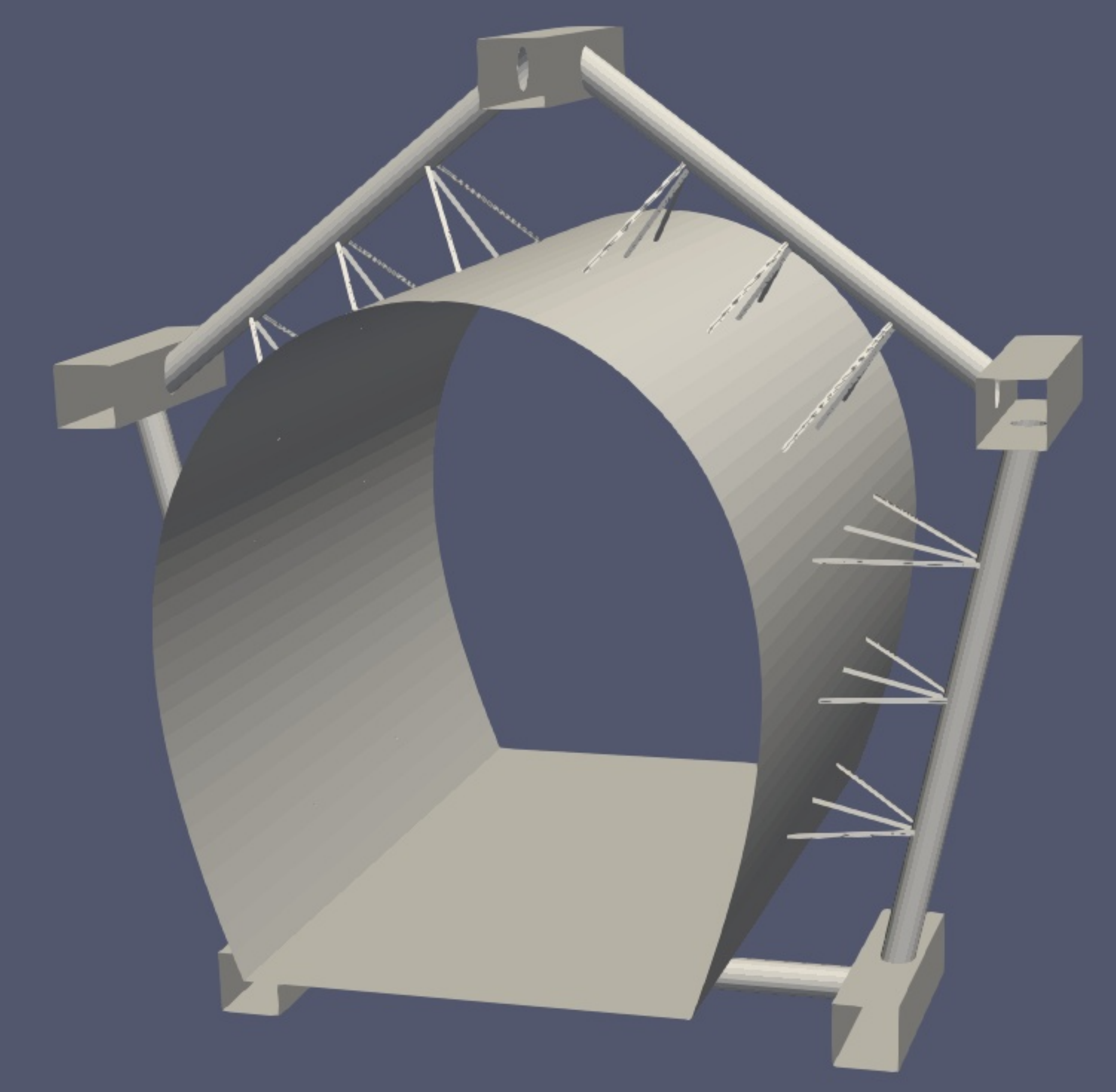}
\put(40,5){Stage 3}
\end{overpic}
\caption{Example of large cavern: Perspective view of excavation stages.}
\label{Presup}
\end{center}
\end{figure}
The practical example relates to the excavation of a large underground cavern with a height of 60 m, a width of 50 m and an extension of 300 m.
The cavern is at a depth of 500 m resulting in a vertical virgin compressive stress of 15 MPa with $k_{0}=0.8$. The rock mass properties are listed in Table \ref{tab:RM}. A Mohr-Coulomb yield condition with a dilation angle of $\psi=0$ is assumed.
\begin{mytable}
  {H}               
  {Rock mass properties}  
  {tab:RM}  
  {cccc}         
  \mytableheader{  E (MPa) & $\nu$ & $c$ (MPa) & $\phi$ (degrees) }  
 23000 & 0.25 & 1.63 & 36\\
\end{mytable}%
The aim of the simulation is to investigate if a pre-installation of ground support would lead to a safer and more economic excavation. The required excavation stages are shown on a 60 m section of the cavern in \myfigref{Presup}. In the first stage a tunnel system consisting of circular tunnels (by raise boring) and rectangular horizontal tunnels (by blasting) is excavated. Next cables are installed from the circular tunnels of such length that they reach to the surface of the cavern to be excavated. Finally the cavern is excavated. The cross-sectional diameter of the cables is 25 mm with the modulus of elasticity of E=210 GPa.
\begin{figure}
\begin{center}
\includegraphics[scale=0.5]{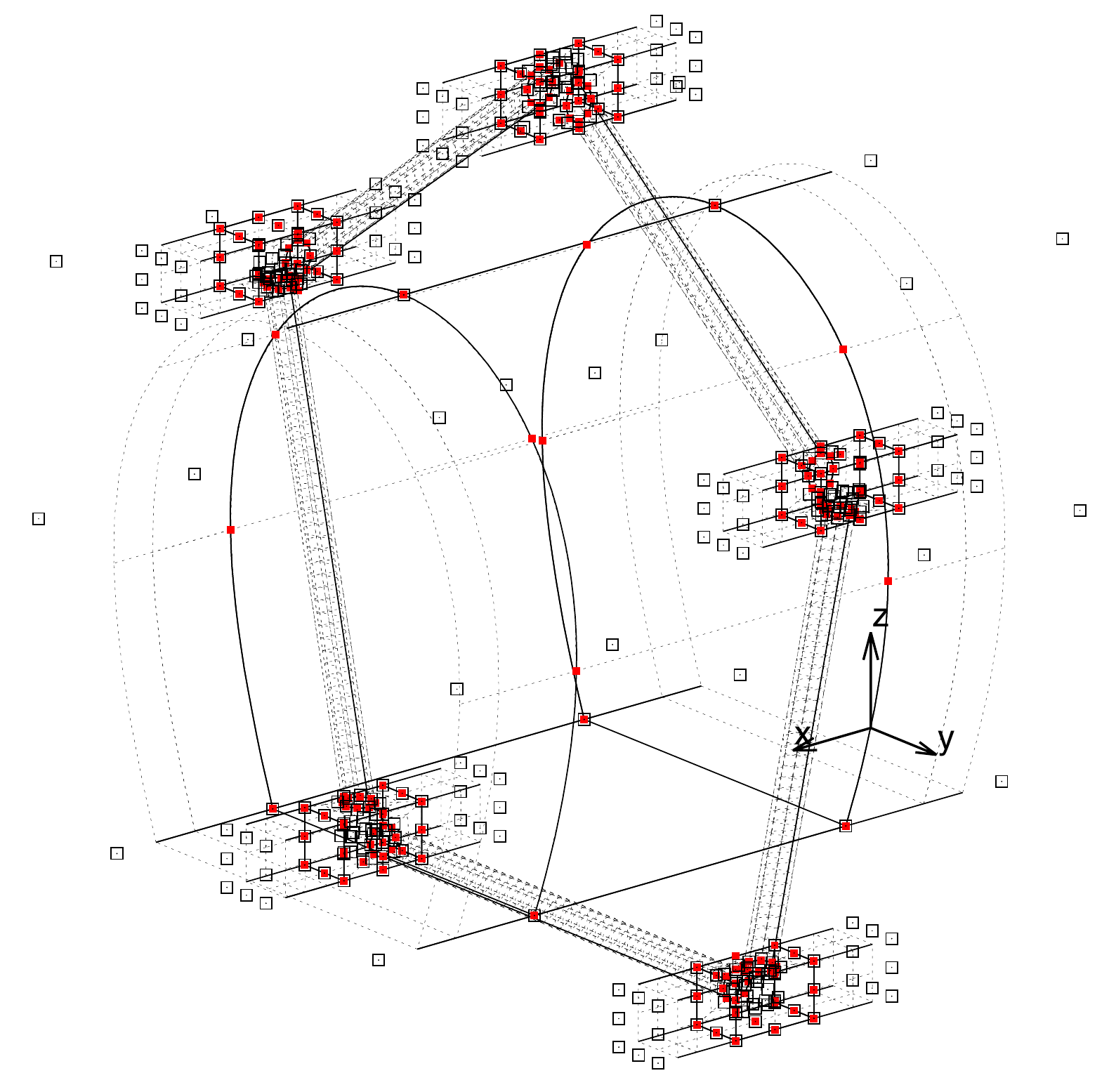}
\caption{IGABEM model of excavation surfaces, showing control points as hollow squares and collocation points as filled red squares.}
\label{Surf}
\end{center}
\end{figure}

\subsection{The simulation model}
For a preliminary analysis we analyse a 60m section of the cavern  At the edges of the section we use infinite plane strain boundary elements. 
Nine cables per circular excavation are considered.
We start with the definition of the excavation surfaces using the method outlined in section \ref{DefNURBS}. For the description of the smooth shape of the cavern walls 4 control points and a basis function of order 3 has been used for each half. 
For the variation of the unknowns we use the same basis functions as for the description of the geometry except that for the bottom surface the basis function was elevated by one order from linear to quadratic across the cavern.
\begin{figure}[H]
\begin{center}
\includegraphics[scale=0.5]{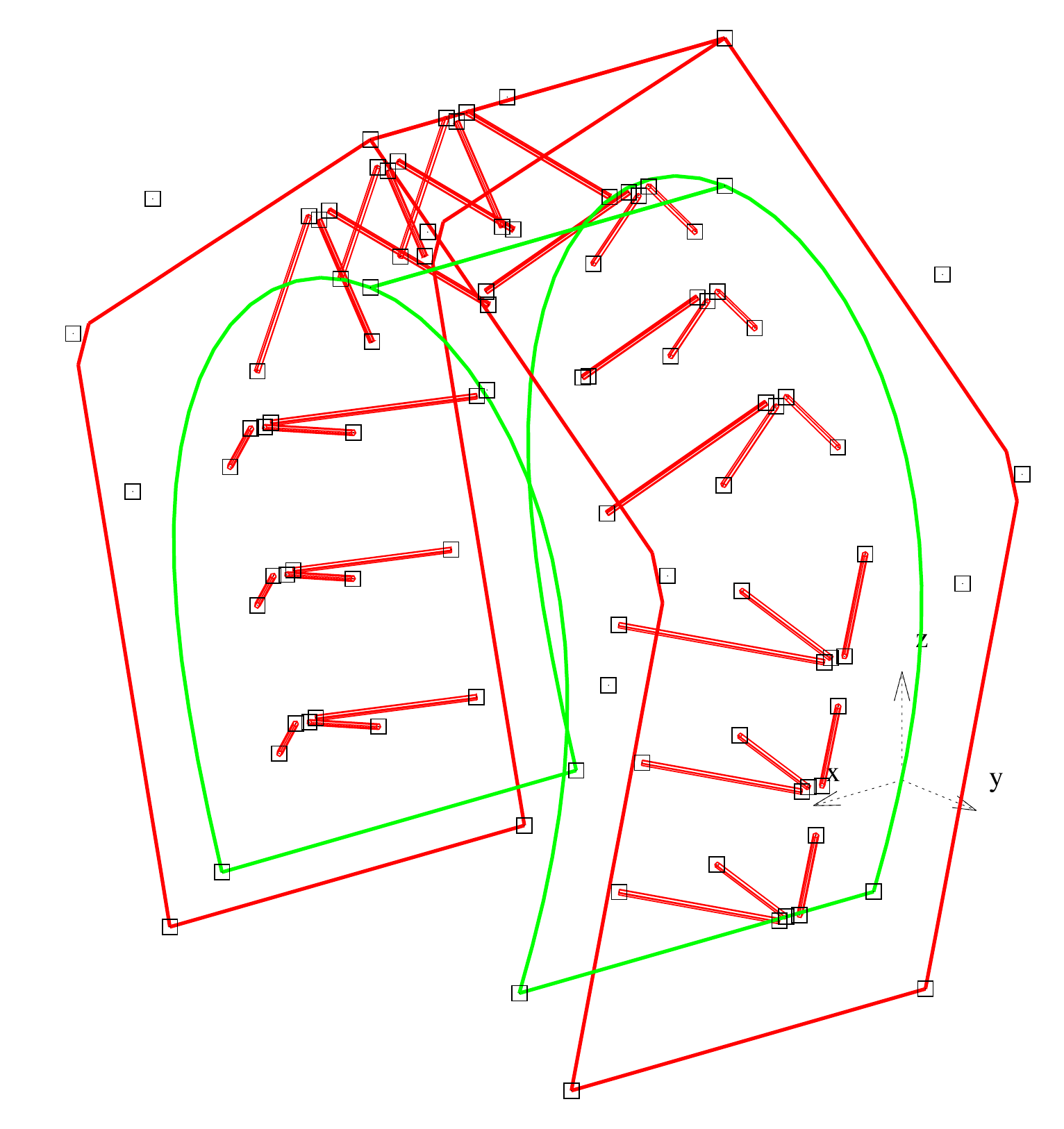}
\caption{Definition of the inclusions showing bolt inclusions and bounding surfaces (red and green lines) defining the part of the domain where non-linear behaviour is modelled.}
\label{Incl}
\end{center}
\end{figure}

The resulting simulation model is shown in \myfigref{Surf} for the final excavation stage and has 612 degrees of freedom.
\textbf{ It should be stressed that no mesh generation is involved here. The dotted lines in \myfigref{Surf} define integration regions, which are automatically determined depending of the location of collocation points and their aspect ratio.}

\begin{figure}
\begin{center}
\includegraphics[scale=0.25]{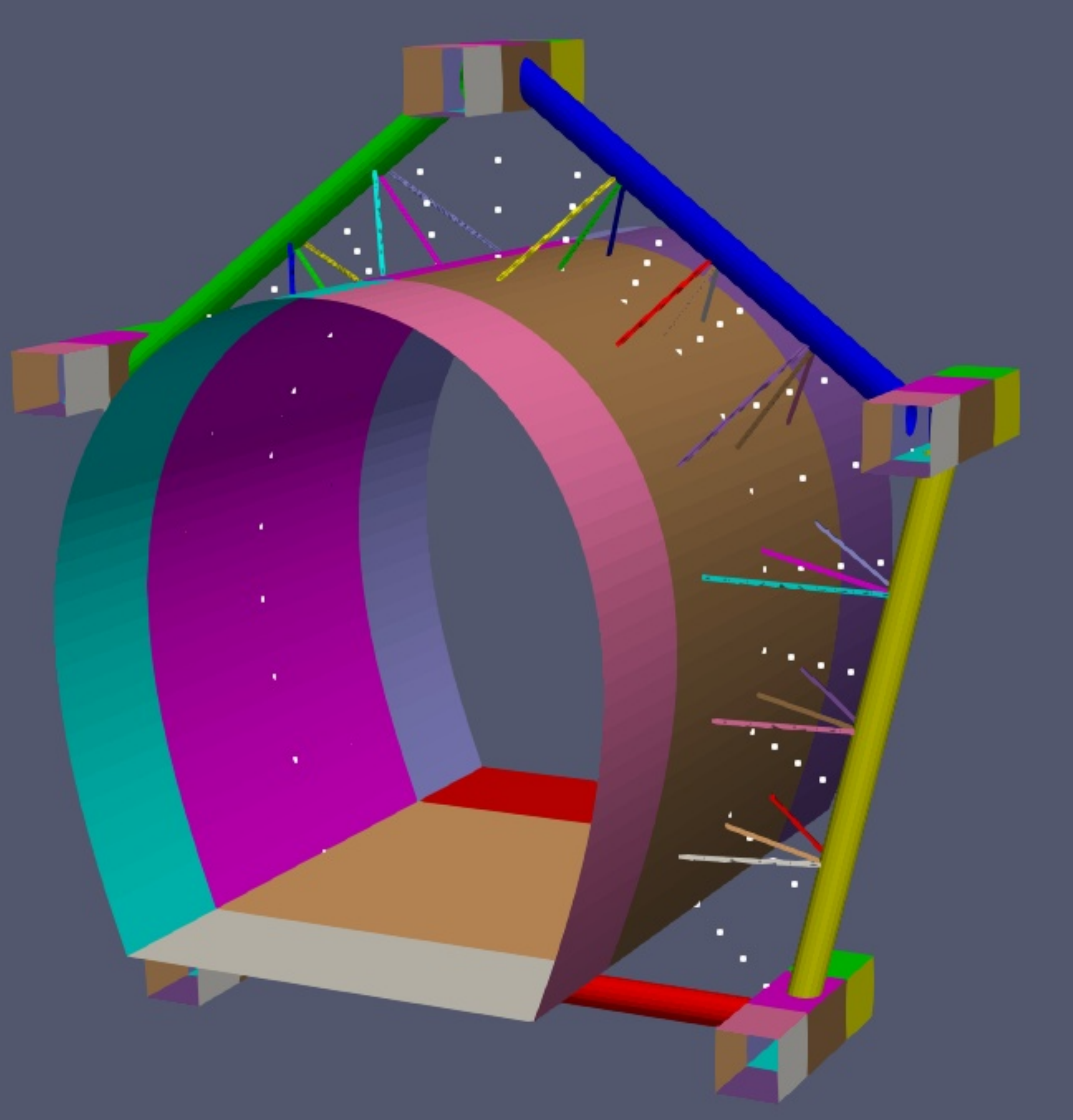}
\caption{Perspective view of geometry definition for the final excavation stage. Boundary Patches are colour coded. The grid points inside the general inclusion for modelling non-linear behaviour are shown as white dots.}
\label{Meshfinal}
\end{center}
\end{figure}

Next the inclusions are specified. This relates the cables and to the part of the domain where it is assumed that non-linear behaviour is taking place. In this study we concentrate on the rock mass behaviour between the tunnels and the excavation surface, so this part is selected.
The definition of the inclusions is shown in \myfigref{Incl} and the final model in \myfigref{Meshfinal}.

\subsection{Preliminary results}
\begin{figure}[H]
\begin{center}
\includegraphics[scale=0.25]{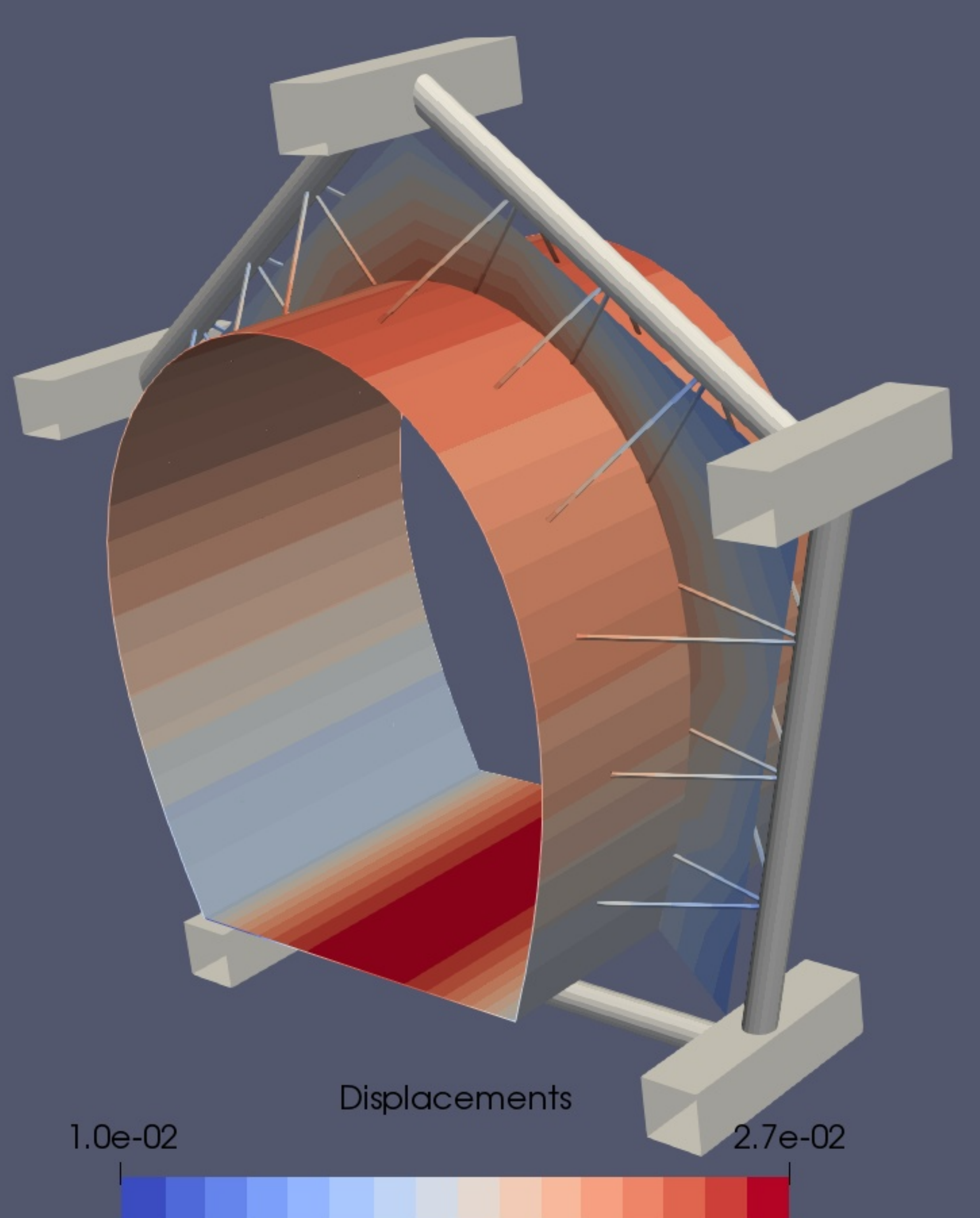}
\caption{Result of the analysis: Contours of absolute displacement plotted on the cavern surface, the cables and a plane inside the general inclusion.}
\label{RibinRockdisp}
\end{center}
\end{figure}

The novel approach to simulation will be used for a detailed study to determine if pre-installation of ground support can make the excavation of large underground caverns more feasible, safe and economic. Here only a preliminary result is presented in \myfigref{RibinRockdisp}. For this case the convergence to 1\% residual was achieved in 10 iterations.

 It is clear that more cables are required to provide adequate pre-support but it should be noted that the number of cables can be increased substantially without significantly increasing the simulation effort. This is because cables are generated automatically and because the numerical effort only increases slightly due to the fact that analytical integration is used. It should be noted that the number of degrees of freedom is not increased by increasing the number of cables.

\section{Summary and Conclusions}
A new approach to simulation for underground excavations using isogeometric methods and NURBS has been presented, which makes it extremely user friendly and efficient.
No mesh generation is required and geometric data can be specified in a CAD data format. Since CAD software also uses NURBS and trimmed NURBS surfaces the connection to CAD is natural.
Parts of the domain that have different material properties or behave in an inelastic way can be considered but a definition of inclusions is necessary. This means that for most practical applications only near field effects are considered in the simulation, i.e inclusions have limited extent and are located near the region of interest. 
If this restriction can be accepted, then the software can do everything a FEM simulation can do, but with a drastic reduction of user effort and execution time.

One of the main contributions of the paper is the presentation of analytically integrated rock bolt inclusions.
Numerical integration is the most compute intensive aspect of any BEM simulation. In case of large patches the number of Gauss points can be high. This is because of the nature of the integrals, which exhibit singular behaviour, i.e. where the value of the integrand  increases rapidly as the source point is approached. For the numerical volume integration, where the integration is in 3 directions the computational effort is particularly high.
However, it is noted that the matrices that require volume integration only depend on geometry and can be precomputed for a particular model.
Also, the computations are able to exploit to a very high degree the use of multiple processor hardware.

At the time of writing the simulation model has been written in MATLAB, had a first application to a real problem and will be used for a detailed study to simulate the effect of pre-support for large caverns.

Although concrete arches and shotcrete can be modelled with the general inclusion approach presented here, it becomes cumbersome, when the thickness is small. Therefore, special shotcrete inclusions with a semi-analytical integration are being developed next. It should be pointed out that it is also possible to develop inclusion geometries, that are not restricted to the cuboid shape presented here (in the sense that only two surfaces can be of general shape).

Finally, it is hoped that this contribution will make the BEM more attractive for simulations of underground excavations.

\section{Acknowledgements}
The authors would like to acknowledge the significant contribution of Scott Sloan to the topic of limit analysis. The first author had the pleasure of spending some time working at the centre of excellence at the University of Newcastle directed by him and was impressed by his leadership and quest for innovation.
Thanks are due to Vaughan Griffiths, for supplying a MATLAB code for the elasto-plastic constitutive matrix.
We are thankful to the Institute of Soil Mechanics, Foundation Engineering and Computational Geotechnics of the TU Graz for facilitating the access to the University version of PLAXIS.

\appendix
\section{Linear inclusion. Volume regular integral}
The analytical solution in Voigt notation (only terms different from zero are listed) is for $l=1$:
 \begin{eqnarray}
 \label{eq_DeltaEprimeReg_14M1}
 \triangle \fund{E^{\prime}}(1,4)(\tilde{y}^{\prime}\neq 0) & = & 2 C\pi R^2 \frac{C_3}{H \tilde{y}^{\prime}}\left[  r_{c0} + \frac{\tilde{z}^{\prime} \triangle z^{\prime}-\tilde{y}^{\prime 2}}{r_{c1}}  \right]   \nonumber \\
 \label{eq_DeltaEprimeReg_16M1}
\triangle \fund{E^{\prime}}(1,6) (\tilde{y}^{\prime}\neq 0) & = & 2 C\pi R^2 \frac{C_3}{H} \left[ \frac{H}{r_{c1}} + \ln \left( \frac{r_{c0}-\tilde{z}^{\prime}}{r_{c1}+\triangle z^{\prime}}\right)  \right]   \nonumber \\ 
\triangle \fund{E^{\prime}}(1,6)(\tilde{y}^{\prime}=0) & = &  2 C\pi R^2 \frac{C_3}{H} 
\left\{ \frac{\tilde{z}^{\prime}}{\lvert \triangle z^{\prime} \rvert} - 
- \frac{\tilde{z}^{\prime}}{\lvert \tilde{z}^{\prime} \rvert} -
\ln \left[ \left(\triangle z^{\prime}\right)^{\bar{\triangle z^{\prime}}}
\left( -\tilde{z}^{\prime} \right)^{\bar{z^{\prime}}} 
\right] \right\} \nonumber \\
 \label{eq_DeltaEprimeReg_21M1}
\triangle \fund{E^{\prime}}(2,1) (\tilde{y}^{\prime}\neq 0) & = & C\pi \frac{R^2}{H \tilde{y}^{\prime}} \left(  \frac{y^{\prime 2}-\tilde{z}^{\prime} \triangle z^{\prime}}{r_{c1}} - r_{c0}  \right)   \nonumber \\
 \label{eq_DeltaEprimeReg_22M1}
\triangle \fund{E^{\prime}}(2,2)(\tilde{y}^{\prime}\neq 0) & = & C\pi \frac{R^2}{H \tilde{y}^{\prime}} \, \left\{  \frac{\tilde{z}^{\prime 2} + 2 C_3 r_{c0}^2}{r_{c0}} + 
\right. \nonumber \\ 
& & + \left. \frac{\triangle z^{\prime} \left[ \tilde{z}^{\prime} \triangle z^{\prime 2}+(H+\tilde{z}^{\prime})\tilde{y}^{\prime 2} \right] -2 C_3  r_{c1}^2 (\tilde{y}^{\prime 2}-\tilde{z}^{\prime} \triangle z^{\prime})}{r_{c1}^3} \right\}  \nonumber \\
 \label{eq_DeltaEprimeReg_23M1}
\triangle \fund{E^{\prime}}(2,3) (\tilde{y}^{\prime}\neq 0) & = &  C\pi R^2 \frac{\tilde{y}^{\prime}}{H} \, \left( \frac{1}{r_{c0}} - \frac{r_{c0}^2 + 2 H^2 - 3H \tilde{z}^{\prime}}{r_{c1}^3}  \right)   \nonumber \\
 \label{eq_DeltaEprimeReg_25M1}
\triangle \fund{E^{\prime}}(2,5) (\tilde{y}^{\prime}\neq 0) & = & 2 C\pi \frac{R^2}{H}  \, \left\{ 
\frac{H\left[(C_3-1)H(\triangle z^{\prime}-\tilde{z}^{\prime})+C_3 r_{c0}^2 - \tilde{z}^{\prime2}\right]}{r_{c1}^3} + \right. \nonumber \\
&  & \left. \tilde{z}^{\prime}(\frac{1}{r_{c1}}-\frac{1}{r_{c0}}) + C_3 \ln \frac{r_{c0}-\tilde{z}^{\prime}}{r_{c1}+\triangle z^{\prime}}  \right\}   \nonumber \\
\triangle \fund{E^{\prime}}(2,5)(\tilde{y}^{\prime}=0) & = &  \triangle \fund{E^{\prime}}_{M_1}(1,6)(\tilde{y}^{\prime}=0)   \\
 \label{eq_DeltaEprimeReg_31M1}
\triangle \fund{E^{\prime}}(3,1) (\tilde{y}^{\prime}\neq 0) & = & C\pi \frac{R^2}{H}
 \left( \ln \frac{r_{c1}+\triangle z^{\prime}}{r_{c0}-\tilde{z}^{\prime}} - \frac{H}{r_{c1}}  \right)   \nonumber \\
\triangle \fund{E^{\prime}}(3,1)(\tilde{y}^{\prime}=0) & = & - \frac{\triangle \fund{E^{\prime}}_{M_1}(1,6)(\tilde{y}^{\prime}=0)}{2 C_3}  \nonumber \\ 
 \label{eq_DeltaEprimeReg_32M1}
\triangle \fund{E^{\prime}}(3,2)(\tilde{y}^{\prime}\neq 0) & = & C\pi \frac{R^2}{H}
 \left[ - \frac{\tilde{z}^{\prime}}{r_{c0}^2}- \frac{\triangle z^{\prime}(2H \triangle z^{\prime}-H\tilde{z}^{\prime}+r_{c0}^2)}{r_{c1}^3} + \ln \frac{\triangle z^{\prime}+r_{c1}}{r_{c0}-\tilde{z}^{\prime}}  \right]   \nonumber \\
 \triangle \fund{E^{\prime}}(3,2)(\tilde{y}^{\prime}=0) & = & \triangle \fund{E^{\prime}}_{M_1}(3,1)(\tilde{y}^{\prime}=0)   \nonumber  \\ 
 \label{eq_DeltaEprimeReg_33M1}
\triangle \fund{E^{\prime}}(3,3) (\tilde{y}^{\prime}\neq 0) & = & C\pi \frac{R^2}{H}
 \left\{ \frac{H}{r_{c1}^3} \left[ (3+2C_3)H^2+2(1+C_3)\tilde{y}^{\prime 2} - 2(3+2C_3) H \tilde{z}^{\prime} + \right. \right. \nonumber \\
 & & \left. \left. +(3+2C_3)\tilde{z}^{\prime 2}\right]+\tilde{z}^{\prime} (\frac{1}{r_{c0}} - \frac{1}{r_{c1}}) - 2(1+C_3) \ln \frac{r_{c1}+\triangle z^{\prime}}{r_{c0}-\tilde{z}^{\prime}} \right\}  \nonumber \\
  \triangle \fund{E^{\prime}}(3,3)(\tilde{y}^{\prime}=0) & = &  \frac{1+C_3}{C_3} \triangle\fund{E^{\prime}}_{M_1}(1,6)(\tilde{y}^{\prime}=0)  \nonumber \\ 
 \label{eq_DeltaEprimeReg_35M1}
\triangle \fund{E^{\prime}}(3,5) (\tilde{y}^{\prime}\neq 0) & = & 2 C\pi \frac{R^2}{H \tilde{y}^{\prime}} 
\left[ \frac{(2+C_3)\tilde{y}^{\prime 2} + (1+C_3) \tilde{z}^{\prime 2}}{r_{c0}} + 
\frac{(1+C_3) \tilde{z}^{\prime} \triangle z^{\prime 3}}{r_{c1}^3} -  \right. \nonumber \\
& & \left. - \frac{(2+C_3)\tilde{y}^{\prime 4}+\tilde{y}^{\prime 2} \triangle z^{\prime} ((3+C_3)\triangle z^{\prime} - C_3 \tilde{z}^{\prime})}{r_{c1}^3}
 \right]  \nonumber
\end{eqnarray}

and for $l=2$:
\begin{eqnarray}
 \label{eq_DeltaEprimeReg_14M2}
 \triangle \fund{E^{\prime}}(1,4)(\tilde{y}^{\prime}\neq 0) & = & 2 C\pi R^2 \frac{C_3}{H \tilde{y}^{\prime} }\left[  r_{c1} + \frac{\tilde{z}^{\prime} \triangle z^{\prime}-\tilde{y}^{\prime 2}}{r_{c0}} \right]    \nonumber \\
  \label{eq_DeltaEprimeReg_16M2}
\triangle \fund{E^{\prime}}(1,6) (\tilde{y}^{\prime}\neq 0) & = & - 2 C\pi R^2 \frac{C_3}{H} \left[ \frac{H}{r_{c0}} + \ln \left( \frac{r_{c0}-\tilde{z}^{\prime}}{r_{c1}+\triangle z^{\prime}}\right)  \right]   \nonumber \\
\triangle \fund{E^{\prime}}(1,6)(\tilde{y}^{\prime}=0) & = &  2 C\pi R^2 \frac{C_3}{H} 
\left\{ \frac{\triangle z^{\prime}}{\lvert \triangle z^{\prime} \rvert} 
- \frac{\triangle z^{\prime}}{\lvert \tilde{z}^{\prime} \rvert} +
\ln \left[ \left(\triangle z^{\prime}\right)^{\bar{\triangle z^{\prime}}}
\left( -\tilde{z}^{\prime} \right)^{\bar{z^{\prime}}} 
\right] \right\}  \nonumber \\
 \label{eq_DeltaEprimeReg_21M2}
\triangle \fund{E^{\prime}}(2,1) (\tilde{y}^{\prime}\neq 0) & = & C\pi  \frac{R^2}{H \tilde{y}^{\prime}} \left(  \frac{y^{\prime 2}-\tilde{z}^{\prime} \triangle z^{\prime}}{r_{c0}} - r_{c1}  \right)   \nonumber \\
 \label{eq_DeltaEprimeReg_22M2}
\triangle \fund{E^{\prime}}(2,2)(\tilde{y}^{\prime}\neq 0) & = & C\pi \frac{R^2}{H \tilde{y}^{\prime}} \, \left[  \frac{(1+2 C_3)H^2+\tilde{z}^{\prime 2} - 2 H \tilde{z}^{\prime}(1+2C_3)+2C_3r_{c0}^2}{r_{c1}} + 
\right. \nonumber \\ 
& & + \left. \frac{\tilde{z}^{\prime} (r_{c0}^2 \triangle z^{\prime} + H \tilde{y}^{\prime 2}) - 2 C_3  r_{c0}^2 (\tilde{y}^{\prime 2}-\tilde{z}^{\prime} \triangle z^{\prime})}{r_{c0}^3} \right] 
 \nonumber \\
 \label{eq_DeltaEprimeReg_23M2}
\triangle \fund{E^{\prime}}(2,3) (\tilde{y}^{\prime}\neq 0) & = &  C\pi R^2 \frac{\tilde{y}^{\prime}}{H} \, \left( \frac{1}{r_{c1}} - \frac{r_{c0}^2 + H \tilde{z}^{\prime}}{r_{c0}^3}  \right)   \nonumber \\
 \label{eq_DeltaEprimeReg_25M2}
\triangle \fund{E^{\prime}}(2,5) (\tilde{y}^{\prime}\neq 0) & = & 2 C\pi \frac{R^2}{H}  \, \left( 
\frac{\tilde{z}^{\prime} - C_3H}{r_{c0}} - \frac{H \tilde{y}^{\prime 2}}{r_{c0}^3} + \frac{\triangle z^{\prime}}{r_{c1}} + C_3 \ln \frac{r_{c1}+\triangle z^{\prime}}{r_{c0}- \tilde{z}^{\prime}}  \right)   \nonumber \\
\triangle \fund{E^{\prime}}(2,5)(\tilde{y}^{\prime}=0) & = &  \triangle \fund{E^{\prime}}_{M_2}(1,6)(\tilde{y}^{\prime}=0)  \\ \nonumber 
 \label{eq_DeltaEprimeReg_31M2}
\triangle \fund{E^{\prime}}(3,1) (\tilde{y}^{\prime}\neq 0) & = & C\pi \frac{R^2}{H}
 \left( - \ln \frac{r_{c1}+\triangle z^{\prime}}{r_{c0}-\tilde{z}^{\prime}} + \frac{H}{r_{c0}}  \right)   \nonumber \\
\triangle \fund{E^{\prime}}(3,1)(\tilde{y}^{\prime}=0) & = & - \frac{\triangle \fund{E^{\prime}}_{M_2}(1,6)(\tilde{y}^{\prime}=0)}{2 C_3}  \nonumber \\ 
 \label{eq_DeltaEprimeReg_32M2}
\triangle \fund{E^{\prime}}(3,2)(\tilde{y}^{\prime}\neq 0) & = & C\pi \frac{R^2}{H}
 \left[ \frac{\triangle z^{\prime}}{r_{c1}}+ \frac{\tilde{z}^{\prime}(H \tilde{z}^{\prime}+r_{c0}^2)}{r_{c0}^3} - \ln \frac{\triangle z^{\prime}+r_{c1}}{r_{c0}-\tilde{z}^{\prime}}  \right]   
 \nonumber \\
 \triangle \fund{E^{\prime}}(3,2)(\tilde{y}^{\prime}=0) & = & \triangle \fund{E^{\prime}}_{M_2}(3,1)(\tilde{y}^{\prime}=0)  \nonumber \\ 
 \label{eq_DeltaEprimeReg_33M2}
\triangle \fund{E^{\prime}}(3,3) (\tilde{y}^{\prime}\neq 0) & = & C\pi \frac{R^2}{H}
 \left[ \frac{H\tilde{y}^{\prime 2}}{r_{c0}^3} - \frac{2(1+C_3)H+\tilde{z}^{\prime}}{r_{c0}} - \frac{\triangle z^{\prime}}{r_{c1}} + \right. \nonumber \\
 & & \left. + 2(1+C_3) \ln \frac{r_{c1}+\triangle z^{\prime}}{r_{c0}-\tilde{z}^{\prime}} \right] 
  \nonumber \\
  \triangle \fund{E^{\prime}}(3,3)(\tilde{y}^{\prime}=0) & = &  \frac{1+C_3}{C_3} \triangle\fund{E^{\prime}}_{M_2}(1,6)(\tilde{y}^{\prime}=0)  \nonumber \\ 
 \label{eq_DeltaEprimeReg_35M2}
\triangle \fund{E^{\prime}}(3,5) (\tilde{y}^{\prime}\neq 0) & = & 2 C\pi \frac{R^2}{H \tilde{y}^{\prime}} 
\left[ \frac{(1+C_3)r_{c0}^2 + \tilde{y}^{\prime 2} + (1+C_3) (H^2-2H\tilde{z}^{\prime})}{r_{c1}} + 
\frac{(1+C_3) \tilde{z}^{\prime 3} \triangle z^{\prime}}{r_{c0}^3} -  \right. \nonumber \\
& & \left. - \frac{(2+C_3)\tilde{y}^{\prime 4}-\tilde{y}^{\prime 2} \tilde{z}^{\prime} (C_3 \triangle z^{\prime} - (3+C_3) \tilde{z}^{\prime})}{r_{c0}^3}
 \right]  \nonumber 
\end{eqnarray}
where:
\begin{equation}
\triangle z^{\prime} = H - \tilde{z}^{\prime} \quad r_{c1} = \sqrt{\tilde{y}^{\prime 2}+\triangle z^{\prime 2}} \quad  r_{c0} = \sqrt{\tilde{y}^{\prime 2}+z^{\prime 2}} \quad \bar{\triangle z^{\prime}}=\triangle z^{\prime}/\lvert \triangle z^{\prime} \rvert \quad \bar{z^{\prime}}=\tilde{z}^{\prime}/\lvert \tilde{z}^{\prime} \rvert (\triangle z^{\prime}
\end{equation}

\section{Linear inclusion. Volume singular integral}
The terms  of $\triangle \fund{E}^{\prime} = \triangle \fund{E}^{\prime}_1+\triangle \fund{E}^{\prime}_2$, in Voigt notation, different from zero are
 for $l=1$:

 \begin{eqnarray}
 \label{eq_DeltaEprime_1625M1}
 \triangle \fund{E^{\prime}}(1,6) = \triangle \fund{E^{\prime}}(2,5) & = & \frac{C\pi}{4H} \left\{ H^2 \left[ 8+8 C_3-(9+8 C_3) \cos{\tilde{\theta}} +  \cos{3\tilde{\theta}} \right. \right. + \nonumber \\
& & \left. \left. R^2 \left( (3-8C_3) \cos{\tilde{\theta}} + \cos{3\tilde{\theta}} + 8C_3 \ln\frac{\cos{\tilde{\theta}/2}}{\sin{\tilde{\theta}/2}} \right) \right] \right\} \nonumber  \\
 \label{eq_DeltaEprime_3132M1}
 \triangle \fund{E^{\prime}}(3,1) = \triangle \fund{E^{\prime}}(3,2) & = & \frac{C\pi}{8H} \left[ R^2 \left(11 \cos{\tilde{\theta}} + \cos{3\tilde{\theta}} + 8 \ln \frac{\sin{\tilde{\theta}/2}}{\cos{\tilde{\theta}/2}}\right) - \right. \\
 & & \left. -4H^2\cos{\tilde{\theta}}\sin^2{\tilde{\theta}} \right]  \nonumber  \\
 \label{eq_DeltaEprime_33M1}
 \triangle \fund{E^{\prime}}(3,3) & = & -C\pi \frac{R^2}{4H} 
 \left[ (11+8C_3)\cos{\tilde{\theta}} + \cos{3\tilde{\theta}} + 8(1+C_3) \ln \frac{\sin{\tilde{\theta}/2}}{\cos{\tilde{\theta}/2}} \right] + \nonumber \\
 & & C\pi H \left( 1+4C_3+2\cos{\tilde{\theta}} + \cos{2\tilde{\theta}} \right) \sin^2{\tilde{\theta}/2} \nonumber
\end{eqnarray}

and for $l=2$

\begin{eqnarray}
 \label{eq_DeltaEprime_1625M2}
 \triangle \fund{E^{\prime}}(1,6) = \triangle \fund{E^{\prime}}(2,5) & = & \frac{C\pi}{4H} \left\{ H^2 \left[ 8+8 C_3-(9+8 C_3) \cos{\tilde{\theta}} +  \cos{3\tilde{\theta}} \right] \right. + \nonumber \\
& & \left. + R \left[ (8C_3-3) R\cos{\tilde{\theta}} - R\cos{3\tilde{\theta}} + 8(H + 2C_3H + \right. \right. \nonumber \\ 
& & \left. \left. C_3 R \ln\frac{\sin{\tilde{\theta}/2}}{\cos{\tilde{\theta}/2}})
 + 4H (-1-4C_3+\cos{\tilde{2\theta}}) \sin{\tilde{\theta}} 
\right] \right\} \nonumber  \\
 \label{eq_DeltaEprime_3132M2}
 \triangle \fund{E^{\prime}}(3,1) = \triangle \fund{E^{\prime}}(3,2) & = & 
 -\frac{C\pi H}{2} \cos{\tilde{\theta}} \sin^2{\tilde{\theta}} - 
 \frac{C\pi R}{8H} \left[ R \left( 11 \cos{\tilde{\theta}} + \cos{3\tilde{\theta}} + \right. \right.
\nonumber   \\
& &  \left. \left. 8 \ln \frac{\sin{\tilde{\theta}/2}}{\cos{\tilde{\theta}/2}} \right)
    -2H \left( 5 \sin{\tilde{\theta}} + \sin{3\tilde{\theta}} - 4 \right) \right]
 	  \\
  \label{eq_DeltaEprime_33M2}
 \triangle \fund{E^{\prime}}(3,3) & = & C\pi H \left( 1+4C_3+2\cos{\tilde{\theta}} + \cos{2\tilde{\theta}} \right) \sin^2{\tilde{\theta}/2} + \nonumber \\
 &  &  C\pi \frac{R}{4H} \left[ (11+8C_3)R\cos{\tilde{\theta}} + R\cos{3\tilde{\theta}} + \right. \nonumber \\
& &  \left. 8  \left(H + 2C_3 H + (1+C_3)R \ln \frac{\sin{\tilde{\theta}/2}}{\cos{\tilde{\theta}/2}} \right) - \right. \nonumber \\
 & & \left.  4H \left(3 + 4C_3 + \cos{2\tilde{\theta}} \right) \sin{\tilde{\theta}} \right]
 \nonumber  
\end{eqnarray}

\bibliographystyle{myplainnat}
\bibliography{bookbib}

\begin{thebibliography}{25}
\providecommand{\natexlab}[1]{#1}
\providecommand{\url}[1]{\texttt{#1}}
\expandafter\ifx\csname urlstyle\endcsname\relax
  \providecommand{\doi}[1]{doi: #1}\else
  \providecommand{\doi}{doi: \begingroup \urlstyle{rm}\Url}\fi

\bibitem[Aliabadi(2002)]{aliabadi}
Aliabadi,~M.
\newblock \emph{The Boundary Element Method, Volume 2: Applications in Solids
  and Structures}.
\newblock Wiley, 2002.

\bibitem[An et~al.(2018)An, Yu, Bui, Wang, and Trinh]{An2018}
An,~Z.; Yu,~T.; Bui,~T.; Wang,~C.; Trinh,~N., Implementation of isogeometric
  boundary element method for 2-d steady heat transfer analysis, \emph{Advances
  in Engineering Software}, 116:\penalty0 36--49, 2018.

\bibitem[Atroshchenko et~al.(2018)Atroshchenko, Tomar, Xu, and
  Bordas]{doi:10.1002/nme.5778}
Atroshchenko,~E.; Tomar,~S.; Xu,~G.; Bordas,~S.P., Weakening the tight coupling
  between geometry and simulation in isogeometric analysis: From sub- and
  super-geometric analysis to geometry-independent field approximation (gift),
  \emph{International Journal for Numerical Methods in Engineering},
  114\penalty0 (10):\penalty0 1131--1159, 2018.

\bibitem[Banerjee and Butterfield(1981)]{Banerjee1981}
Banerjee,~P.K.; Butterfield,~R.
\newblock \emph{Boundary element methods in engineering science}.
\newblock McGraw-Hill, 1981.

\bibitem[Banerjee(1994)]{Banerjee}
Banerjee,~P.
\newblock \emph{The Boundar Element Method in Engineering}.
\newblock McGraw-Hill, 1994.

\bibitem[Banerjee and Raveendra(1986)]{BanerjeeRaveendra1986}
Banerjee,~P.; Raveendra,~S., Advanced boundary element of two- and
  three-dimensional problems of elastoplasticity, \emph{International Journal
  for Numerical Methods in Engineering}, 23\penalty0 (6):\penalty0 985--1002,
  1986.

\bibitem[Beer(2015)]{Beer2015b}
Beer,~G., Mapped infinite patches for the {NURBS} based boundary element
  analysis in geomechanics, \emph{Computers and Geotechnics}, 66:\penalty0
  66--74, 2015.

\bibitem[Beer et~al.(2019)Beer, Marussig, and Duenser]{BeerMarussig}
Beer,~G.; Marussig,~B.; Duenser,~C.
\newblock \emph{The isogeometric Boundary Element method}, volume~90 of
  \emph{Lecture Notes in Applied and Computational Mechanics}.
\newblock Springer Nature, 2019.

\bibitem[Brebbia and Walker(1980)]{Brebbia1980}
Brebbia,~C.A.; Walker,~S.
\newblock \emph{Boundary element techniques in engineering}.
\newblock Newnes-Butterworths, 1980.

\bibitem[Brebbia et~al.(1984)Brebbia, Telles, and Wrobel]{Brebbia1983}
Brebbia,~C.A.; Telles,~J.C.; Wrobel,~L.
\newblock \emph{Boundary element techniques}.
\newblock Springer, 1984.

\bibitem[Fang et~al.(2020)Fang, An, Yu, and Bui]{Fang2020}
Fang,~W.; An,~Z.; Yu,~T.; Bui,~T., Isogeometric boundary element analysis for
  two-dimensional thermoelasticity with variable temperature, \emph{Engineering
  Analysis with Boundary Elements}, 110:\penalty0 80--94, 2020.

\bibitem[Gao and Davies(2011)]{Gao2011}
Gao,~X.; Davies,~T.
\newblock \emph{Boundary Element Programming in Mechanics}.
\newblock Cambridge University Press, Cambridge, UK, 2011.

\bibitem[Hoek and Brown(1980)]{Hoek}
Hoek,~E.; Brown,~T.
\newblock \emph{Underground excavations in rock}.
\newblock CRC Press, 1980.

\bibitem[Hughes et~al.(2005)Hughes, Cottrell, and Bazilevs]{Hughes2005a}
Hughes,~T.J.R.; Cottrell,~J.A.; Bazilevs,~Y., Isogeometric analysis: {CAD},
  finite elements, {NURBS}, exact geometry and mesh refinement, \emph{Computer
  Methods in Applied Mechanics and Engineering}, 194\penalty0
  (39--41):\penalty0 4135--4195, October 2005.

\bibitem[Kirsch(1898)]{Kirsch}
Kirsch,~, Die theorie der elastizit{\"a}t und die bed{\"u}rfnisse der
  festigkeitslehre., \emph{Zeitschrift des Vereines deutscher Ingenieure}, 42,
  1898.

\bibitem[Marussig et~al.(2015)Marussig, Zechner, Beer, and Fries]{Marussig2015}
Marussig,~B.; Zechner,~J.; Beer,~G.; Fries,~T.-P., Fast isogeometric boundary
  element method based on independent field approximation, \emph{Computer
  Methods in Applied Mechanics and Engineering}, 284:\penalty0 458--488, 2015.

\bibitem[Marussig(2016)]{marussig2016b}
Marussig,~B.
\newblock \emph{Seamless Integration of Design and Analysis through Boundary
  Integral Equations}.
\newblock Monographic Series TU Graz: Structural Analysis Verlag der
  Technischen Universit{\"a}t Graz, 2016.

\bibitem[Pan and Chou(1976)]{PanChow1976}
Pan,~Y.-C.; Chou,~T.-W., Point force solution for an infinite transversely
  isotropic solid, \emph{Journal of Applied Mechanics}, 43:\penalty0 608--612,
  1976.

\bibitem[Scott et~al.(2013)Scott, Simpson, Evans, Lipton, Bordas, Hughes, and
  Sederberg]{scott2013isogeometric}
Scott,~M.A.; Simpson,~R.N.; Evans,~J.A.; Lipton,~S.; Bordas,~S.P.A.;
  Hughes,~T.J.R.; Sederberg,~T.W., Isogeometric boundary element analysis using
  unstructured {T}-splines, \emph{Computer Methods in Applied Mechanics and
  Engineering}, 254:\penalty0 197--221, 2013.

\bibitem[Simpson et~al.(2012)Simpson, Bordas, Trevelyan, and
  Rabczuk]{simpson2012two}
Simpson,~R.N.; Bordas,~S.P.; Trevelyan,~J.; Rabczuk,~T., A two-dimensional
  isogeometric boundary element method for elastostatic analysis,
  \emph{Computer Methods in Applied Mechanics and Engineering}, 209:\penalty0
  87--100, 2012.

\bibitem[Simpson et~al.(2013)Simpson, Bordas, Lian, and
  Trevelyan]{simpson2013isogeometric}
Simpson,~R.N.; Bordas,~S.P.; Lian,~H.; Trevelyan,~J., An isogeometric boundary
  element method for elastostatic analysis: {2D} implementation aspects,
  \emph{Computers \& Structures}, 118:\penalty0 2--12, 2013.

\bibitem[Smith et~al.(2013)Smith, Griffiths, and Margetts]{Smith}
Smith,~I.M.; Griffiths,~D.V.; Margetts,~L.
\newblock \emph{Programming the Finite Element Method}.
\newblock Wiley, 2013.

\bibitem[Sun et~al.(2020)Sun, Gong, and Dong]{Sun2020}
Sun,~F.; Gong,~Y.; Dong,~C., A novel fast direct solver for {3D} elastic
  inclusion problems with the isogeometric boundary element method,
  \emph{Journal of Computational and Applied Mathematics}, 377:\penalty0
  112904, 2020.

\bibitem[Tanaka et~al.(2001)Tanaka, Matsumoto, and Takakuwa]{Tanaka1}
Tanaka,~M.; Matsumoto,~T.; Takakuwa,~S., Dual reciprocity {BEM} for
  time-stepping approach to the transient heat conduction problem in nonlinear
  materials, \emph{Computer Methods in Applied Mechanics and Engineering},
  195(37-40):\penalty0 4953--4961, 2001.

\bibitem[Wendland(1990)]{Wendlandbook}
Wendland,~W.L.
\newblock On the coupling of finite elements and boundary elements.
\newblock In: Kuhn,~G.; Mang,~H., editors, \emph{Discretization Methods in
  Structural Mechanics}, pages 405--414, Berlin, Heidelberg, 1990 Springer
  Berlin Heidelberg.

\end{thebibliography}

\end{document}